\documentclass[12pt]{amsart}
\usepackage[a4paper]{anysize}
\usepackage{amsthm}
\usepackage{amssymb}
\usepackage{times}
\usepackage{mathrsfs}

\include{xy}
\xyoption{v2}
\xyoption{all}
\xyoption{2cell}
\UseAllTwocells

    \newtheorem{theorem}{Theorem}[subsection]
    \newtheorem{proposition}[theorem]{Proposition}
    \newtheorem{lemma}[theorem]{Lemma}
    \newtheorem{corollary}[theorem]{Corollary}
    
    \theoremstyle{definition}
    \newtheorem{example}[theorem]{Example}
    \newtheorem{definition}[theorem]{Definition}
    \newtheorem{remark}[theorem]{Remark}
    \theoremstyle{remark}
    \newtheorem{claim}[theorem]{Claim}

    \def\set{\setcounter{equation}
             {\value{theorem}}\addtocounter{theorem}{1}}
    \numberwithin{equation}{subsection}
    \def\sset{\setcounter{subsubsection}
             {\value{theorem}}\addtocounter{theorem}{1}}

\def\A{{\mathbb A}}

\def\C{{\mathbb C}}

\def\cC{{{\mathscr C}}}

\def\cD{{{\mathscr D}}}

\def\D{{\mathbb D}}

\def\sD{\mathsf{D}}
\def\E{{\mathbb E}}
\def\cE{{{\mathscr E}}}

\def\F{{\mathbb F}}
\def\cF{{{\mathscr F}}}

\def\sF{\mathsf{F}}
\def\G{{\mathbb G}}
\def\cG{{{\mathscr G}}}

\def\cH{{{\mathscr H}}}

\def\cK{{{\mathscr K}}}
\def\cL{{{\mathscr L}}}

\def\cM{{{\mathscr M}}}

\def\N{{\mathbb N}}

\def\cO{{{\mathscr O}}}

\def\P{{\mathbb P}}

\def\Q{{\mathbb Q}}
\def\R{{\mathbb R}}

\def\Z{{\mathbb Z}}

\def\cV{{{\mathscr V}}}
\def\cK{{{\mathscr K}}}

\newcommand{\colim}[1]
{\underset{#1}{\mathrm{colim}\,}}

\newcommand{\Pscolim}[1]
{\underset{#1}{\mathrm{2\text{-}colim}\,}}

\newcommand{\liminv}[1]
{\underset{#1}{\lim}}

\renewcommand{\labelenumi}{(\roman{enumi})}

\newenvironment{pfclaim}[1][$\Diamond$]
{\def\claimQED{{#1}}\noindent {\em Proof of the claim. }}
{\leavevmode\unskip\penalty9999 \hbox{}\nobreak\hfill
    \quad\hbox{\claimQED}{\smallskip}}

\def\ad{\mathrm{ad}}
\def\gr{\mathrm{gr}}
\def\localg{\mathrm{loc.alg}}
\def\alg{\mathrm{alg}}

\def\Gal{\mathrm{Gal}}
\def\TopGrp{\mathbf{Top.Grp}}
\def\Spa{\mathrm{Spa}}
\def\Spec{\mathrm{Spec}}

\def\Hom{\mathrm{Hom}}

\def\Ind{\mathrm{Ind}}

\def\Ker{\mathrm{Ker}}

\def\Aut{\mathrm{Aut}}

\def\et{\mathrm{\acute{e}t}}

\def\Mod{\text{-}\mathbf{Mod}}
\def\bMod{\text{-}\mathbf{b.Mod}}

\def\Set{\mathbf{Set}}

\def\rk{\mathrm{rk}}
\def\one{\mathbf{1}}
\def\bar#1{\overline{#1}}
\def\hat{\widehat}
\def\tilde{\widetilde}
\def\fd{\mathfrak d}
\def\fm{\mathfrak m}

\def\eps{\varepsilon}

\renewcommand\emptyset{\varnothing}
\renewcommand\phi{\varphi}
\renewcommand\theta{\vartheta}

\def\cHom{\mathscr{H}\!om}

\def\cIsom{\mathscr{I}\!som}

\def\fn{\mathfrak{n}}

\def\length{\mathrm{length}}

\def\isom{\stackrel{\sim}{\to}}

\def\Ob{\mathrm{Ob}}

\def\Fil{\mathrm{Fil}}

\def\bCov{\mathbf{Cov}}

\def\bmu{\boldsymbol{\mu}}

\def\alphaenu{\renewcommand{\labelenumi}{(\alph{enumi})}}

\def\La{\langle}
\def\Ra{\rangle}

\def\ssw{\mathsf{sw}}
\def\bLambda{\mathbf{\Lambda}}

\makeindex

\begin{document}
\title{Hasse-Arf filtrations in $p$-adic analytic geometry}

\author{Lorenzo Ramero}

\maketitle


\tableofcontents

\section{Introduction}
This work is a sequel to \cite{Ram}, and it also draws
on my older paper \cite{Ram0}. Its aim is to prove that
the break decomposition for a $\Lambda$-module with bounded
ramification defined on a punctured (non-archimedean) disc
-- which was the main result of \cite{Ram} -- is already
defined on a punctured disc of possibly smaller radius
(here $\Lambda$ is a suitable ind-finite ring : see
\eqref{subsec_assumption-Lambda}).

The significance of this refinement is better conveyed, if
we restate it in terms of the category
$$
\Lambda[\pi(z)]\bMod
$$
of {\em germs} of $\Lambda$-modules with bounded ramification
(around a given closed point $z\in\P^1$ : see
\eqref{subsec_see-here}); namely, it implies that
$\Lambda[\pi(z)]\bMod$ is a tensor category {\em with a
Hasse-Arf filtration}, in the sense of Y.Andr\'e's paper \cite{And}.
Thus, the tannakian machinery of \cite{And} is now available,
to complete the classification of germs of modules with bounded
ramification : see section \ref{sec_break-zero} for some indication
about the work that remains to be done, in order to achieve this goal.

The plan of the proof is to chop such a module
along its breaks, using certain operators (the local Fourier
transforms), whose behaviour can be analyzed very precisely
by a local-to-global method (the principle of stationary phase).
Of course, all this is directly inspired by the article \cite{Lau2}
of G.Laumon, where the same sort of constructions are made, but
in the world of $k$-schemes, where $k$ is a perfect field of positive
characteristic.

Chapter \ref{chap_harmonic} implements the first step of this
programme; namely we prove the main finiteness results concerning
the (global) Fourier transform of modules with bounded ramification
everywhere on the affine line (on a non-archimedean field of
characteristic zero), then we construct our local Fourier
transforms and establish the principle of stationary phase.

Chapter \ref{chap_monodro} completes the programme : our
break decomposition is achieved in theorem \ref{th_main};
however, this is not the end of the story, since currently
we do not even know how to classify the modules of break
zero : what we know about this case is gathered in section
\ref{sec_break-zero}.

The global Fourier transform used here is the same that was
introduced in \cite{Ram0}, and a construction of local Fourier
transforms was also already attempted there, with modest
success. However, \cite{Ram0} was only concerned with a
special class of examples (the modules with meromorphic
ramification), defined by an ad hoc condition, whereas
now we are dealing with modules with bounded ramification,
which form {\em the} natural class for which a general study
of local monodromy can reasonably be carried out.

More generally, one may wonder whether the familiar tools of
perverse $\ell$-adic sheaves and weight filtrations that Deligne
introduced in his proof of the Weil conjectures, can be transplanted
in the setting of non-archimedean analytic geometry. We are still
several key insights away from being able to answer such kind
of questions; however, it is already clear that the condition of
bounded ramification must intervene in the definition of any
category of perverse sheaves which is stable under the usual
cohomological operations, and under Fourier transform. For
instance, the cohomology of a Zariski constructible module
$M$ on a quasi-compact analytic curve $C$ has finite length,
if and only if $M$ has bounded ramification at every point of $C$.

The article is structured as follows. In section \ref{sec_loc-alg-cov}
we introduce a wide class of \'etale coverings of analytic (more
precisely, adic in the sense of R.Huber) spaces; the definition
is actually found already in \cite{Ram0}, but the current treatment
is more general, systematic and tidier.

Section \ref{sec_vanish} is devoted to the construction and
basic properties of the vanishing cycles which later are used
to define the local Fourier transforms. Again, this parallels
what is done in \cite{Ram0}, but many unnecessary complications
have been removed, and some inaccuracy has been repaired.
Section \ref{sec_Fourier} recalls the (global) Fourier transform
of \cite{Ram0}; the main new theorem here (theorem \ref{th_perversity})
morally states that the Fourier transform of a sheaf with bounded
ramification is a perverse sheaf with bounded ramification everywhere.
(Except that we do not really try to define what a perverse sheaf
is.)

Finally, the local Fourier transforms are introduced in section
\ref{sec_Station}, and our version of the principle of stationary
phase is proved (theorem \ref{th_Stationary-phase}). These functors
act on certain categories, whose objects should be thought of as
representations of inertia subgroups (of the fundamental
group of a punctured disc). In truth, in our situation these
inertia subgroups are rather elusive, so the actual definitions
are somewhat more complicated than in the algebraic geometric
case : see \eqref{subsec_more-trivialities}. The basic properties
of the local Fourier transforms are then established in section
\ref{sec_study} : they are wholly analogous to the ones found in
\cite{Lau2}. The proofs are similar, but usually more difficult;
partly, this is because we do not have a functorial local-to-global
extension of representations of the inertia subgroups (in positive
characteristic, the existence of such an extension is a theorem of
Gabber : see \cite{Kat}). Instead, we rely on a theorem of Garuti
(\cite{Gar}), which provides a non-functorial extension, on a
Zariski-open subset of the affine line (see \eqref{subsec_Garuti}).
This suffices, up to some extra contorsions, to derive everything
we need; however, it would be interesting to know whether Garuti's
methods can be strengthened enough, to give a functorial extension
which would be truly analogous to Gabber's theorem.

Section \ref{sec_main-theor} contains the proof of our break
decomposition; with its help we may then also complete the analysis
of the local Fourier transforms.

We conclude with an application to the question of the
{\em localization of the determinant of cohomology}. This refers
to a problem, first proposed by Deligne in \cite{Del}, which can be
stated as follows. Let $C$ be a quasi-projective curve defined over
a subfield $K_0$ of an algebraically closed field $K$, and set
$G:=\Gal(K/K_0)$; let also $F$ be an $\ell$-adic {\em lisse} sheaf
on $C$; then one wants to decompose the $G$-representation
$$
\det(R\Gamma_c(C\times_{K_0}K,F))
$$
as a tensor product of {\em local $\eps$-factors}. The latter should
be certain characters $\eps(x,F)$ of $G$, attached to the points $x$
of the compatification of $C$, and depending functorially on the local
monodromy of $F$ around $x$. Deligne was motivated by the related
problem of the factorization of the $\eps$-constant appearing in
the functional equation of the $L$-function associated to $F$ when
$K_0$ is a finite field. And indeed, the latter arithmetic problem
was later solved by Laumon in \cite{Lau2}, by reducing it to the
cohomological localization problem, and proving a version of the
sought tensor decomposition in case where $C$ is an open subset
of the affine line over a finite field.

More recently, there has been a renewed interest in this problem :
for instance, the article \cite{BBE} establishes an analogue of
Laumon's theorem for $\cD$-modules over a field $K_0$ of
characteristic zero. Also, in Beilinson's paper \cite{Bei1} one
finds a real-analytic counterpart, that actually works for manifolds
of arbitrary dimension. And the very recent preprint \cite{Bei2},
also by Beilinson, establishes a general formalism that allows to
compare the deRham $\eps$-factors from \cite{BBE} with the
topological ones of \cite{Bei1}.

In section \ref{subsec_application} we establish a $p$-adic
analytic analogue of Laumon's formula, for the case where $K_0$
is a local field of zero characteristic. This is, in some sense,
an obvious application of our work, since the proof of \cite{Lau2}
uses the full apparatus of local Fourier transforms and stationary
phase. However, I should stress that only now our theory has
matured enough for such an application to lie within our grasp :
neither the results of our previous \cite{Ram}, nor the rudimentary
local Fourier transforms of \cite{Ram0} would have sufficed in order
to adapt the arguments of \cite{Lau2}.

We prove our tensor decomposition only for local systems of
$\bar\F_\ell$-modules (rather than $\ell$-adic sheaves, but see
remark \ref{rem_end-of-story}); on the other hand, we prove it
for general local systems with bounded ramification on the analytic
\'etale site of a Zariski open subset $U$ of the affine line
(heuristically, these sheaves correspond to holonomic $\cD$-modules
with arbitrary meromorphic singularities).
This includes the category of local systems coming from the algebraic
\'etale site of $U$ (which correspond to $\cD$-modules with regular
singularities), and in this case we rejoin Deligne's original
problem; our results appear to be new even for this special case,
and even for this restricted class of sheaves, our proof employs
essentially all the arsenal amassed in the previous sections.

An unexpected feature of our decomposition, is that it takes
place in the category of Galois characters that are
{\em semilinear\/} with respect to a certain action of the
Galois group of $K$ on $\bar\F_\ell$. We refer to remark
\ref{rem_unusual} for some observations concerning this kind
of hybrid, partly $p$-adic, partly $\ell$-adic representations.

\medskip

{\bf Acknowledgements} :\ \ 
I thank Roland Huber and Vadim Vologodsky for useful discussions.
This work has been written during a few stays at the
Max-Planck-Institut in Bonn, a visit at the Center for Advanced
Studies in Beer-Sheba, and a three-month invitation at the
Dipartimento di Matematica of the Universit\`a di Padova.

\section{Harmonic analysis for modules with bounded ramification}
\label{chap_harmonic}

\subsection{Locally algebraic coverings}\label{sec_loc-alg-cov}
Throughout this work, $(K,|\cdot|)$ is a fixed algebraically closed
valued field, complete for a rank one valuation and of residue
characteristic $p>0$, and we denote by $K^+$ the valuation ring of
$|\cdot|$. We shall use freely the language of adic spaces from
\cite{Hu2}, and we set :
$$
S:=\Spa(K,K^+).
$$
For a $K$-scheme $X$ locally of finite type, we let $X^\ad$ be the
adic space associated to $X$ (this is denoted $X\times_{\Spec\,K}S$
in \cite[Prop.3.8]{Hu4}).
\begin{definition}
A morphism $f:Y\to X$ of analytic $S$-adic spaces is said to be
a {\em locally algebraic covering}, if the following holds. For
every quasi-compact open subset $U\subset X$, the preimage
$f^{-1}U$ decomposes as the (possibly infinite) disjoint union
of open and closed subspaces, each of which is a finite \'etale
covering of $U$. We denote by
$$
\bCov^\alg(X)
\qquad
\text{(resp. $\bCov^\localg(X)$)}
$$
the full subcategory of the category of $X$-adic spaces, whose
objects are the finite \'etale (resp. the locally algebraic)
coverings.
\end{definition}

\begin{remark}\label{rem_loc-alg-covs}
(i)\ \
It is immediate from the definition, that every locally
algebraic covering is a separated morphism.

(ii)\ \
It is easily seen that a fibre product $Y_1\times_XY_2$ of
locally algebraic covering spaces of $X$, is also a locally
algebraic covering of $X$. If $f:Y\to X$ is a locally algebraic
covering, and $X'\to X$ is any morphism of $S$-adic spaces, then
$f\times_X\one_{X'}:Y\times_XX'\to X'$ is a locally algebraic covering.

(iii)\ \
Likewise, if $f:Z\to Y$ and $g:Y\to X$ are locally algebraic
coverings, then the same holds for $g\circ f:Z\to X$. And if
$(Y_i\to X~|~i\in I)$ is any family of locally algebraic coverings,
then the disjoint union $\coprod_{i\in I}Y_i\to X$ is a locally
algebraic covering as well. (These two properties do not hold for
general \'etale coverings, defined as in \cite[Def.2.1]{deJ} : in
{\em loc. cit.} the language of Berkovich's analytic spaces is used,
but the definition and the proofs of several of the preliminary
results can be repeated {\em verbatim} for adic spaces).

(iv)\ \ 
Morever, let $f_i:Y_i\to X$, for $i=1,2$, be two locally
algebraic coverings, and $g:Y_1\to Y_2$ a morphism of $X$-adic
spaces; then $g$ is a locally algebraic covering of $Y_2$.
Indeed, let $U\subset Y_2$ be a quasi-compact open subset;
after replacing $X$ by $f_2U$, and $Y_i$ by $f^{-1}_if_2U$
($i=1,2$) we may assume that $X$ is quasi-compact, hence
$Y_1=\bigcup_{\alpha\in I}Z_\alpha$ for a family of finite
and \'etale $X$-adic spaces $Z_\alpha$. Denote by
$g_\alpha:Z_\alpha\to Y_2$ the restriction of $g$; we may write
$g_\alpha=p_\alpha\circ s_\alpha$, where
$p_\alpha:Z_\alpha\times_XY_2\to Y_2$ is the projection (which
is a finite \'etale morphism) and
$s_\alpha:Z_\alpha\to Z_\alpha\times_XY_2$ is the graph of
$g_\alpha$, which is a closed immersion, in view of (i) and
a standard argument. Hence $g_\alpha$ is a finite \'etale
morphism. However, $g^{-1}U=\bigcup_{\alpha\in I}g_\alpha^{-1}U$,
and the foregoing shows that each $g_\alpha^{-1}U$ is a finite
\'etale $U$-adic space, whence the assertion.

(v)\ \ 
If $X$ is quasi-compact, we have a natural equivalence of
categories :
$$
\Ind(\bCov^\alg(X))\isom\bCov^\localg(X)
$$
where, for a category $\cC$, we denote by $\Ind(\cC)$ the category
of ind-objects of $\cC$, defined as in \cite[Exp.I, \S8.2.4]{SGA4-1}.

(vi)\ \ 
On the other hand, local algebraicity is {\em not} a property
local on $X$ : for instance, let $E$ be an elliptic curve over
$K$ with bad reduction over the residue field of $K^+$, and
consider the analytic uniformization $p:\G_m\to E$. This is
not a locally algebric covering of $E$, but there exist a
covering $E=X_1\cup X_2$ by open (quasi-compact) subspaces,
such that the two restrictions $p^{-1}X_i\to X_i$ are both
locally algebraic.
\end{remark}

\begin{lemma}\label{lem_quot-equiv-rel}
Let $Y\to X$ be a locally algebraic covering, and $R$ an equivalence
relation on $Y$ which is represented by a union of open and closed
subspaces of $Y\times_XY$. Then the (categorical) quotient sheaf
$Y/R$ is represented by a locally algebraic covering.
\end{lemma}
\begin{proof} Let $U'\subset U\subset X$ be any inclusion of open
subsets of $X$, and suppose that both restrictions $(Y/R)_{|U}$
and $(Y/R)_{|U'}$ are represented by morphisms $Z\to U$, respectively
$Z'\to U'$; we deduce easily a natural isomorphism
$Z\times_UU'\isom Z'$, compatible with inclusion of smaller open
subsets $U''\subset U'$. It follows that the lemma holds for the
datum $(X,Y,R)$ if and only if it holds for every datum of
the form $(U,f^{-1}U,R\times_XU)$ where $U\subset X$ is an
arbitrary quasi-compact open subset. Hence we may assume that $X$ is
quasi-compact. In that case, $Y=\bigcup_{i\in I}Y_i$ for a family
$(Y_i~|~i\in I)$ of finite \'etale $X$-spaces, such that $Y_i$
is open and closed in $Y$, for every $i\in I$. For every finite
subset $S\subset I$, let $Y_S:=\bigcup_{i\in S}Y_i$; it follows easily
that $Y$ is the (categorical) colimit of the filtered system of
open and closed immersions $Y_{S'}\to Y_S$, where $S'\subset S$
range over the finite subsets of $I$. For every such $S$, let
$R_S:=R\cap(Y_S\times_XY_S)$. Then $R$ is likewise the colimit of the
open and closed immersions $R_{S'}\to R_S$ for $S$, $S'$ as above.
Finally, $Y/R$ (resp. $Y_S/R_S$) is the coequalizer of the two
projections $R\to Y$ (resp. $R_S\to Y_S$); by general nonsense,
we deduce that $Y/R$ is naturally isomorphic to the colimit of
the filtered system of open and closed immersions
$Y_{S'}/R_{S'}\to Y_S/R_S$. We are therefore reduced to showing that
each quotient $Y_S/R_S$ is representable, hence we may assume from
start that $Y\to X$ is a finite \'etale covering. In this case,
the lemma can be checked locally on $X$, hence we may further
assume that $X$ is affinoid, say $X=\Spa(A,A^+)$; then there
exists a finite \'etale $A$-algebra $B$ such that
$Y=\Spec\,B\times_{\Spec\,A}X$. Likewise, $R=Z\times_{\Spec\,A}X$
for some open and closed equivalence relation
$Z\subset\Spec\,B\otimes_AB$. We may then represent $Y/R$ by the
$X$-adic space $(\Spec\,B/Z)\times_{\Spec\,A}X$.
\end{proof}

\sset\subsubsection{}\label{subsec_fix-a-max}
Fix a {\em maximal geometric point\/} of $X$, {\em i.e.} a morphism
$\xi:\Spa(E,E^+)\to X$ of $S$-adic spaces, where $(E,|\cdot|)$ is a
complete and algebraically closed valued field extension of $K$,
with valuation group of rank one.
To $\xi$ we associate a {\em fibre functor}
$$
\sF_{X,\xi}:\bCov^\localg(X)\to\Set
\quad :\quad
(Y\to X)\mapsto\Hom_X(\Spa(E,E^+),Y)
$$
(unless we are dealing with more than one adic space, we shall usually
drop the subscript $X$, and write just $\sF_\xi$).
The {\em locally algebraic fundamental group\/} of $X$, pointed
at $\xi$, is defined as usual, as the automorphism group of $\sF_\xi$ :
$$
\pi_1^\localg(X,\xi):=\Aut(\sF_\xi).
$$
It is endowed with a natural topology, as in \cite[\S2]{deJ}.
Namely, for every pair of the form $(Y,\bar y)$, with $Y$ a
locally algebraic covering of $X$ and $\bar y\in\sF_\xi(Y)$, let
$H(Y,\bar y)\subset\pi_1^\localg(X,\xi)$ be the stabilizer
of $\bar y$ (for the action on $\sF_\xi(Y)$); then the family $\cH$
of such subgroups is stable under finite intersections (by remark
\ref{rem_loc-alg-covs}(ii)) and under conjugation by arbitrary
elements $\gamma\in\pi_1^\localg(X,\xi)$, since
$\gamma H(Y,\bar y)\gamma^{-1}=H(Y,\gamma\cdot\bar y)$.
Hence, there is a unique topology on $\pi_1^\localg(X,\xi)$
for which the family $\cH$ forms a fundamental system of open
neighborhoods of the identity element.

Arguing as in the proof of \cite[Lemma 2.7]{deJ} we see that
$\pi_1^\localg(X,\xi)$ is Hausdorff and prodiscrete, and more
precisely, the natural map :
$$
\pi_1^\localg(X,\xi)\to\lim_{H\in\cH}\pi_1^\localg(X,\xi)/H
$$
is an isomorphism of topological groups (where the target is
the (projective) limit computed in the category of topological
groups). Also, the fibre functor can be upgraded, as usual, to
a functor 
\set\begin{equation}\label{eq_upgrade}
\sF_\xi:\bCov^\localg(X)\to\pi_1^\localg(X,\xi)\text{-}\Set
\end{equation}
with values in the category of (discrete) sets with continuous
left action of the locally algebraic fundamental group. The continuity
condition means that the stabilizer of any point of such a set,
is an open subgroup.

\begin{proposition}\label{prop_pi-sets}
Let $X$ be a connected analytic $S$-adic space, and $\xi$ a
maximal geometric point of $X$. Then the functor \eqref{eq_upgrade}
is an equivalence.
\end{proposition}
\begin{proof} Taking into account remark \ref{rem_loc-alg-covs}(iii)
and lemma \ref{lem_quot-equiv-rel}, the proof of \cite[Th.2.10]{deJ}
can be taken over {\em verbatim} : the crucial point is the following

\begin{claim}\label{cl_iso-of-fibre-fctrs}
In the situation of the theorem, suppose that $\xi$ and $\xi'$ are
two maximal geometric points of $X$. Then there exists an isomorphism
of functors $\sF_\xi\isom\sF_{\xi'}$.
\end{claim}
\begin{pfclaim}[] For every point $x\in X$, let $U(x)\subset X$
be the union of all the quasi-compact connected open neighborhoods
of $x$ in $X$. It is easily seen that $U(x)\cap U(x')\neq\emptyset$
if and only if $U(x)=U(x')$; since $X$ is connected, it follows
that $U(x)=X$ for every $x\in X$. Thus, we may find a connected
quasi-compact open subset $U\subset X$ containing the supports
of both $\xi$ and $\xi'$. The restrictions of $\sF_{U,\xi}$ and
$\sF_{U,\xi'}$ yield fibre functors 
$$
\sF_{U,\xi}^\alg,\sF^\alg_{U,\xi'}:\bCov^\alg(U)\to\Set
$$
which fulfill the axiomatic conditions for a Galois theory,
as in \cite[Exp.V, \S4]{SGA1}; there follows an isomorphism
of functors $\sF^\alg_{U,\xi}\isom\sF^\alg_{U,\xi'}$
(\cite[Exp.V, Cor.5.7]{SGA1}). In view of remark
\ref{rem_loc-alg-covs}(v), we deduce an isomorphism between
the two fibre functors
$$
\sF_{U,\xi},\sF_{U,\xi'}:\bCov^\localg(U)\to\Set.
$$
However, $\sF_{X,\xi}=\sF_{U,\xi}\circ j^*$ and
$\sF_{X,\xi'}=\sF_{U,\xi'}\circ j^*$, where
$j^*:\bCov^\localg(X)\to\bCov^\localg(U)$ is the restriction
functor (for the open immersion $j:U\to X$). The claim follows.
\end{pfclaim}
\end{proof}

\sset\subsubsection{}
Let $(G,\cdot)$ be a group; we denote by $G^o$ the group {\em opposite}
to $G$, whose elements are the same as those of $G$, and whose
composition law is given by the rule : $(g_1,g_2)\mapsto g_2\cdot g_1$
for all $g_1,g_2\in G$. A {\em right $G$-torsor\/} (resp. a {\em right
$G$-principal homogeneous space}) on $X$ is, as usual, a sheaf
$F$ on $X_\et$ (resp. a morphism $Y\to X$ of analytic adic spaces),
with a group homomorphism $G^o\to\Aut(F)$ (resp. $G^o\to\Aut_X(Y)$),
for which there exists a covering family $(U_i\to X~|~i\in I)$
in $X_\et$, such that the restriction $F_{|U_i}$ (resp.
$(U_i\times_XY\to U_i)$), together with the induced right $G$-actions,
is isomorphic to the constant sheaf $G_{U_i}$ (resp. is $U_i$-isomorphic
to $U_i\times G$) for every $i\in I$. As usual, any right $G$-principal
homogeneous space on $X$ represents a right $G$-torsor, unique up to
isomorphism.

\begin{lemma}
Let $f:Y\to X$ be a $G$-principal space on $X_\et$. Then the
morphism of $Y$-adic spaces
\set\begin{equation}\label{eq_themorph}
Y\times G\to Y\times_XY \qquad : \qquad (y,g)\mapsto(y,yg)
\end{equation}
is an isomorphism. Especially, $f$ is separated.
\end{lemma}
\begin{proof} Let $F$ denote the $G$-torsor on $X_\et$ represented
by $Y$; then \eqref{eq_themorph} represents the morphism of sheaves
$\phi:F\times G\to F\times F$ on $X_\et$ given by the rule :
$(s,g)\mapsto(s,sg)$.
By assumption, there exists a covering family $(U_i\to X~|~i\in I)$
such that $\phi_{|U_i}$ is an isomorphism for every $i\in I$; hence
$\phi$ is an isomorphism, and then the same holds for
\eqref{eq_themorph}. It follows that the diagonal immersion
$\Delta_{Y/X}\to Y\times_XY$ can be identified to the closed
immersion $Y\to Y\times G$ which maps $Y$ onto the closed subspace
$Y\times\{1\}$ (where $1\in G$ is the neutral element); hence $f$
is separated.
\end{proof}

One defines {\em left $G$-torsors}, and {\em left $G$-principal
homogeneous spaces} in the same way, except that instead of $G^o$
one has the group $G$. The set of isomorphism classes of right
$G$-torsors on $X$ is denoted 
$$
H^1(X_\et,G).
$$
If $G$ is abelian, this is a group, naturally isomorphic to the
\v{C}ech cohomology group $\check{H}^1(X_\et,G)$. By general nonsense,
the latter is isomorphic to the cohomology group $H^1(X_\et,G)$
calculated by injective resolutions (\cite[Ch.II, Cor.3.6]{Art}).

If $G\to G'$ is any group homomorphism, and $F$ is a $G$-torsor
on $X$, we obtain a $G'$-torsor 
$$
F\overset{G}{\times} G'
$$
such that $F\overset{G}{\times}G'(Z)=(F(Z)\times G')/G$ for any
\'etale $X$-adic space $Z$ (where $G$ acts on $F(Z)\times G'$ by
the rule : $(g,(s,g'))\mapsto(sg,g^{-1}g')$). This operation
defines a natural transformation :
$$
H^1(X_\et,G)\to H^1(X_\et,G').
$$

\begin{proposition}\label{prop_represent-ind-finite}
Let $G$ be an ind-finite group ({\em i.e.} $G$ is the union of
the filtered family $\cF$ of its finite subgroups), $X$ a
$S$-adic space, and $F$ a right $G$-torsor on $X$. Then :
\begin{enumerate}
\item
$F$ is represented by a right $G$-principal homogeneous space
$Y\to X$.
\item
If $X$ is quasi-separated, $Y\to X$ is a locally algebraic covering.
\item
If $X$ is quasi-compact and quasi-separated, the natural map :
$$
\colim{H\in\cF}H^1(X_\et,H)\to H^1(X_\et,G)
$$
is a bijection.
\end{enumerate}
\end{proposition}
\begin{proof}(iii): This is proved as in \cite[Lemma 4.1.2]{Ram},
which is a special case. Since one assumption is missing in
{\em loc.cit.}, we repeat briefly the argument here.

One covers $X$ by finitely many affinoids $U_1,\dots,U_n$ which
admit finite \'etale coverings $V_i\to U_i$, such that $F_{|V_i}$
is the trivial $G$-torsor on $V_i$, for every $i\leq n$.
The latter is endowed with a natural descent datum, given by a
cocycle consisting of a locally constant function
$g_i:V_i\times_{U_i}V_i\to G$.
By quasi-compactness, such a function assumes only finitely
many values; we may therefore find a finite subgroup $H\subset G$
containing the image of all maps $g_i$, and it follows that
the isomorphism class of the restriction $F_{|U_i}$ lies in
the image of $H^1(U_{i,\et},H)$.
Lastly, the gluing datum for the $G$-torsor $F$, relative
to the covering $(U_i~|~i=1,\dots,n)$ amounts to another
cocycle $(g_{ij}:U_i\cap U_j\to G~|~i,j\leq n)$ of the same
type. Since $X$ is quasi-separated, $U_i\cap U_j$ is again
quasi-compact, hence again the class of $Y$ is the image of
some class in $H^1(X_\et,H')$, for some finite group
$H'$ containing $H$. This shows the surjectivity of our
map; the injectivity is proved in a similar way : details
left to the reader.

(i): The assertion is local on $X$, hence we may assume that
$X$ is quasi-compact and quasi-separated, in which case, by (iii),
$F\simeq F'\overset{H}{\times}G$, for some finite subgroup
$H\subset G$, and a $H$-torsor $F'$; according to \cite[\S2.2.3]{Hu2},
$F$ is representable by a $H$-principal homogeneous space
$Y_H\to X$. It follows easily that $F$ is representable by 
the $G$-principal homogeneous space $Y_H\overset{H}{\times}G$.

(ii): If $X$ is quasi-separated, the foregoing proof of (i)
works on every quasi-compact open subset of $X$, and a simple
inspection shows that the construction yields a locally
algebraic covering in this case.
\end{proof}

\sset\subsubsection{}\label{subsec_G-sets}
Let $X$ be a connected quasi-separated $S$-adic space, $\xi$
a maximal geometric point of $X$, $G$ an ind-finite group, and
$f:Y\to X$ a right $G$-principal homogeneous space.
By propositions \ref{prop_represent-ind-finite}(ii) and
\ref{prop_pi-sets}, the isomorphism class of $Y$ is determined
by the fibre $f^{-1}(\xi)$ together with its natural left
$\pi^\localg_1(X,\xi)$-action and right $G$-actions, which
commute with each other. Conversely, every set $S$ endowed
with commuting left $\pi^\localg_1(X,\xi)$-action and right
$G$-action, such that moreover the $G$-action is free and
transitive, arises as the fibre of a right $G$-torsor.

In this situation, it is easily seen that the stabilizer -- for
the left $\pi^\localg_1(X,\xi)$-action -- of any point
$y\in f^{-1}(\xi)$ is an open normal subgroup $H$, and we obtain a
group homomorphism :
\set\begin{equation}\label{eq_left-qnd-right}
\pi^\localg_1(X,\xi)/H\to G
\end{equation}
by assignining to the class $\bar\sigma$ of any
$\sigma\in\pi^\localg_1(X,\xi)$ the unique $g\in G$ such that
$\sigma(y_0)=y_0\cdot g$.

The effect of replacing $y_0$ by a different point, is to
change \eqref{eq_left-qnd-right} by an inner automorphism
of $G$; summing up, we deduce a natural bijection :
$$
H^1(X_\et,G)/\isom\Hom_\TopGrp(\pi^\localg_1(X,\xi),G)/\!\!\sim
$$
onto the set of equivalence classes of continuous group homomorphisms
as in \eqref{eq_left-qnd-right} (where $G$ is endowed with the
discrete topology), for the equivalence relation $\sim$ induced
by inner conjugation on $G$. Finally, under this bijection, the
isomorphism classes of connected $G$-torsors correspond precisely
to the equivalence classes of continuous surjective group
homomorphisms.

\sset\subsubsection{}\label{subsec_was-remark}
Let $\Lambda$ be an ind-finite ring, $X$ a connected $S$-adic
space, $\xi$ a maximal geometric point of $X$, and $F$ a locally
constant $\Lambda$-module of finite type on $X_\et$.
Set $M:=F_\xi$; denote by
$$
\cIsom_\Lambda(M_X,F)
$$
the sheaf on $X_\et$ whose sections on any \'etale $X$-adic space
$Z$ are the isomorphisms $M_Z\isom F_{|Z}$ of $\Lambda$-modules
on $Z_\et$. Clearly $\cIsom_\Lambda(M_X,F)$ is a right $G$-torsor, for
$G:=\Aut_\Lambda(M)$; by proposition \ref{prop_represent-ind-finite}(ii)
it follows that $F$ is trivialized on a locally algebraic covering
$g:Y\to X$, corresponding to a certain open normal subgroup of
$\pi^\localg_1(X,\xi)$. The natural map $\Gamma(Y,F)\to M$ is an
isomorphism, therefore the natural action of $\pi^\localg_1(X,\xi)$
on $\Gamma(Y,F)$ (cp. \eqref{subsec_vanish}) can be transferred
to $M$. This action is {\em continuous}, {\em i.e.} the stabilizer
of every element of $M$ is an open subgroup. In this way we obtain
an equivalence :
\set\begin{equation}\label{eq_belated}
\Lambda_X\Mod_\mathrm{loc}\isom
\Lambda[\pi^\localg_1(X,\xi)]\Mod_\mathrm{f.cont}
\qquad
F\mapsto F_\xi
\end{equation}
from the category of locally constant $\Lambda$-modules of finite
type on $X_\et$, to the category of $\Lambda$-modules of finite type,
endowed with a continuous linear action of $\pi^\localg_1(X,\xi)$.

\subsection{Vanishing cycles}\label{sec_vanish}
For every $\eps\in|K^\times|$, denote by $\D(z,\eps)$ the disc of
radius $\eps$, centered at a $K$-rational point $z\in(\P^1_K)^\ad$.
Then $\D(z,\eps)$ is an open analytic subspace of $(\P^1_K)^\ad$, and
$$
\D(z,\eps)^*:=\D(z,\eps)\setminus\{z\}
$$
is an open subspace of $\D(z,\eps)$. Choose a maximal geometric point
$\xi\in\D(z,\eps)^*$ and set 
$$
\pi(z,\eps):=\pi_1^\localg(\D(z,\eps)^*,\xi).
$$
Moreover, for every open subgroup $H\subset\pi(z,\eps)$,
choose a locally algebraic covering
$$
\phi_H:C_H\to\D(z,\eps)^*
$$
with an isomorphism $\sF_\xi(C_H)\isom\pi(z,\eps)/H$ of
$\pi(z,\eps)$-sets (where $\sF_\xi$ is the fibre functor
as in \eqref{eq_upgrade}). The family
$$
(C_H~|~H\subset\pi(z,\eps))
$$
is a cofiltered system of analytic adic spaces (in the terminology
of \cite[Exp.V, \S5]{SGA1}, this is a {\em fundamental pro-covering\/}
of $\D(z,\eps)^*$).

\sset\subsubsection{}\label{subsec_vanish}
Let now $\Lambda$ be an ind-finite ring such that $p\in\Lambda^\times$,
and $f:X\to\D(z,\eps)$ a morphism of analytic $S$-adic spaces.
For every open subgroup $H\subset\pi(z,\eps)$, let
$j_H:C_H\times_{\D(z,\eps)}X\to X$ be the projection;
set $X_0:=f^{-1}(z)$, and let $i:X_0\to X$ be the closed immersion.
To any lower bounded complex $\cF^\bullet$ of $\Lambda$-modules on
$X_\et$, we associate the complex of $\Lambda$-modules
$$
R\Psi_{\eta_z,\eps}\cF^\bullet:=
\colim{H\subset\pi(z,\eps)}\,i^*Rj_{H*}j_H^*\cF^\bullet
\qquad
\text{on $X_{0,\et}$.}
$$
We construct as follows an action of $\pi(z,\eps)$ on the
cohomology sheaves $R^i\Psi_{\eta_z,\eps}\cF^\bullet$ of this complex.
To begin with, let $F$ be any sheaf on $X_\et$, and set 
$$
\tilde F:=\colim{H\subset\pi(z,\eps)}\,j_{H*}j_H^*F.
$$
Let $g\in\pi(z,\eps)$ be any element, and $H\subset\pi(z,\eps)$ an
open subgroup; the right translation action of $g$ on $\pi(z,\eps)$
induces a bijection of pointed sets :
$$
\pi(z,\eps)/H\isom\pi(z,\eps)/g^{-1}Hg
\quad :\quad
[\gamma H]\mapsto[(\gamma g)g^{-1}Hg]
$$
which is clearly equivariant for the left action of $\pi(z,\eps)$
on both sets, and therefore corresponds to a unique isomorphism
of $\D(z,\eps)^*$-adic spaces :
$$
g_H:C_H\isom C_{g^{-1}Hg}
$$
whence an isomorphism
\set\begin{equation}\label{eq_for-g-and-H}
F(U\times_{\D(z,\eps)}g_H):
F(U\times_{\D(z,\eps)}C_{g^{-1}Hg})\isom F(U\times_{\D(z,\eps)}C_H)
\end{equation}
for every object $U\to X$ of $X_\et$. Moreover, if $H\subset H'$
for any other open subgroup $H'\subset\pi(z,\eps)$, we have a
commutative diagram of $S$-adic spaces :
$$
\xymatrix{ C_H \ar[r]^-{g_H} \ar[d] & C_{g^{-1}Hg} \ar[d] \\
           C_{H'} \ar[r]^-{g_{H'}}  & C_{g^{-1}H'g} 
}$$
whose vertical arrows are induced by the quotient map
$\pi(z,\eps)/H\to\pi(z,\eps)/H'$.
Therefore the automorphisms \eqref{eq_for-g-and-H}, for fixed
$g$ and varying $H$, organize into a directed system.
Suppose that $U$ is affinoid; then, taking colimits, we obtain
an automorphism :
$$
\tilde g(U):=\colim{H\subset\pi(z,\eps)}\,
F(U\times_{\D(z,\eps)}g_H) : \tilde F(U)\isom\tilde F(U)
$$
(since $U$ is quasi-compact and quasi-separated, the functor
$\Gamma(U,-)$ commutes with filtered colimits). Clearly the
maps $\tilde g(U)$ for variable $U$, assemble to a well defined
automorphism 
$$
\tilde g:\tilde F\isom\tilde F
$$
which, for variable $g$, defines a {\em left\/} $\pi(z,\eps)$-action
on $\tilde F$ (the controvariance of $F$ transforms the right
translation action into a left action).

In case $F$ is a sheaf of groups, it is clear that $\pi(z,\eps)$
acts by group automorphisms on $\tilde F$. To define the action
on $R^i\Psi_{\eta_z,\eps}\cF^\bullet$, it suffices now to apply the
above constructions to the terms of an injective resolution
of $\cF^\bullet$.

\begin{remark}\label{rem_cont-action}
With the notation of \eqref{subsec_vanish}, notice that
$g_H$ is the identity automorphism of $C_H$, whenever $g\in H$.
It follows that, for every geometric point $\bar x$ of $X_0$,
the $\pi(z,\eps)$-action on the stalk $\tilde F_{\bar x}$ is
continuous in the sense of \eqref{subsec_was-remark}.
Hence, for every $i\in\Z$, the stalks of
$R^i\Psi_{\eta_z,\eps}\cF^\bullet$ are objects of the category
$$
\Lambda[\pi(z,\eps)]\Mod_\mathrm{cont}
$$
whose objects are all the $\Lambda$-modules endowed with a
continuous action of $\pi(z,\eps)$.
\end{remark}

\begin{proposition}\label{prop_vanish-properties}
In the situation of \eqref{subsec_vanish}, let $g:Y\to X$ be a
morphism of analytic $S$-adic spaces, set $Y_0:=g^{-1}X_0$, and
denote by $g_0:Y_0\to X_0$ the restriction of $g$. Let also
$\cF^\bullet$ (resp. $\cG^\bullet$) be a lower bounded complex
of $\Lambda$-modules on $X_\et$ (resp. on $Y_\et$). We have :
\begin{enumerate}
\item
If $g$ is smooth, there exists a natural isomorphism :
$$
g_0^*R\Psi_{\eta_z,\eps}\cF^\bullet\isom
R\Psi_{\eta_z,\eps}g^*\cF^\bullet
$$
in the derived category $\sD^+(Y_{0,\et},\Lambda)$ of lower
bounded complexes of $\Lambda$-modules on $Y_{0,\et}$.
\item
If $g$ is weakly of finite type and quasi-separated, there
exists a natural isomorphism :
$$
R\Psi_{\eta_z,\eps}(Rg_*\cG^\bullet)\isom
Rg_{0*}R\Psi_{\eta_z,\eps}\cG^\bullet
\qquad
\text{in $\sD^+(X_{0,\et},\Lambda)$.}
$$
\item
There is a spectral sequence :
$$
E^{pq}_2:=R^p\Psi_{\eta_z,\eps}(\cH^q\cF^\bullet)\Rightarrow
R^{p+q}\Psi_{\eta_z,\eps}\cF^\bullet.
$$
\item
The isomorphisms in {\em(i)} and {\em(ii)} induce
$\pi(z,\eps)$-equivariant isomorphisms on the cohomology
sheaves of the respective terms, and the spectral sequence
of {\em (iii)} is $\pi(z,\eps)$-equivariant.
\end{enumerate}
\end{proposition}
\begin{proof}(i): This follows easily from the base change
theorem \cite[Th.4.1.1(a)]{Hu2}, applied to the cartesian
diagram :
\set\begin{equation}\label{eq_smooth-bc}
{\diagram  C_H\times_{\D(z,\eps)}Y \ar[r]
           \ar[d]_{C_H\times_{\D(z,\eps)}g} & Y \ar[d]^g \\
           C_H\times_{\D(z,\eps)} X \ar[r]^-{j_H} & X.
\enddiagram}
\end{equation}
(ii): Likewise, this follows from \eqref{eq_smooth-bc} and
\cite[Th.4.1.1(c)]{Hu2} (since $j_H$ is generalizing, and
$g$ is weakly of finite type and quasi-separated), together
with the cartesian diagram :
$$
\xymatrix{ Y_0 \ar[r] \ar[d]_{g_0} & Y \ar[d]^g \\
           X_0 \ar[r]^-i & X
}$$
and again \cite[Th.4.1.1(c)]{Hu2} (since $g$ is weakly of finite
type and quasi-separated, and $i$ is generalizing). Notice here,
that the functor $Rg_*$ commutes with filtered colimits, since
$g$ is quasi-compact and quasi-separated
(\cite[Lemma 2.3.13(ii)]{Hu2}).

Finally, (iii) and (iv) are immediate from the construction.
\end{proof}

\begin{lemma}\label{lem_complete-later}
In the situation of \eqref{sec_vanish}, let $G$ be any ind-finite
group. We have :
$$
\colim{H\subset\pi(z,\eps)}H^1(C_{H,\et},G)=0.
$$
\end{lemma}
\begin{proof} Let $c\in H^1(C_{H,\et},G)$ for some
$C_H$ as above. The class $c$ represents a right $G$-torsor
$F$ on $C_{H,\et}$, and by proposition
\ref{prop_represent-ind-finite}(ii) we may find a connected
locally algebraic covering $Y\to C_H$ trivializing $F$.
In view of remark \ref{rem_loc-alg-covs}(iii), there exists
an open subgroup $L\subset H$ such that $C_L$ dominates
$Y$, which means that the image of $c$ vanishes in
$H^1(C_{L,\et},G)$.
\end{proof}

\sset\subsubsection{}\label{subsec_smaller-radius}
Let now $\eps'\in|K^\times|$ be a value with $\eps'<\eps$,
choose a maximal geometric point $\xi'\in\D(z,\eps')^*$, and
let $\pi(z,\eps')$ be the automorphism group of the corresponding
fibre functor $\sF_{\xi'}$ on the locally algebraic coverings of
$\D(z,\eps')^*$.

As usual, the open immersion $j:\D(z,\eps')^*\to\D(z,\eps)^*$
induces a group homomorphism
\set\begin{equation}\label{eq_choose-up-to-inner}
\rho:\pi(z,\eps')\to\pi(z,\eps)
\end{equation}
well-defined up to inner automorphisms; more precisely, $\rho$
is induced by a choice of isomorphism of fibre functors on
$\bCov^\localg(\D(z,\eps)^*)$ :
\set\begin{equation}\label{eq_choose-an-iso}
\sF_\xi\isom\sF_{\xi'}\circ j^*
\end{equation}
where $j^*:\bCov^\localg(\D(z,\eps)^*)\to\bCov^\localg(\D(z,\eps')^*)$
is the restriction functor.
Now, pick a fundamental pro-covering of $\D(z,\eps')^*$ :
$$
(C'_L~|~L\subset\pi(z,\eps'))
$$
as in \eqref{sec_vanish}. For every open subgroup
$H\subset\pi(z,\eps)$, the preimage $\rho^{-1}H\subset\pi(z,\eps')$
is an open subgroup, and the induced injective map
$$
\pi(z,\eps')/\rho^{-1}H\to\pi(z,\eps)/H
$$
corresponds to an open and closed immersion of $\D(z,\eps')$-adic
spaces :
\set\begin{equation}\label{eq_bH}
C'_{\rho^{-1}H}\to C_H\times_{\D(z,\eps)}\D(z,\eps').
\end{equation}
If $F$ is any sheaf on $X_\et$, define $\tilde F$ as in
\eqref{subsec_vanish}; let also $X':=f^{-1}\D(z,\eps')$,
and for every open subgroup $L\subset\pi(z,\eps')$, denote
by $j'_L:C'_L\times_{\D(z,\delta)}X(\eps')\to X(\eps')$
the projection. Set $F':=F_{|X'}$, and define :
$$
\tilde F':=\colim{L\subset\pi(z,\eps')}\,j'_{L*}j^{\prime*}_LF'.
$$
Notice that, for every open subgroup $H\subset\pi(z,\eps)$,
the morphism $j'_{\rho^{-1}H}$ can be written as the composition
of $b_H:=X'\times_{\D(z,\eps')}\eqref{eq_bH}$ and the
restriction $j''_H:C_H\times_{\D(z,\eps)}X'\to X'$ of $j_H$.
Then the unit of adjunction $\one\to b_{H*}b_H^*$ yields
a natural morphism :
$$
\tilde F_{|X'}\isom\colim{H\subset\pi(z,\eps)}\,
j''_{H*}b_{H*}b_H^*j^{\prime\prime*}_HF'\to\tilde F'
$$
which is equivariant for the left action of $\pi(z,\eps')$
defined above (where $\pi(z,\eps')$ acts on $\tilde F_{|\D(z,\eps')}$
by restriction of the $\pi(z,\eps)$-action, via
\eqref{eq_choose-up-to-inner}).
Finally, let $\cF^\bullet$ be a complex of $\Lambda$-modules on
$X_\et$; by applying this construction to the terms of an injective
resolution of $\cF^\bullet$, we deduce natural morphisms :
$$
R^i\Psi_{\eta_z,\eps}\cF^\bullet\to R^i\Psi_{\eta_z,\eps'}\cF^\bullet_{|X'}
\qquad
\text{for every $i\geq 0$}.
$$

\sset\subsubsection{}\label{subsec_pass-to-lim}
We wish now to ``take the limit'' for $\eps\to 0$, of the
previous construction.  To begin with, fix a strictly
decreasing sequence $(\eps_n~|~n\in\N)$ with $\eps_n\in|K^\times|$
for every $n\in\N$, and such that $\eps_n\to 0$ for $n\to+\infty$.
For every $n\in\N$, fix a geometric base point $\xi_n$ of
$\D(z,\eps_n)^*$ and a group homomorphism 
\set\begin{equation}\label{eq_restrict-radius}
\pi(z,\eps_{n+1})\to\pi(z,\eps_n)
\end{equation}
as in \eqref{eq_choose-up-to-inner}, as well as a fundamental
pro-covering of $\D(z,\eps_n)^*$ :
$$
(C^{(n)}_H~|~H\subset\pi(z,\eps_n)).
$$
Suppose that we have a morphism of analytic $S$-adic spaces
$f:X\to\D(z,\eps_0)$, and a complex $\cF^\bullet$ of
$\Lambda$-modules on $X_\et$; for every $n\in\N$, set
$$
X_n:=f^{-1}\D(z,\eps_n).
$$
In light of the discussion in \eqref{subsec_smaller-radius},
the foregoing choices determine a system of morphisms of
$\Lambda$-modules on $X_{0,\et}$ :
\set\begin{equation}\label{eq_vanish-cycles-at-0}
R^i\Psi_{\eta_z,\eps_n}\cF^\bullet_{|X_n}\to
R^i\Psi_{\eta_z,\eps_{n+1}}\cF^\bullet_{|X_{n+1}}
\qquad
\text{for every $i,n\in\N$}.
\end{equation}
For each $n\in\N$, the group $\pi(z,\eps_n)$ acts (on the left)
on the $n$-th term of this system, and the maps are
$\pi(z,\eps_{n+1})$-equivariant (for the action on the $n$-th
term obtained by restriction via \eqref{eq_restrict-radius}).
If we forget the group actions, we have for each $i\in\N$ a
direct system of $\Lambda$-modules on $X_{0,\et}$, indexed by
$n\in\N$, whose colimit shall be denoted simply
$$
R^i\Psi_{\eta_z}\cF^\bullet.
$$
With the notation of \eqref{subsec_vanish}, the unit of adjuction
$\one\to j_{H*}j^*_H$ induces, after taking the colimit over
the subgroups $H\subset\pi(z,\eps)$, a natural morphism :
$$
\cF^\bullet\to R\Psi_{\eta_z}\cF^\bullet
\qquad
\text{for every complex $\cF^\bullet$ of $\Lambda$-modules}
$$
whose cone shall be denoted :
$$
R\Phi_{\eta_z}\cF^\bullet.
$$

\begin{corollary}\label{cor_vanish-properties}
Assertions {\em (i)\/}, {\em (ii)\/} and {\em (iii)\/} of proposition
{\em\ref{prop_vanish-properties}} still hold with
$R\Psi_{\eta_z,\eps}$ replaced everywhere by either $R\Psi_{\eta_z}$
or $R\Phi_{\eta_z}$.
\end{corollary}
\begin{proof} To show assertion (i) for $R\Psi_{\eta_z}$, it
suffices to remark that, by general nonsense, the functors
$g^*$ and $g_0^*$ commute with arbitrary colimits; for assertion
(ii), likewise the functors and $Rg_*$ and $Rg_{0*}$ commute
with arbitrary filtered colimits (\cite[Lemma 2.3.13(ii)]{Hu2}).
Assertion (iii) is immediate. Finally, the assertions
concerning $R\Phi_{\eta_z}$ follow from the case of
$R\Psi_{\eta_z}$, since $g^*$ and $Rg_*$ are triangulated functors.
\end{proof}

Of course, the limit of the projective system
\eqref{eq_restrict-radius} still acts on
$R^i\Psi_{\eta_z}\cF^\bullet$; however, not much is known
concerning this projective limit. Nevertheless, in certain
cases it is possible to salvage a useful residual group action
on $R^0\Psi_{\eta_z}\cF^\bullet$ and $R^0\Phi_{\eta_z}\cF^\bullet$,
as explained in the following paragraph.

\sset\subsubsection{}\label{subsec_more-trivialities}
Namely, let $F$ be a $\Lambda$-module on $\D(z,\eps_0)_\et$,
and suppose that there exists $N\in\N$, such that the restriction
$F'$ of $F$ to $\D(z,\eps_N)^*_\et$ is a locally constant
$\Lambda$-module of finite type. In that case, $F'$ is determined,
up to natural isomorphism, by the stalk $F_{\xi_N}$ endowed with
its natural continuous $\pi(z,\eps_N)$-action. The latter is
obtained as follows. According to \eqref{subsec_was-remark},
$F'$ is trivialized by some locally algebraic Galois covering
$Y\to\D(z,\eps_N)^*$, corresponding to a certain open normal
subgroup $L\subset\pi(z,\eps_N)$. Then we have a natural
isomorphism $\Gamma(Y,F)\isom F_{\xi_N}$, and the quotient
group $\pi(z,\eps_N)/L$ acts on $\Gamma(Y,F)$ as explained
in \eqref{subsec_vanish}.

On the other hand, by inspecting the definitions, we see that
the $\Lambda$-module $R^0\Psi_{\eta_z}F[0]$ is computed as the
colimit
$$
\colim{n\in\N}\,\colim{H\subset\pi(z,\eps_n)}
\Gamma(C^{(n)}_{H,\et},F).
$$
The first colimit does not change if we replace $\N$ by the
cofinal subset of all integers $n\geq N$. Likewise, by
\eqref{subsec_was-remark}, we may moreover replace
the filtered system of all open subgroups $H\subset\pi(z,\eps_n)$
by the cofinal system of all such $H$ with the property that $F$
restricts to a constant sheaf on $C^{(n)}_{H,\et}$.
After these changes, notice that all the transition maps
in the filtered system are isomorphisms, since each $C^{(n)}_H$
is connected. There follows a natural isomorphism :
$$
R^0\Psi_{\eta_z}F[0]\isom F_{\xi_N}
$$
especially we obtain a natural continuous $\pi(z,\eps_N)$-action
on the above vanishing cycle.  It is then convenient to introduce
the notation
$$
F_{\eta_z}
$$
for the datum of $R^0\Psi_{\eta_z} F[0]$, viewed as an object
of $\Lambda[\pi(z,\eps_N)]\Mod_\mathrm{f.cont}$. Finally, to
eliminate the dependence on the choice of $N$, notice that
the group homomorphisms \eqref{eq_restrict-radius} induce a
direct system of categories and functors :
\set\begin{equation}\label{eq_dir-syst-cats}
\Lambda[\pi(z,\eps_n)]\Mod_\mathrm{f.cont}\to
\Lambda[\pi(z,\eps_{n+1})]\Mod_\mathrm{f.cont}
\qquad
\text{for every $n\in\N$.}
\end{equation}
We may therefore consider the 2-colimit
$$
\Lambda[\pi(z)]\Mod:=
\Pscolim{n\in\N}\Lambda[\pi(z,\eps_n)]\Mod_\mathrm{f.cont}.
$$
Then, under the foregoing assumptions on $F$, we may regard
$F_{\eta_z}$ as an object of this category, which should be
thought of as a replacement for the category of continuous
representations of the inertia subgroup which one encounters
in the study of the local monodromy (around a given closed
point) of \'etale sheaves on algebraic curves.

Likewise, we may represent $R^0\Phi_{\eta_z}F[0]$ as the
cokernel of the natural map :
$$
F_z\to F_{\eta_z}
$$
which is a morphism in $\Lambda[\pi(z)]\Mod$, provided we
endow $F_z$ with the trivial $\pi(z,\eps_0)$-action.
These mostly trivial observations shall be amplified later,
starting with remark \eqref{rem_not-much}(i).

\begin{remark}\label{rem_one-more-obviety}
Following standard conventions, we shall use the
notation $F_{\eta_z}$ also for $\Lambda$-modules defined
only on $\D(z,\eps)^*_\et$ : rigorously, this means that we
replace $F$ by, for instance, its extension by zero over the
whole of $\D(z,\eps)$.
\end{remark}

\begin{proposition}\label{prop_vanish-vanishes}
In the situation of \eqref{subsec_pass-to-lim}, suppose that
$X\to\D(z,1)$ is a smooth morphism, and let $M$ be a
$\Lambda$-module of finite length. Then we have :
$$
R^i\Psi_{\eta_z}M_X=
\left\{\begin{array}{ll}
              M_{X_0} & \text{if\/ $i=0$} \\
              0       & \text{otherwise}.
\end{array}\right.
$$
\end{proposition}
\begin{proof} In light of corollary \ref{cor_vanish-properties},
we may assume that $X=\D(z,1)$, in which case, notice that the
open subspaces $\D(z,\eps)\subset\D(z,1)$ form a cofinal system
of \'etale neighborhoods of $z\in\D(z,1)$.
Since the stalk of a presheaf is isomorphic to the stalk of the
associated sheaf, we are thus reduced to showing that the
$\Lambda$-module
$$
\colim{n\in\N}\,\colim{H\subset\pi(z,\eps_n)}\,H^i(C^{(n)}_{H,\et},M)
$$
equals $M$ for $i=0$, and vanishes for $i>0$. The assertion for
$i=0$ just translates the fact that $C^{(n)}_H$ is a connected
adic space, for every $n\in\N$ and every open subgroup
$H\subset\pi(z,\eps_n)$. The case where $i=1$ follows from the
more precise lemma \ref{lem_complete-later}.

The assertion for $i\geq 2$ can be deduced from
\cite[Cor.7.5.6]{Hu2} by a reduction to the case where
$\Lambda$ is a finite ring.
However, it is quicker to reduce first to the case where $\Lambda$
has finite length as a $\Lambda$-module, then to the case where
$\Lambda$ is a field (by an easy induction on the length of
$\Lambda$), and then apply the Poincar\'e duality of lemma
\ref{lem_PD-lambda}(iii) : the details shall be left to the reader.
\end{proof}

\begin{corollary}\label{cor_vanish-smooth}
Suppose that $X\to\D(z,1)$ is a smooth morphism, and let $F$
be a locally constant $\Lambda$-module of finite type on
$X_\et$. Then :
$$
R\Phi_{\eta_z}F=0
\qquad
\text{in $\sD^+(X_0,\Lambda)$.}
$$
\end{corollary}
\begin{proof} By \eqref{subsec_was-remark}, $F$ is trivialized on
a locally algebraic covering $g:Y\to X$. Let $g_0:g^{-1}X_0\to X_0$
be the restriction of $g$; clearly it suffices to show that
$g^*_0R\Phi_{\eta_z}F$ vanishes. The latter is an immediate
consequence of propositions \ref{prop_vanish-properties}(i)
and \ref{prop_vanish-vanishes}, and corollary \ref{cor_vanish-properties}.
\end{proof}

\subsection{Fourier transform}\label{sec_Fourier}
We wish to recall briefly the construction and the main properties
of the cohomological Fourier transform introduced in \cite{Ram0}.
To begin with, we remark that the logarithm power series :
$$
\log(1+x):=\sum_{n\in\N}(-1)^n\cdot\frac{x^{n+1}}{n+1}
$$
defines a locally algebraic covering
$$
\log:\D(0,1^-)\to(\A^1_K)^\ad
\qquad
\text{where}\quad \D(0,1^-):=\bigcup_{0<\eps<1}\D(0,\eps).
$$
Moreover, $\log$ is a morphism of {\em analytic adic groups\/},
if we endow $(\A^1_K)^\ad$ with the additive group law of
$(\G_a)^\ad$, and $\D(0,1^-)$ with the restriction of the multiplicative
group law of $(\G_m)^\ad$ (more precisely, with the group law
given by the rule : $(x,y)\mapsto x+y+xy$, which is the ``translation''
of the usual multiplication, whose neutral element is $0\in\D(0,1^-)$).

The kernel of $\log$ is the group
$$
\mu_{p^\infty}:=\bigcup_{n\in\N}\mu_{p^n}
$$
of $p$-primary roots of of $1$ in $K$; therefore $\D(0,1^-)$ represents
a $\mu_{p^\infty}$-torsor on $(\A^1_K)^\ad_\et$, which we denote $\cL$.

\sset\subsubsection{}\label{subsec_Lubin-Tate-torsor}
Suppose now that the multiplicative group $\Lambda^\times$ contains
a subgroup isomorphic to $\mu_{p^\infty}$ (for instance, this holds
for $\Lambda:=\bar\F_\ell$, if $\ell$ is any prime number different
from $p$). Then, the choice of a non-trivial group homomorphism :
$$
\psi:\mu_{p^\infty}\to\Lambda^\times
$$
determines a locally constant $\Lambda$-module, free of rank one on
$(\A^1_K)^\ad_\et$ :
$$
\cL_\psi:=\cL\overset{\mu_{p^\infty}}{\times}\Lambda
$$
whose restriction to $\D(0,1^-)_\et$ is a constant $\Lambda$-module.

Let $X$ be any $S$-adic space, and $f\in\cO_X(X)$ be any global
section; we may regard $f$ as a morphism of $S$-adic spaces
$f:X\to(\A^1_K)^\ad$, and then we may set :
$$
\cL_\psi\La f\Ra:=f^*\cL_\psi
$$
which is a locally constant $\Lambda$-module on $X_\et$. If
$f,g\in\Gamma(X,\cO_X)$ are any two analytic functions, we have
a natural isomorphism of $\Lambda$-modules :
\set\begin{equation}\label{eq_sum-for-LT}
\cL_\psi\La f+g\Ra\isom\cL_\psi\La f\Ra\otimes_\Lambda\cL_\psi\La g\Ra.
\end{equation}

\sset\subsubsection{}\label{subsec_intro-Fourier}
Let $A$ and $A'$ be two copies of the affine line $(\A^1_K)^\ad$,
with global coordinates $x$ and respectively $x'$.
Let also $D$ and $D'$ be two copies of the projective line
$(\P^1_K)^\ad$, and fix open immersions
$$
\alpha:A\to D
\qquad
\alpha':A'\to D'
$$
so that $D=A\cup\{\infty\}$ and $D'=A'\cup\{\infty'\}$.
We use $\cL_\psi$ as the ``integral kernel'' of our operator; namely,
the {\em Fourier transform\/} is the triangulated functor
$$
\cF_\psi:\sD^b(A_\et,\Lambda)\to\sD^b(A'_\et,\Lambda)
\qquad
F^\bullet\mapsto Rp'_!(p^*F^\bullet\otimes_\Lambda\cL_\psi\La m\Ra)[1]
$$
where
$$
A\xleftarrow{\ p\ }A\times_SA'\xrightarrow{\ p'\ }A'
$$
are the two projections, and $m:A\times_SA'\to(\A^1_K)^\ad$
is the dual pairing, given by the rule : $(x,x')\mapsto x\cdot x'$.
Especially, let $0'\in A'$ be the origin of the coordinate $x'$
({\em i.e.} $x'(0')=0$ in the residue field of the point $0'$);
then notice the natural isomorphism :
\set\begin{equation}\label{eq_stalk-at-0-dual}
\cF(F^\bullet)_{0'}\isom R\Gamma_{\!c}\,F^\bullet
\qquad
\text{for every $F^\bullet\in\Ob(\sD^b(A_\et,\Lambda))$.}
\end{equation}
It is shown in \cite[Th.7.1.2]{Ram0} that $\cF_\psi$ is
an equivalence of categories; more precisely, let $A''$ be a third
copy of the affine line, with global coordinate $x''$, and dual
pairing
$$
m':A'\times_SA''\to(\A^1_K)^\ad
\qquad
(x',x'')\mapsto x'\cdot x''.
$$
We have a Fourier transform
$\cF'_\psi:\sD^b(A'_\et,\Lambda)\to\sD^b(A''_\et,\Lambda)$ (namely,
the ``integral operator'' whose kernel is $\cL_\psi\La m'\Ra$),
and there is a natural isomorphism of functors :
\set\begin{equation}\label{eq_FT-is-selfdual}
\cF'_\psi\circ\cF_\psi F^\bullet\isom a_* F^\bullet(-1)
\end{equation}
where $(-1)$ denotes Tate twist and $a:A\isom A''$ is $(-1)$-times
the double duality isomorphism, given by the rule : $x\mapsto -x''$.
An immediate consequence of \eqref{eq_stalk-at-0-dual} and
\eqref{eq_FT-is-selfdual} is the natural isomorphism :
\set\begin{equation}\label{eq_stalk-at-0}
F^\bullet_0(-1)[-1]\isom R\Gamma_{\!c}\,\cF_\psi(F^\bullet)
\qquad
\text{for every $F^\bullet\in\Ob(\sD^b(A_\et,\Lambda))$}
\end{equation}
where $F^\bullet_0$ denotes the stalk of $F^\bullet$ over the origin
$0\in A$ of the coordinate $x$.

\sset\subsubsection{}
Moreover, $\cF$ commutes with Poincar\'e-Verdier duality; the latter
assertion is a consequence of the following. Define a second functor :
$$
\cF_{\psi,*}:\sD^b(A_\et,\Lambda)\to\sD^b(A'_\et,\Lambda)
\qquad
F^\bullet\mapsto Rp'_*(p^*F^\bullet\otimes_\Lambda\cL_\psi\La m\Ra)[1].
$$

\begin{proposition}\label{prop_perverse-FT}
The natural transformation
$$
\cF_\psi\to\cF_{\psi,*}
$$
(deduced from the natural morphism $Rp'_!\to Rp'_*$) is an isomorphism
of functors.
\end{proposition}
\begin{proof} This is shown in \cite[Th.7.1.6]{Ram0}. A ``formal''
proof can also be given as in a recent preprint by Boyarchenko and
Drinfeld (\cite[Th.G.6]{Bo-Dr}) : in order to repeat the argument
in our context, one needs a version of Poincar\'e duality for
$\Lambda$-modules. We sketch the arguments in the following paragraph.
\end{proof}

\sset\subsubsection{}\label{subsec_Poincare-dual}
To begin with, let $R$ be any torsion ring, and $f:X\to Y$ a morphism
of analytic pseudo-adic spaces which is locally of $^+$weakly finite
type, separated, taut, and of finite transcendental dimension; by
\cite[Th.7.1.1]{Hu2}, the functor
$$
Rf_!:\sD^+(X_\et,R)\to\sD^+(Y_\et,R)
$$
admits a right adjoint
$$
Rf^!_{(R)}\sD^+(Y_\et,R)\to\sD^+(X_\et,R).
$$
On the other hand, for any ring homomorphism $R\to R'$, denote
by $\rho_X:\sD^+(X_\et,R')\to\sD^+(X_\et,R)$
the restriction of scalars, and likewise define $\rho_Y$; then
we have :

\begin{lemma}\label{lem_PD-lambda}
With the notation of \eqref{subsec_Poincare-dual}, the following
holds :
\begin{enumerate}
\item
There exists a natural isomorphism of functors :
$$
\rho_X\circ Rf^!_{(R')}\isom Rf^!_{(R)}\circ\rho_Y.
$$
\item
Suppose additionally that $f$ is smooth of pure dimension $d$,
and that $R$ is a $\Z/n\Z$-algebra, where $n\in\N$ is invertible
in $\cO^+_Y$; then there is a natural isomorphism :
$$
f^*F^\bullet(d)[2d]\isom Rf_{(R)}^!F^\bullet
\qquad
\text{for every $F^\bullet\in\Ob(\sD^+(Y_\et,\Lambda))$.}
$$
\item
In the situation of\/ {\em (ii)}, assume that $Y=S$, and that $R$
is an injective $R$-module ({\em e.g.} $R$ is a field). Let $F$
be a constructible locally constant $R$-module on $X_\et$;
then there are natural isomorphisms :
$$
H^{2d-q}(X_\et,\cHom_R(F,R_X)(d))\isom\Hom_R(H^q_c(X_\et,F),R)
\qquad
\text{for every $q\in\Z$}.
$$
\end{enumerate}
\end{lemma}
\begin{proof} (i): This isomorphism is not explicitly stated in
\cite{Hu2}, but it can be proved as in
\cite[Exp.XVIII, Prop.3.1.12.1]{SGA4-3}.

(ii): Such isomorphism is established in \cite[Th.7.5.3]{Hu2}, for
$R=\Z/n\Z$; the general case follows from this and from (i).

(iii) is a standard consequence of (ii) : cp. \cite[Cor.7.5.6]{Hu2}.
\end{proof}

\sset\subsubsection{}\label{subsec_assumption-Lambda}
Clearly our $\Lambda$ is a $\Z/n\Z$-algebra, for a suitable
$n$ fulfilling the conditions of lemma \ref{lem_PD-lambda}(ii)
(take $n>0$ such that $n\cdot 1=0$ in $\Lambda$). On the other
hand, it is not clear to me whether, under our current assumptions,
$\Lambda$ is always an injective $\Lambda$-module. In order
to apply lemma \ref{lem_PD-lambda}(iii) to a given $\Lambda$-module
$F$, we shall therefore try to find a finite descending filtration
$\Fil^\bullet F$ of $F$, whose graded terms $\gr^i F$ are
annihilated by residue fields of $\Lambda$. This is always
possible, provided $\Lambda$ is noetherian (since in this case
$\Lambda$ shall be also artinian). For this reason,
{\em henceforth we assume that $\Lambda$ is a noetherian
ind-finite ring such that $p\in\Lambda^\times$}.

\sset\subsubsection{}\label{subsec_especially-inter}
Let $X$ be any scheme locally of finite type over $\Spec\,K$,
and $F$ a $\Lambda$-module on the \'etale site of $X^\ad$; recall
that $F$ is {\em Zariski constructible}, if for every $x\in X^\ad$
there exists a subset $Z\subset X$ locally closed and constructible
for the Zariski topology of $X$, such that $x\in Z^\ad$, and such
that the restriction of $F$ is a locally constant $\Lambda$-module
of finite type on the \'etale site of $Z^\ad$.
We are especially interested in the Fourier transform of Zariski
constructible $\Lambda$-modules with {\em bounded ramification\/}
on $A_\et$. To explain the latter condition, consider more generally,
for a given $K$-rational point $z$ of $\P^1_K$, a locally constant
$\Lambda$-module $F$ of finite type on $\D(z,\eps)^*$; in
\cite[Cor.4.1.16]{Ram} we attach to $F$ its {\em Swan conductor\/}
$$
\ssw^\natural_z(F,0^+)\in\N\cup\{+\infty\}
$$
and following \cite[Def.4.2.1]{Ram}, we say that $F$ has {\em bounded
ramification around $z$}, if $\ssw^\natural_z(F,0^+)$ is a (finite)
integer. Regrettably, the definition in {\em loc.cit.} assumes
moreover that the stalks of $F$ are free $\Lambda$-modules, but
by inspecting the proofs, one sees easily that all the results of
\cite[\S4.2]{Ram} hold when the stalks are only assumed to be
$\Lambda$-modules of finite length : the point is that the Swan
conductor depends only on the class of $F$ in the Grothendieck
group of $\Lambda$-modules on $D(z,\eps)^*$, hence we may reduce
to the case where $F$ is annihilated by a maximal ideal of $\Lambda$,
by considering a descending filtration
$F\supset F_1\supset F_2\supset\cdots$, where $F_{i+1}=\fn_i F_i$
for every $i\geq 1$, with $\fn_i\subset\Lambda$ a maximal ideal.
(Such a filtration is compatible with the break decompositions
of {\em loc.cit.})

\sset\subsubsection{}\label{subsec_tame-ram}
For instance, suppose that $F$ has {\em finite monodromy},
{\em i.e.} $F$ is trivialized on a finite \'etale covering
$\D(z,\eps)^*$; in that case $\ssw^\natural_z(F,0^+)=0$, especially
$F$ has bounded ramification. In such situation, the $p$-adic
Riemann existence theorem (\cite[Th.2.4.3]{Ram}) ensures that
-- after replacing $\eps$ by a smaller radius -- $F$ is trivialized
on a covering {\em of Kummer type} :
$$
\D(z,\eps^{1/d})^*\to\D(z,\eps)^* \quad :\quad x\mapsto x^d.
$$
For this reason, we shall also say that such $F$ is {\em tamely
ramified}.

\sset\subsubsection{}\label{subsec_Groth-Ogg-Sha}
Let $C$ be a smooth connected projective curve defined over $\Spec\,K$,
and $F$ a $\Lambda$-module on the \'etale site of $C^\ad$; assume
that $F$ is Zariski-constructible, so all the stalks of $F$ are
$\Lambda$-modules of finite length, and there exists an open subset
$U\subset C$ such that $F$ restricts to a locally constant
$\Lambda$-module on the \'etale site of $U^\ad$. For every
$x\in C\setminus U$, we may choose an open immersion
$\phi_x:\D(0,1)\to C^\ad$ of $S$-adic spaces, such that
$\phi(0)=x$. Then we may define the {\em Swan conductor of $F$ at $x$}
as the integer
$$
\ssw_x^\natural(F,0^+):=\ssw_0^\natural(\phi^*_xF,0^+).
$$
One verifies easily that this quantity is independent of the choice
of $\phi_x$, hence we say that {\em $F$ has bounded ramification
around $x$}, if $\ssw_x^\natural(F,0^+)\in\N$. One may calculate
the Euler-Poincar\'e charateristic of $F$, {\em i.e.} the integer
$$
\chi_c(C,F):=\sum_{i=0}^2 (-1)^i\cdot\length_\Lambda H^i_c(C_\et,F)
$$
in terms of these local invariants. Indeed, denote by $l(F)$ the
{\em generic length\/} of $F$, defined as $\length_\Lambda(F_\xi)$
for any (hence all) geometric point $\xi$ of $U^\ad$; moreover,
for every $x\in C\setminus U$, let :
$$
a_x(F):=l(F)+\ssw_x^\natural(F,0^+)-\length_\Lambda F_{\bar x}
$$
where $\bar x$ is any geometric point of $C^\ad$ localized at $x$;
this is the {\em total drop of the rank\/} of $F$ at the point $x$.
With this notation, we have the following
{\em Grothendieck-Ogg-Shafarevich formula} :

\begin{lemma}\label{lem_Groth-Ogg-Shaf}
In the situation of \eqref{subsec_Groth-Ogg-Sha}, suppose that $F$
has bounded ramification at all the points of $C\setminus U$, and
that $\Lambda$ is a local ring with residue field $\bar\Lambda$. Then :
$$
\chi_c(C,F)=
\chi_c(C,\bar\Lambda_C)\cdot l(F)-\sum_{x\in C\setminus U} a_x(F).
$$
\end{lemma}
\begin{proof} We may easily reduce to the case where
$F=j_!(F_{|U^\ad})$, where $j:U^\ad\to C^\ad$ is the open immersion.
Moreover, using the additivity of the Swan conductor and of the
Euler-Poincar\'e characteristic, we may assume that $\Lambda$
is a field (for a general $\Lambda$, one considers the
$\fm_\Lambda$-filtration of $F$, where $\fm_\Lambda\subset\Lambda$
is the maximal ideal). Then we may argue as in
\cite[Rem.10.6]{Hu3} : the details shall be left to the reader.
\end{proof}

\sset\subsubsection{}\label{subsec_brk-decomp}
In the situation of \eqref{subsec_especially-inter}, suppose that
$F$ has bounded ramification around $z$; then \cite[Th.4.2.40]{Ram}
ensures the existence of a connected open subset $U\subset\D(z,\eps)^*$
such that $U\cap\D(z,\delta)^*\neq\emptyset$ for every
$\delta\in|K^\times|$, and a canonical {\em break decomposition}
of $\Lambda_U$-modules :
$$
F_{|U}=\bigoplus_{\gamma\in\Gamma}F(\gamma)
$$
indexed by the group $\Gamma:=\Q\times|K^\times|$, which we endow
with the lexicographic ordering given by the rule :
$$
(q,c)\leq(q',c')
\quad \Leftrightarrow \quad
\text{either $q<q'$ or else $q=q'$ and $c\geq c'$}
$$
(there is a slight confusion in \cite[\S4.2.39]{Ram} concerning
this ordering). The projection :
$$
\Gamma\to\Q \qquad (q,c)\mapsto(q,c)^\natural:=q
$$
is a homomorphism of ordered groups; we have $F(\gamma)=0$ unless
$\gamma^\natural\geq 0$. The break decomposition is functorial : if
$\phi:F\to G$ is a morphism of $\Lambda$-modules with bounded
ramification, and $U$ as above is chosen small enough so that
both $F_{|U}$ and $G_{|U}$ admit break decompositions, then $\phi_{|U}$
restricts to morphisms of $\Lambda_U$-modules $F(\gamma)\to G(\gamma)$
for every $\gamma\in\Gamma$.

Moreover, the break decomposition is well behaved with respect
to tensor products : namely, if $F$ and $G$ have both bounded
ramification, then the same holds for $F\otimes_\Lambda G$, and
for a sufficiently small $U$ as above, and every
$\gamma,\gamma'\in\Gamma$, we have :
\set\begin{equation}\label{eq_tannakian-property}
F(\gamma)\otimes_\Lambda G(\gamma')\subset\left\{
\begin{array}{ll}
(F\otimes_\Lambda G)(\max(\gamma,\gamma')) &
\text{if $\gamma\neq\gamma'$} \\
\bigoplus_{\rho\leq\gamma}(F\otimes_\Lambda G)(\rho) &
\text{otherwise.}
\end{array}\right.
\end{equation}
Also, the break decomposition of $\Hom_\Lambda(F,G)$ is likewise
related to those of $F$ and $G$.

Furthermore, if $F$ is a locally free $\Lambda$-module,
$\gamma^\natural\cdot\rk_\Lambda F(\gamma)\in\N$ for every
$\gamma\in\Gamma$ (Hasse-Arf property); when $F$ is only
constructible, we have the identity :
\set\begin{equation}\label{eq_break-and-swan}
\ssw^\natural_z(F,0^+)=
\sum_{\gamma\in\Gamma}\gamma^\natural\cdot\length_\Lambda F(\gamma).
\end{equation}
The elements $\gamma\in\Gamma$ such that $F(\gamma)\neq 0$ are
called the {\em breaks\/} of $F$ around $z$, and the integer
$\length_\Lambda F(\gamma)$ is called the {\em multiplicity}
of $\gamma$. In view of \eqref{eq_break-and-swan}, it is
convenient to define
$$
F^\natural(r):=\bigoplus_{\gamma^\natural=r}F(\gamma)
\qquad
\text{for every $r\in\Q_+$}.
$$
The $r\in\Q_+$ such that $F^\natural(r)\neq 0$ are called the
$\natural$-breaks of $F$.

\sset\subsubsection{}\label{subsec_drop}
For instance, the $\Lambda$-module $\cL_\psi$ of
\eqref{subsec_Lubin-Tate-torsor} has bounded ramification at
the point $\infty\in D$. Since $\cL_\psi$ has $\Lambda$-rank one,
obviously it admits a single break $\gamma$ around $\infty$. More
generally, let $f\in\Gamma(\A^1_K,\cO_{\A^1_K})$ be an algebraic
function, and denote by $o_\infty(f)\in\N$ the pole order of $f$
at $\infty$; then $\cL_\psi\La f\Ra$ has bounded ramification
around $\infty$, with a unique break of the form
$$
(o_\infty(f),c)\in\N\times|K^\times|\subset\Gamma
$$
where the value $c$ depends only on the valuation of the leading
coefficient of $f$, and on the largest power of $p$ dividing
$o_\infty(f)$. Indeed, this follows from
\cite[Lemma 4.2.6(i),(ii)]{Ram} if $f$ is a monomial
$f(T)=a\cdot T^n$ (for a chosen coordinate $T$, so that
$\A^1_K=\Spec\,K[T]$), and the general case is easily deduced, by
writing $f=f_n+g$, where $f_n$ is a monomial and $\deg_Tg<n:=\deg_Tf$,
and by applying induction on $n$, together with \eqref{eq_sum-for-LT}
and \eqref{eq_tannakian-property}.

Now, let $F$ be any $\Lambda$-module on $\D(\infty,\eps)^*$, with
bounded ramification around $\infty$; set 
$$
e(F):=\sum_{r\geq 1}r\cdot\length_\Lambda F^\natural(r)+
\sum_{r<1}\length_\Lambda F^\natural(r).
$$
For any $a\in K$, let
$F_{[a]}:=F\otimes_\Lambda\cL_\psi\La m_a\Ra_{|\D(\infty,\eps)^*}$,
where $m_a(T):=aT\in\Gamma(\A^1_K,\cO_{\A^1_K})$; taking
\eqref{eq_tannakian-property} and \eqref{eq_break-and-swan}
into account, we see that :
$$
\ssw^\natural_\infty(F_{[a]},0^+)\leq e(F)
\qquad
\text{for every $a\in K$.}
$$
We consider the subset
$$
d_K(F):=\{a\in K~|~\ssw^\natural_\infty(F_{[a]},0^+)<e(F)\}.
$$
We shall regard $d_K(F)$ as a subset of $A'(K)$, via the identification
$K=A'(K)$ induced by the global coordinate $x'$.

\begin{remark}\label{rem_when-0-drops}
Notice that $a\in d_K(F)$ if and only if $F_{[a]}^\natural(r)\neq 0$
for some $r<1$ (on a suitable $U\subset\D(\infty,\eps)^*$ as in
\eqref{subsec_brk-decomp}). Especially, $0\in d_K(F)$ if and only
if $F^\natural(r)\neq 0$ for some $r<1$.
\end{remark}

\begin{proposition}\label{prop_drop-is-finite}
With the notation of \eqref{subsec_drop}, the following holds :
\begin{enumerate}
\item
$d_K(F)$ is a finite set.
\item
There exist :
\begin{enumerate}
\item
a connected open subset $U\subset\D(\infty,\eps)^*$,
such that $U\cap\D(\infty,\delta)^*\neq\emptyset$ for every
$\delta\in|K^\times|$
\item
and a natural decomposition of $\Lambda_U$-modules :
$$
F_{|U}=\left(\bigoplus_{r>1}F^\natural(r)\right)
\oplus\left(\bigoplus_{a\in d_K(F)}\bigoplus_{r<1}
F^\natural_{[a]}(r)\otimes_\Lambda\cL_\psi\La m_{-a}\Ra_{|U}\right)
\oplus G
$$
where 
$$
G:=F^\natural(1)\cap\bigcap_{a\in K}
(F^\natural_{[a]}(1)\otimes_\Lambda\cL_\psi\La m_{-a}\Ra_{|U}).
$$
\end{enumerate}
\item
If\/ $0\to F_1\to F\to F_2\to 0$ is any short exact sequence of
locally constant $\Lambda$-modules on $\D(\infty,\eps)^*_\et$, we have :
$$
d_K(F)=d_K(F_1)\cup d_K(F_2).
$$
\item
If $F$ is a locally free $\Lambda$-module of finite type on
$\D(\infty,\eps)^*_\et$, then :
$$
d_K(F^\vee)=\{-a~|~a\in d_K(F)\}
\qquad
\text{where $F^\vee:=\cHom_\Lambda(F,\Lambda_{\D(\infty,\eps)})$}.
$$
\end{enumerate}
\end{proposition}
\begin{proof} For any finite subset $\Sigma\subset K$, let
$$
G_\Sigma:=F^\natural(1)\cap
\bigcap_{a\in\Sigma}(F^\natural_{[a]}(1)\otimes_\Lambda
\cL_\psi\La m_{-a}\Ra_{|U})
$$
where the intersection is taken over a sufficiently small
open subset $U$, fulfilling condition (a).

\begin{claim}\label{cl_set-S}
For every $b\in K\setminus\Sigma$, we have (on a suitable
open subset $U$) :
$$
G_\Sigma=\left(
\bigoplus_{r<1}
F^\natural_{[b]}(r)\otimes_\Lambda\cL_\psi\La m_{-b}\Ra_{|U}\right)
\oplus G_{\Sigma\cup\{b\}}.
$$
\end{claim}
\begin{pfclaim} On the one hand, we may write :
\set\begin{equation}\label{eq_first}
F^\natural(1)=\left(
\bigoplus_{r<1}
F^\natural_{[b]}(r)\otimes_\Lambda\cL_\psi\La m_{-b}\Ra_{|U}\right)
\oplus
\left((F^\natural_{[b]}(1)\otimes_\Lambda\cL_\psi\La m_{-b}\Ra_{|U})
\cap F^\natural(1)\right).
\end{equation}
On the other hand, notice that :
\set\begin{equation}\label{eq_second}
F^\natural_{[b]}(r)\otimes_\Lambda\cL_\psi\La m_{-b}\Ra_{|U}
\subset G_\Sigma
\qquad
\text{for every $r<1$}.
\end{equation}
Indeed, \eqref{eq_second} is equivalent to the family of inclusions :
$$
F^\natural_{[b]}(r)\otimes_\Lambda\cL_\psi\La m_{a-b}\Ra_{|U}
\subset F_{[a]}^\natural(1)
\qquad
\text{for every $a\in\Sigma$ and every $r<1$}
$$
which hold, since $a-b\neq 0$ for every $a\in\Sigma$. The claim
follows straightforwardly from \eqref{eq_first} and \eqref{eq_second}.
\end{pfclaim}

Using claim \ref{cl_set-S} we deduce, by induction on the cardinality
of $\Sigma$, the decomposition :
$$
F_{|U}=\left(\bigoplus_{r>1}F^\natural(r)\right)\oplus
\left(\bigoplus_{a\in\Sigma}\bigoplus_{r<1}
F^\natural_{[a]}(r)\otimes_\Lambda\cL_\psi\La m_{-a}\Ra_{|U}\right)
\oplus G_\Sigma
$$
for every finite subset $\Sigma\subset K$ such that $0\in\Sigma$
(on a suitable $U$, dependent on $\Sigma$; see remark
\ref{rem_when-0-drops}). However, it is easily seen that $a\in d_K(F)$
if and only if there exists $r<1$ such that $F_{[a]}^\natural(r)\neq 0$.
Since $F$ is a locally constant $\Lambda$-module of finite type,
assertion (i) is an immediate consequence, and we also see that
$G=G_{d_K(F)}$.

Finally, (iii) is an immediate consequence of the additivity and
non-negativity properties of the Swan conductor and (iv) is
straigthforward from \eqref{eq_sum-for-LT}.
\end{proof}

\sset\subsubsection{}
Let now $K\subset L$ be an extension of algebraically closed
non-archimedean fields with rank $1$ valuations; denote by
$$
\pi_L:\D(\infty,\eps)\times_S\Spa(L,L^+)\to\D(\infty,\eps)
$$
the projection. In view of \cite[Lemma 3.3.8]{Ram} we have
$\ssw^\natural_\infty(\pi_L^*F,0^+)=
\ssw^\natural_\infty(F,0^+)$, hence $\pi^*_LF$
has bounded ramification, and we may consider the set
$d_L(\pi_L^*F)$. With this notation we have :

\begin{corollary}\label{cor_drop-is-finite}
$d_L(\pi^*_LF)=d_K(F)$.
\end{corollary}
\begin{proof} Let $A'_L:=A'\times_S\Spa(L,L^+)$ and denote by
$\pi'_L:A'_L\to A'$ the projection; under the natural
identifications $K=A'(K)$ and $L=A'_L(L)$, the inclusion
$K\subset L$ corresponds to the map 
$$
A'(K)\to A'_L(L)
\quad :\quad
x\mapsto\pi_L^{\prime-1}(x)
$$
and the assertion means precisely that
$d_L(\pi^*_LF)=\pi_L^{\prime-1}d_K(F)$. Taking into account
proposition \ref{prop_drop-is-finite}(i), it then suffices to
show that, for every $a\in L\setminus K$, there exist an infinite
subset $\Sigma\subset K$ such that
\set\begin{equation}\label{eq_swan-sigma}
\ssw^\natural_\infty((\pi_L^*F)_{[a+x]})=
\ssw^\natural_\infty((\pi_L^*F)_{[a]})
\qquad
\text{for every $x\in\Sigma$.}
\end{equation}
To this aim, we may assume that $L$ is the completion of the
algebraic closure of $K(a)$ (\cite[Lemma 3.3.8]{Ram}); in this
case, we will exhibit, more precisely, an infinite subset
$\Sigma\subset K$, and for every $x\in\Sigma$, an isometric
$K$-automorphism $\sigma_x:L\to L$ such that $\sigma_x(a)=a+x$.
Then $\sigma_x$ will induce an automorphism $\sigma'_x$ of $A'_L$
such that
$$
\pi'_L\circ\sigma'_x=\pi'_L
\qquad
\text{and}
\qquad
\sigma^{\prime*}_x\cL_\psi\La m_a\Ra\simeq\cL_\psi\La m_{a+x}\Ra
$$
from which \eqref{eq_swan-sigma} follows straightforwardly.

\begin{claim}\label{cl_extend-isometry}
Let $(E,|\cdot|_E)$ be any complete non-archimedean valued field,
$(E^\mathrm{a},|\cdot|_{E^\mathrm{a}})$ a fixed algebraic closure of
$E$, and $\sigma$ an isometric automorphism of $E$.
Then $\sigma$ extends to an isometric automorphism of $E^\mathrm{a}$.
\end{claim}
\begin{pfclaim} This is presumably well known : by Zorn's lemma,
we may find a maximal subextension $F\subset E^\mathrm{a}$ such that
$\sigma$ extends to an isometry $\sigma_F:F\isom F'$ (where
$F'\subset E^\mathrm{a}$ as well). We claim that $F=E^\mathrm{a}$;
indeed, if $x\in E^\mathrm{a}\!\setminus\!F$, let $P$ be the minimal
polynomial of $x$ over $F$, $P^\sigma$ the polynomial obtained
by applying $\sigma_F$ to all the coefficients of $P$, and $x'$
any root of $P^\sigma$; we may extend $\sigma_F$ to an isomorphism
$\tilde\sigma_F:F[x]\to F[x']$ by mapping $x$ to $x'$. Notice now
that the spectral norms of $F[x]$ and $F[x']$ agree with the
restriction of $|\cdot|_{E^\mathrm{a}}$ (\cite[Lemma 1.1.17]{Ram}),
and on the other hand, $\tilde\sigma_F$ clearly preserves these
spectral norms; therefore $\tilde\sigma_F$ is an isometry
$(F[x],|\cdot|_{E^\mathrm{a}})\to(F[x'],|\cdot|_{E^\mathrm{a}})$.
By maximality, we must have $F[x]=F$, a contradiction.
We leave to the reader the verification that
$\sigma^\mathrm{a}(E^\mathrm{a})=E^\mathrm{a}$.
\end{pfclaim}

Now, the norm $|\cdot|_L$ of $L$ corresponds to a maximal point of
$\Spec\,K[a]\times_{\Spec\,K}S$. According to the classification
in \cite[\S1.4.4]{Berk}, such point can only be of type (2), (3)
or (4). We will show show that, in either of these cases, there
exists $\rho_a\in|K^\times|$ such that the $K$-automorphism $\omega_x$
of $K[a]$ given by the rule $a\mapsto a+x$ is an isometry for the
norm $|\cdot|_L$, whenever $x\in K$ and $|x|<\rho_a$. In view
of claim \ref{cl_extend-isometry}, this $\omega_x$ will extend to
an automorphism $\sigma_x$ of $L$ as sought.

Indeed, points of type (2) and (3) arise from the sup norm on a
certain disc $\D(x_0,\rho)$ : in these cases, it is easily seen
that $\rho_a:=\rho$ will do. Likewise, a point of type (4) arises
from a family $\{\D(x_i,\rho_i)~|~i\in I\}$ of embedded discs;
the infimum of the radii $\rho_i$ must be $>0$, otherwise $a$
would be of type (1); then we may take $\rho_a$ equal to this
infimum.
\end{proof}

\sset\subsubsection{}\label{subsec_set-up_-for-th}
Let $U\subset A$ be a Zariski open subset ({\em i.e.} $U=V^\ad$
for some open subset $V$ of the scheme $\A^1_K$), $F$ a
constructible and locally constant $\Lambda$-module on $U_\et$,
and $j:U\to A$ the open immersion. Set $F':=\cF_\psi(j_*F)$
and $U':=A'\setminus d_K(F)$.

\begin{theorem}\label{th_perversity}
In the situation of \eqref{subsec_set-up_-for-th}, suppose
moreover that $F$ has bounded ramification at all the
(finitely many) points of $D\setminus U$. Then :
\begin{enumerate}
\item
$\cH^1(F')$ is a constructible $\Lambda$-module, and its support
is contained in $d_K(F)$.
\item
$\cH^0(F')$ restricts to a locally constant $\Lambda$-module on
$U'_\et$, and has bounded ramification at every point of\/
$d_K(F)\cup\{\infty'\}$.
\item
$\cH^i(F')=0$\ \ for every $i\in\Z\setminus\{0,1\}$.
\item
$R^{-1}\Phi_{\eta_z}F'=0$\ \ for every $z\in A'(K)$.
\end{enumerate}
\end{theorem}
\begin{proof} Assertion (iii) is obvious.

(i): For $F$ and $F'$ as in the theorem, let $\Sigma(F)$ denote the
support of $\cH^1(F')$; by \cite[Cor.8.2.4]{Hu2}, $\cH^i(F')$ is
an overconvergent $\Lambda$-module for every $i\in\Z$, especially,
$\Sigma(F)$ is the set of all specializations (in $A'$) of the subset
$\Sigma_0(F)$ of all the maximal points of $A'$ lying in $\Sigma(F)$.
Now, let $a\in\Sigma_0(F)$ be any point, and $(L,|\cdot|_L)$ the
algebraic closure of the residue field of $a$; after base change
along the natural morphism $\Spa(L,L^+)\to S$, we may
assume that $L=K$ (corollary \ref{cor_drop-is-finite}), in which
case the condition $\cH^1(F')_a\neq 0$ means that
\set\begin{equation}\label{eq_a-in-support}
H^2_c(U,F\otimes_\Lambda\cL_\psi\La m_a\Ra_{|U})\neq 0.
\end{equation}
Let $0\to F_1\to F\to F_2\to 0$ be a short exact sequence of
$\Lambda$-modules on $U_\et$; from \eqref{eq_a-in-support} it
is clear that $\Sigma(F)\subset\Sigma(F_1)\cup\Sigma(F_2)$.
In view of proposition \ref{prop_drop-is-finite}(iii), we reduce
to prove assertion (i) for $F_1$ and $F_2$.
Especially, if $\fn\subset\Lambda$ is any maximal ideal,
assertion (i) for $F$ follows from the same assertion for $\fn F$
and $F/\fn F$; hence, an easy induction on the length of $\Lambda$
reduces to the case where $\Lambda$ is a field.
In this case, lemma \ref{lem_PD-lambda}(iii) says that
\eqref{eq_a-in-support} is equivalent to :
$$
H^0(U,F^\vee\otimes_\Lambda\cL_\psi\La m_{-a}\Ra)\neq 0
\qquad
\text{where $F^\vee:=\cHom_\Lambda(F,\Lambda_U)$}
$$
especially $(F^\vee)_{[-a]}^\natural(0)\neq 0$, and
therefore $-a\in d_K(F^\vee)$ (remark \ref{rem_when-0-drops}), whence
$a\in d_K(F)$, according to proposition \ref{prop_drop-is-finite}(iv).
This shows that $\Sigma_0(F)\subset d_K(F)$, whence (i), in view of
proposition \ref{prop_drop-is-finite}(i).

(ii): Let $(\rho_n~|~n\in\N)$ be an increasing sequence of
values in $|K^\times|$, such that $\rho_n\to+\infty$ for
$n\to+\infty$, and such that $A\setminus U\subset\D(0,\rho_0)$.
For every $n\in\N$, let $j_n:\D(0,\rho_n)\to A$ be the
open immersion, and set
$$
H_n:=\cH^0(\cF_\psi(j_{n!}j^*_nF)).
$$

\begin{claim}\label{cl_building-the-sys}
$H_n\subset H_{n+1}$ for every $n\in\N$.
\end{claim}
\begin{pfclaim} The assumption can be checked on the stalks, hence
let $a\in A'$ be any point; set
$$
V_n:=p^{\prime-1}(a)\cap(\D(0,\rho_n)\times_SA')
\qquad\text{and}\qquad
Z_n:=V_{n+1}\setminus V_n.
$$
The assumption on $\rho_0$ implies that the restriction of
$L:=p^*F\otimes_\Lambda\cL_\psi\La m\Ra$ is locally constant on
the \'etale site of the pseudo-adic space
$(A\times_SA',Z_n)$.
However, let $i_a:Z_n\to V_{n+1}$ (resp. $j_a:V_n\to V_{n+1}$)
be the closed (resp. open) immersion of pseudo-adic spaces; we
have a short exact sequence of $\Lambda$-modules on
$(A\times_SA',V_{n+1})_\et$ :
$$
0\to j_{a!}L_{|V_n}\to L_{|V_{n+1}}\to i_{a*}L_{|Z_n}\to 0
$$
whence, by \cite[Cor.5.4.8]{Hu2}, a cohomology exact sequence :
$$
H^0_c(Z_n,L)\to H_{n,a}\to H_{n+1,a}.
$$
Since $Z_n$ is not proper over $S$, the left-most term vanishes,
and the claim follows.
\end{pfclaim}

\begin{claim}\label{cl_deduce-loc-constance}
Let $X$ be an adic space, locally of finite type over $S$,
and $G_\bullet:=(G_n~|~n\in\N)$ a direct system of constructible
$\Lambda$-modules on $X_\et$, with injective transition maps
$G_n\to G_{n+1}$ (for every $n\in\N$). Suppose that :
\begin{enumerate}
\alphaenu
\item
the colimit $G$ of the system $G_\bullet$ is overconvergent;
\item
there exists $l\in\N$ such that $\length_\Lambda G_\xi=l$ for every
geometric point $\xi$ of $X$.
\end{enumerate}
Then $G$ is a locally constant $\Lambda$-module of finite type on $X_\et$.
\end{claim}
\begin{pfclaim} Let $x\in X$ be any point, and $\bar x$ a geometric
point of $X$ localized at $x$; since $G_{\bar x}$ has finite length,
we may find $n\in\N$ large enough, so that $G_{n,\bar x}=G_{\bar x}$.
By assumption, there exists a locally closed constructible subset
$Z\subset X$ containing $x$, and such that $G_n$ restricts to
a locally constant $\Lambda$-module on $Z_\et$; after shrinking
$Z$, we may then assume that the stalks $G_{n,\xi}$ are isomorphic
to $G_{n,\bar x}$, for every geometric point $\xi$ of $Z$.
Let $\xi$ be any geometric point of $Z$; by assumption
$G_{n,\xi}\subset G_{m,\xi}$ for every $m\geq n$, therefore
$l=\length_\Lambda G_{n,\xi}\leq\length_\Lambda G_{m,\xi}\leq l$.
It follows that $G_{n,\xi}=G_{m,\xi}$, hence $G_{n|Z}=G_{m|Z}$
for every $m\geq n$, so finally $G_{|Z}=G_{n|Z}$. This shows that
$G$ is constructible, and then the claim follows from
\cite[Lemma 2.7.11]{Hu2}.
\end{pfclaim}

Clearly $\cH^0(F')$ is the colimit of the direct system
$(H_n~|~n\in\N)$, which consists of constructible $\Lambda$-modules,
according to \cite[Th.6.2.2]{Hu2}. On the other
hand, (i), (iii) and lemma \ref{lem_Groth-Ogg-Shaf} easily imply
that the length of the stalk $\cH^0(F')_\xi$ is constant for $\xi$
ranging over the geometric points of $U'$ whose support is
a maximal point, hence also for $\xi$ ranging over all the geometric
points of $U'$, since -- as it has already been noted --
$\cH^0(F')$ is overconvergent. Then the first assertion of (ii)
follows from claims \ref{cl_building-the-sys} and
\ref{cl_deduce-loc-constance}. Finally, (i), (iii) and
\eqref{eq_stalk-at-0} imply that the Euler-Poincar\'e characteristic
of $\cH^0(F')$ is finite, hence its Swan conductor at every point
of $d_K(F)\cup\{\infty'\}$ is finite, in view of lemma
\ref{lem_Groth-Ogg-Shaf}, therefore $\cH^0(F')$ has bounded
ramification everywhere.

(iv): In view of (iii) we have $R^{-1-i}\Phi_{\eta_z}\cH^i(F')=0$
for every $z\in A'(K)$ and every $i\neq 0$, hence we are reduced
to showing that $R^{-1}\Phi_{\eta_z}\cH^0(F')=0$ (proposition
\ref{prop_vanish-properties}(iii) and corollary
\ref{cor_vanish-properties}). In view of (ii), (and corollary
\ref{cor_vanish-smooth}) the assertion is already known for
every $z\in U'(K)$. Hence, we may assume that $z\in d_K(F)$, in
which case we have to prove that the natural map
$\cH^0(F')_z\to\cH^0(F')_{\eta_z}$ is injective. The latter
will follow, once we have shown that $H^0_c(A'_\et,\cH^0(F'))=0$.
Again by (iii), this is the same as showing that $H^0_c(A'_\et,F')$
vanishes, which holds, due to \eqref{eq_stalk-at-0}.
\end{proof}

\begin{remark}\label{rem_perversity}
(i)\ \ For $a\in K=A(K)$, let $i_a:S\to A$ be the
closed immersion with image $\{a\}$, and $\pi:A\to S$ the
structure morphism. If $G^\bullet$ is any object of
$\sD^+(S_\et,\Lambda)$, a direct computation yields a natural
isomorphism :
$$
\cF_\psi(i_{a*}G^\bullet)\isom
\pi^*G^\bullet\otimes_\Lambda\cL_\psi\La a\cdot x'\Ra[1].
$$

(ii)\ \ In the situation of \eqref{subsec_set-up_-for-th}, let $H$
be a Zariski constructible $\Lambda$-module on $A_\et$, such that
$j^*H$ is locally constant on $U_\et$, and has bounded ramification
at the points of $D\setminus U$. Notice that the kernel and cokernel
of the unit of adjunction $H\to j_*j^*H$ are direct sum of sheaves
of the form $i_{a*}G$ (for various $\Lambda$-modules $G$, and with
$a$ running over the points of $A\setminus U$).
Since $\cF_\psi$ is a triangulated functor, it follows easily
from (i) that theorem \ref{th_perversity}(i),(ii) still holds
with $F'$ replaced by $\cF_\psi(H)$.

(iii)\ \ On the other hand, in the situation of (ii), we have
in general $\cH^i(\cF_\psi(H))=0$ only for $i\in\Z\setminus\{-1,0,1\}$.
More precisely, $\cH^{-1}(\cF_\psi(H))=0$ if and only if
$\Gamma_c(A_\et,H)=0$, and in this case, also assertion (iv)
of theorem \ref{th_perversity} still holds.
\end{remark}

We point out the following result, even though it will not
be used in the rest of this work.

\begin{proposition}\label{prop_perversity}
Let $j:U\to A$ be a Zariski open immersion, and $F$ a locally
constant $\Lambda$-module on $U_\et$, with bounded ramification
at all points of $D\setminus U$. Suppose that $F$ does not admit
any subquotient of the type
$\cL_\psi\La xa\Ra_{|U}\otimes_\Lambda\bar\Lambda$, where
$\bar\Lambda$ is any residue field of $\Lambda$, and $a\in K$
is any element. Then there exists a Zariski open immersion
$j':U'\to A'$, and a locally constant $\Lambda$-module $G$
on $U'_\et$ such that $\cF_\psi(j_!F)=j'_*G[0]$.
\end{proposition}
\begin{proof} First, we prove that $\cH^i\cF_\psi(j_!F)=0$
for $i\neq 0$. Indeed, by virtue of remark \ref{rem_perversity}(iii),
the assertion is already known for all $i\neq 0,1$; hence it
suffices to show :

\begin{claim}\label{cl_vanish-H-2}
Under the assumptions of the lemma, we have
$H^2_c(U_\et,F\otimes_\Lambda\cL_\psi\La xa\Ra)=0$ for every
$a\in K$.
\end{claim}
\begin{pfclaim} Notice that $F$ satisfies the assumptions
of the lemma, if and only if the same holds for
$F\otimes_\Lambda\cL_\psi\La xa\Ra$, for every $a\in K$.
Hence, we are reduced to showing that $H^2_c(U_\et,F)=0$.
Since $\Lambda$ is artinian, we may assume that $\Lambda$
is local, say with maximal ideal $\fm$. Set $\bar F:=F/\fm F$;
under our assumptions, $\bar F$ does not admit constant
quotients, hence
$H^0(U_\et,\cHom_\Lambda(\bar F,\Lambda/\fm\Lambda))=0$,
and therefore $H^2_c(U_\et,\bar F)=0$, by Poincar\'e duality
(lemma \ref{lem_PD-lambda}(iii)). This easily implies that
the natural map $\alpha:H^2_c(U_\et,\fm F)\to H^2_c(U_\et,F)$
is surjective. Now, let $t_1,\dots,t_k$ be a
system of generators $\fm$; we may then define a map
of $\Lambda$-modules $\beta:F^{\oplus k}\to F$ in the obvious
way, with image equal to $\fm F$. It follows easily that
the image of $\alpha$ equals the image of the induced map
$H^2_c(\beta):H^2_c(U_\et,F)^{\oplus k}\to H^2_c(U_\et,F)$;
however, the latter is just $\fm H^2_c(U_\et,F)$. By Nakayama's
lemma, the claim follows.
\end{pfclaim}

By theorem \ref{th_perversity}(ii), there exists a Zariski open
immersion $j':U'\to A'$, such that $G:=j^{\prime*}\cH^0\cF_\psi(j_!F)$
is a locally constant $\Lambda$-module, with bounded ramification
at all the points of $D'\setminus U'$. Let $C$ denote the cone
of the natural morphism $\cF_\psi(j_!F)\to j'_*G[0]$. Then
$\cH^iC=0$ for every $i\neq -1,0$, and $C$ is supported on
finitely many $K$-rational points. In view of remark
\ref{rem_perversity}(i), it follows easily that $\cF'_\psi(C)$
is concentrated in degrees $-2,-1$, and we have an
exact sequence :
$$
\cH^{-2}\cF'_\psi(j'_*G[0])\to\cH^{-2}\cF'_\psi(C)\to
\cH^{-1}(\cF'_\psi\circ\cF_\psi(j_!F))
$$
whose first (resp. third) term vanishes, due to theorem
\ref{th_perversity}(iii) (resp. due to \eqref{eq_FT-is-selfdual});
hence the middle term vanishes as well. Likewise, we have
an exact sequence :
$$
\cH^{-1}\cF'_\psi(j'_*G[0])\to\cH^{-1}\cF'_\psi(C)\xrightarrow{\beta}
j_!F\to\cH^0\cF'_\psi(j'_*G[0])\to 0.
$$
However, remark \ref{rem_perversity}(i) also implies that
$\cH^{-1}\cF'_\psi(C)$ is a quotient of a sum of $\Lambda$-modules
of the type $\cL_\psi\La xa\Ra$, for various $a\in K$; under
our assumptions, it follows that the map $\beta$ vanishes;
furthermore, $\cH^{-1}\cF'_\psi(j'_*G[0])=0$, again by theorem
\ref{th_perversity}(iii). Therefore $\cH^{-1}\cF'_\psi(C)$ must
vanish as well, so that $\cF'_\psi(C)=0$, and finally $C=0$,
whence the assertion.
\end{proof}

\subsection{Stationary phase}\label{sec_Station}
Keep the notation of \eqref{subsec_intro-Fourier}, and let
$\eps\in|K^\times|$ be any value; set
$$
X(\eps):=
D\times_S\D(\infty',\eps)\subset D\times_SD'
$$
and denote by $\bar p_\eps:X(\eps)\to D$,
$\bar p{}'_\eps:X(\eps)\to\D(\infty',\eps)$ the projections; here
the radius $\eps$ is meant relative to the global coordinate $x'$
on $A'$; in other words :
\set\begin{equation}\label{eq_disc-at-infty}
\D(\infty',\eps)=\{\infty'\}\cup\{b\in A'~|~|x'(b)|\geq 1/\eps\}.
\end{equation}
The vanishing cycle construction of \eqref{subsec_vanish} applies
to the current situation, and yields functors :
$$
R^i\Psi_{\eta_{\infty'},\eps}:
\sD^+(X(\eps)_\et,\Lambda)\to\Lambda[\pi(\infty',\eps)]_D\Mod
\qquad
\text{for every $i\in\N$}.
$$
(Where $\Lambda[\pi(\infty',\eps)]_D\Mod$ denotes the category of
$\Lambda$-modules on $D_\et$, endowed with a linear action of
$\pi(\infty',\eps)$; notice that $\bar p{}_\eps^{\prime-1}(\infty')$
is naturally isomorphic to $D$.) Let also
$$
\bar\cL_\psi\La m\Ra:=(\alpha\times_S\alpha')_!\cL_\psi\La m\Ra
$$
which is a $\Lambda$-module on the \'etale site of $D\times_SD'$.
The following result is the counterpart of the ``universal local
acyclicity'' of \cite[Th.1.3.1.2]{Lau2}.

\begin{theorem}\label{th_cycles-vanish}
Let $F$ be a $\Lambda$-module on $D_\et$, and $j:U\to D$ an open
immersion such that $j(U)\subset A$ and $j^*F$ is a constructible
locally constant $\Lambda$-module on $U_\et$. Then
$$
j^*R\Psi_{\eta_{\infty'},\eps}
(\bar p{}_\eps^*F\otimes_\Lambda\bar\cL_\psi\La m\Ra_{|X(\eps)})=0
\qquad
\text{in\ \  $\sD(U_\et,\Lambda)$.}
$$
\end{theorem}
\begin{proof} Suppose that $0\to F_1\to F\to F_2\to 0$ is a
short exact sequence of $\Lambda$-modules on $A_\et$, such that
$j^*F_i$ is locally constant, for $i=1,2$; since
$R\Psi_{\eta_{\infty'},\eps}$ is a triangulated functor, it is
clear that the stated vanishings for $F_1$ and $F_2$ imply the
same vanishing for $F$. Since $\Lambda$ is a $\Lambda$-module
of finite length, an easy induction allows then to reduce to
the case where $\Lambda$ has length one, {\em i.e.} we may
assume that $\Lambda$ is a field.

Set $G:=\bar p{}_\eps^*F\otimes_\Lambda\bar\cL_\psi\La m\Ra_{|X(\eps)}$;
choose a fundamental pro-covering
$(C_H~|~H\subset\pi(\infty',\eps))$ of $\D(\infty',\eps)$,
and for every open subgroup $H\subset\pi(\infty',\eps)$ let
$G_H:=Rj_{H*}j_H^*G$ (notation of \eqref{subsec_vanish}).
We will show more precisely that $j^*G_H=0$ for every such $H$.

Notice first that $j_H$ is partially proper, and both $\bar p{}_\eps^*F$
and $\bar\cL_\psi\La m\Ra$ restrict to overconvergent $\Lambda$-modules
on the \'etale site of $\bar p{}_\eps^{-1}U$, hence $\cH^ij^*G_H$ is
overconvergent for every $i\in\Z$, and it suffices to show that :
$$
\cH^i(G_H)_z=0
\qquad
\text{for every maximal point $z\in U$ and every $i\in\N$}.
$$
Fix such a maximal point $z\in U$; in view of
\cite[Prop.2.4.4, Lemma 2.5.12, Prop.2.5.13(i)]{Hu2}, the point
$(z,\infty')\in U\times_S\D(\infty',\eps)$ admits a fundamental
system of \'etale neighborhoods of the form
$$
V\times_S\D(\infty',\delta)
\qquad
\text{for every $V\in\cV$ and every $\delta\in|K^\times|$ with
$\delta<\eps$}
$$
where $\cV$ is a fundamental system of quasi-compact \'etale
neighborhoods of $z$ in $U_\et$. Moreover, since $z$ is maximal
in $U$, we may assume that the structure morphism $q_V:V\to U$
induces a finite \'etale covering $V\to q_V(V)$, for every
$V\in\cV$ (\cite[Prop.1.5.4, Cor.1.7.4]{Hu2}). For every $\delta$
as above, let
$C_H(\delta):=C_H\times_{\D(\infty',\eps)}\D(\infty',\delta)$.
It follows easily that :
$$
\cH^i(G_H)_z\simeq\colim{0<\delta\leq\eps}\,\colim{V\in\cV}\,
H^i(V\times_S C_H(\delta),j_H^*G).
$$
Up to replacing $\cV$ by a cofinal system, we may assume that
$F_{|V}$ is a constant $\Lambda$-module for every $V\in\cV$,
in which case -- since $\Lambda$ is a field -- $F_{|V}$ is the
direct sum of finitely many copies of $\Lambda_U$.
We are then easily reduced to the case where $F_{|U}=\Lambda_U$,
therefore $G=\bar\cL_\psi\La m\Ra_{|X(\eps)}$.

Fix any $V\in\cV$; we shall show more precisely that there exists
$\delta_0>0$ small enough, so that :
$$
H^i(V\times_SC_H(\delta),\cL_\psi\La m\Ra)=0
\qquad
\text{for every $\delta\in|K^\times|$ with $\delta\leq\delta_0$
and every $i\in\N$}.
$$
For every $\delta$ as above, denote by
$q_\delta:V\times_SC_H(\delta)\to C_H(\delta)$ the projection;
clearly it suffices to show that $Rq_{\delta*}\cL_\psi\La m\Ra=0$
for every sufficiently small $\delta>0$. However, $q_\delta$ is
a smooth, separated and quasi-compact morphism, therefore, for
every $i\in\N$ the $\Lambda$-module $R^iq_{\delta*}\cL_\psi\La m\Ra$
is constructible on the \'etale site of $C_H(\delta)$
(\cite[Th.6.2.2]{Hu2}); the latter is locally of finite type over
$S$, hence the subset of its maximal points is everywhere dense
for the constructible topology (\cite[Cor.4.2]{Hu1}), so we reduce
to showing that $(R^iq_{\delta*}\cL_\psi\La m\Ra)_a=0$ for every
sufficiently small $\delta>0$, every $i\in\N$ and every maximal point
$a\in C_H(\delta)$. Given $\delta>0$ as above, and a maximal
point $a\in C_H(\delta)$, let $(L,|\cdot|_L)$ be the algebraic closure
of the residue field of $a$, with its rank one valuation; after base
change along the natural morphism $\Spa(L,L^+)\to S$, we may assume
that $L=K$, in which case the image $b\in\D(\infty',\delta)^*$ of $a$
is a $K$-rational point. Moreover, let
$i_a:q_\delta^{-1}(a)\to V\times_SC_H(\delta)$ be the immersion;
\cite[Th.4.1.1(c)]{Hu2} says that
$$
(R^iq_{\delta*}\cL_\psi\La m\Ra)_a\simeq
H^i(q_\delta^{-1}(a)_\et,i_a^*\cL_\psi\La m\Ra)\simeq
H^i(V_\et,\cL_\psi\La m_b\Ra)
$$
where $m_b:A\to(\A^1_K)^\ad$ corresponds to the global section
$bx\in\cO_A(A)$ (recall that $x$ is a fixed global coordinate on $A$;
also, $b$ is here identified to an element of $K^\times$, via
the global coordinate $x'$, especially, $|b|\geq 1/\delta$).
By lemma \ref{lem_PD-lambda}(iii), it is clear that
\set\begin{equation}\label{eq_for-i>1}
H^i(V_\et,\cL_\psi\La m_b\Ra)=0
\qquad
\text{for every $i>1$ and every $\delta>0$ as above}.
\end{equation}
To deal with the cases $i=0,1$, we shall compute explicitly the
Euler-Poincar\'e characteristic
$$
\chi(V,\cL_\psi\La m_b\Ra):=
\sum_{i=0}^2(-1)^i\cdot\length_\Lambda H^i(V_\et,\cL_\psi\La m_b\Ra).
$$
Let us remark that
\set\begin{equation}\label{eq_proj-formula}
Rq_{V*}\cL_\psi\La m_b\Ra_{|q_V(V)}\simeq
q_{V*}\cL_\psi\La m_b\Ra_{|q_V(V)}\simeq
\cL_\psi\La m_b\Ra_{|q_V(V)}\otimes_\Lambda q_{V*}\Lambda_V.
\end{equation}
We may write $q_V(V)=\D(y_0,\rho_0)-\bigcup_{i=1}^n\E(y_i,\rho_i)$
where $y_0,\dots,y_n\in U$ are certain $K$-rational points,
and $\rho_0,\dots,\rho_n\in|K^\times|$, with $\rho_i\leq\rho_0$ and :
\set\begin{equation}\label{eq_def-E_i}
\E(y_i,\rho_i):=\{y\in A~|~|(x-y_i)(y)|<\rho_i\}
\qquad
\text{for every $i=1,\dots,n$}.
\end{equation}
Clearly we may assume that
$\E(y_i,\rho_i)\cap\E(y_j,\rho_j)=\emptyset$ for $i\neq j$, in
which case 
\set\begin{equation}\label{eq_EP-for-qV}
\chi(q_V(V),\Lambda)=\chi(\D(y_0,\rho_0),\Lambda)-
\sum_{i=1}^n\chi(\E(y_i,\rho_i),\Lambda)=1-n
\end{equation}
(where -- for any locally closed constructible subset $M$ of $A$ --
we denote by $\chi(M,\Lambda)$ the Euler-Poincar\'e characteristic
of the constant $\Lambda$-module $\Lambda_M$ on $M_\et$).
On the other hand, let $\bar V$ be the topological closure of
$V$ in $A$; then $\bar V\setminus V=\{w_0,\dots,w_n\}$, where
$\{w_0\}\cup\D(y_0,\rho_0)$ is the topological closure of
$\D(y_0,\rho_0)$, and $\E(y_i,\rho_i)\setminus\{w_i\}$ is the
interior of $\E(y_i,\rho_i)$, for every $i=1,\dots,n$. To the
$\Lambda$-module $q_{V*}\Lambda_V$ and each point $w_i$, the
theory of \cite{Hu3} associates a sequence of {\em breaks\/}
$$
\beta_{i1},\dots,\beta_{im_i}\in\Gamma_{\!w_i}\otimes_\Z\Q
$$
where $\Gamma_{\!w_i}$ is the valuation group of the residue
field of $w_i$; likewise, every stalk $\cL_\psi\La m_b\Ra_{w_i}$
is a rank one $\Lambda$-module, hence it admits a single
break $\alpha_i(b)$. Define the group homomorphism
$$
\sharp_\Q:\Gamma_{w_i}\otimes_\Z\Q\to\Q
$$
as in \cite[\S1]{Hu3}. With this notation, we have the following :

\begin{claim}\label{cl_lower-bound}
There exists $\delta_0\in|K^\times|$ such that :
\begin{enumerate}
\item
$\beta_{ij}>\alpha_i(b)$ for every $i=0,\dots,n$, and every
$j=1,\dots,m_i$
\item
$\sharp_\Q\alpha_0(b)=1$ and $\sharp_\Q\alpha_i(b)=-1$ for every
$i=1,\dots,n$
\end{enumerate}
for every $b\in K^\times$ with $|b|\geq 1/\delta_0$.
\end{claim}
\begin{pfclaim} Everything follows easily from \cite[Lemma 9.4]{Hu3}.
\end{pfclaim}

In light of \eqref{eq_proj-formula}, claim \ref{cl_lower-bound}(i)
and \cite[\S8.5(ii)]{Hu3}, we deduce that -- whenever
$|b|\geq 1/\delta_0$ -- the value $\alpha_i(b)$ is the unique break
of the stalk $q_{V*}\cL_\psi\La m_b\Ra_{w_i}$. Moreover, claim
\ref{cl_lower-bound}(ii) already ensures that :
\set\begin{equation}\label{eq_for-i=0}
H^0(V_\et,\cL_\psi\La m_b\Ra)=0
\qquad
\text{for every $b\in K^\times$ with $|b|\geq 1/\delta_0$}.
\end{equation}
Furthermore, let $d$ be the degree of the finite morphism $q_V$;
claim \ref{cl_lower-bound}(ii), \eqref{eq_EP-for-qV}, and
\cite[Cor.8.4(ii), Lemma 10.1, Cor.10.3]{Hu3} yield the identity :
$$
\begin{array}{r@{\:=\:}l}
\chi(V,\cL_\psi\La m_b\Ra) & \chi(q_V(V),q_{V*}\cL_\psi\La m_b\Ra) \\
& d\cdot\chi(q_V(V),\Lambda)-d\cdot\sum_{i=0}^n\sharp_\Q\alpha_i(b) \\
& d\cdot(1-n)-d\cdot(1-n)=0
\end{array}
$$
whenever $|b|\geq 1/\delta_0$. In view of \eqref{eq_for-i>1} and
\eqref{eq_for-i=0}, we deduce that $H^1(V_\et,\cL_\psi\La m_b\Ra)$
vanishes as well, and the theorem follows.
\end{proof}

\sset\subsubsection{}\label{subsec_as-expla}
For every $z\in D(K)=K\cup\{\infty\}$, as explained in
\eqref{subsec_was-remark}, the choice of a maximal geometric
point $\xi_z$ of $\D(z,\eps)^*$ induces an equivalence :
$$
\tau_z:\Lambda[\pi(z,\eps)]\Mod_\mathrm{f.cont}\isom
\Lambda_{\D(z,\eps)^*}\Mod_\mathrm{loc}
$$
whose quasi-inverse is the stalk functor $F\mapsto F_{\xi_z}$.
For every continuous $\Lambda[\pi(z,\eps)]$-module of finite type
$M$, denote by $M_!$ the $\Lambda$-module on $D_\et$ obtained as
extension by zero of $\tau_zM$.

Let $i_z:S\to\bar p{}_\eps^{\prime-1}(\infty')$ be the unique
$S$-morphism whose image is the point $(z,\infty')$. In view
of remark \ref{rem_cont-action}, we may define the
{\em $\eps$-local Fourier transform from $z$ to $\infty'$}
as the functor :
$$
\cF^{(z,\infty')}_{\psi,\eps}:
\Lambda[\pi(z,\eps)]\Mod_\mathrm{f.cont}\to
\Lambda[\pi(\infty',\eps)]\Mod_\mathrm{cont}
$$
given by the rule :
$$
M\mapsto
i_z^*R^1\Psi_{\eta_{\infty'},\eps}
(\bar p{}^*_\eps(M_!)\otimes_\Lambda\bar\cL_\psi\La m\Ra_{|X(\eps)}).
$$
Finally, in the situation of \eqref{subsec_set-up_-for-th}, denote
by $F'_\eps$ the complex of $\Lambda$-modules on $\D(\infty',\eps)_\et$
obtained as extension by zero of $\cF_\psi(j_!F)_{|\D(\infty,\eps)^*}$.
With this notation, we have :

\begin{theorem}\label{th_Stationary-phase} {\em (Stationary Phase)}\ \ 
In the situation of \eqref{subsec_set-up_-for-th},
suppose moreover that $F$ has bounded ramification at all the
(finitely many) points of $D\setminus U$. Then for every
$\eps\in|K^\times|$ there is a natural $\pi(\infty',\eps)$-equivariant
decomposition :
$$
R^0\Psi_{\eta_{\infty'},\eps}\cH^0(F'_\eps)\isom
\bigoplus_{z\in D\setminus U}\cF^{(z,\infty')}_{\psi,\eps}(F_{\xi_z}).
$$
\end{theorem}
\begin{proof} Resume the notation of \eqref{sec_Station}, and
let as well $\bar\jmath:U\to D$  be the open immersion.
By inspecting the definitions we get a natural isomorphism :
$$
F'_\eps\isom
R\bar p{}'_{\eps*}(\bar p{}^*_\eps
(\bar\jmath_!F)\otimes_\Lambda\bar\cL_\psi\La m\Ra_{|X(\eps)})[1]
$$
from which, by proposition \ref{prop_vanish-properties}(ii),(iv),
there follows an equivariant isomorphism :
\set\begin{equation}\label{eq_two-vanish}
R\Psi_{\eta_{\infty'},\eps}F'_\eps\isom
R\Gamma(D,R\Psi_{\eta_{\infty'},\eps}
(\bar p{}^*_\eps(
\bar\jmath_!F)\otimes_\Lambda\bar\cL_\psi\La m\Ra_{|X(\eps)})[1]).
\end{equation}
To evaluate the left-hand side of \eqref{eq_two-vanish}, we use
the spectral sequence of proposition \ref{prop_vanish-properties}(iii):
combining with remark \ref{rem_perversity}(iii) we obtain the
equivariant isomorphism 
$$
R^0\Psi_{\eta_{\infty'},\eps}F'_\eps\isom
R^0\Psi_{\eta_{\infty'},\eps}\cH^0(F'_\eps).
$$
Lastly, by theorem \ref{th_cycles-vanish}, the vanishing
cycles appearing on the right-hand side of \eqref{eq_two-vanish}
are supported on $D\setminus U$, whence the theorem (details
left to the reader).
\end{proof}

\begin{remark}\label{rem_not-much}
(i)\ \
Not much can be said concerning the $\eps$-local Fourier transforms
appearing in theorem \ref{th_Stationary-phase}. The situation
improves when we take the limit for $\eps\to 0$.
Namely, choose a sequence $(\eps_n~|~n\in\N)$ and geometric points
$(\xi_n~|~n\in\N)$ as in \eqref{subsec_pass-to-lim}, and define 
$$
\cF_\psi^{(z,\infty')}(M):=
\colim{n\in\N}\cF_{\psi,\eps_n}^{(z,\infty')}(M_n)
\qquad
\text{for every continuous $\Lambda[\pi(z,\eps_0)]$-module $M$}
$$
where $M_n$ denotes the $\Lambda[\pi(z,\eps_n)]$-module which
is the image of $M_0:=M$, under the composition of the functors
\eqref{eq_dir-syst-cats}. It is clear that $\cF^{(z,\infty')}_\psi(M)$
depends only on the image of $M$ in $\Lambda[\pi(z)]\Mod$;
hence -- and in view of theorem \ref{th_Stationary-phase} -- we
deduce a natural decomposition (notation of
\eqref{subsec_more-trivialities}) :
\set\begin{equation}\label{eq_stat-phase}
\cH^0(\cF_\psi(j_!F))_{\eta_{\infty'}}\isom
\bigoplus_{z\in D\setminus U}
\cF^{(z,\infty')}_\psi(F_{\eta_z})
\end{equation}
for every $\Lambda$-module with bounded ramification $F$ on
$U_\et$ (see remark \ref{rem_one-more-obviety}). However,
remark \ref{rem_perversity}(ii) and the discussion in
\eqref{subsec_more-trivialities} imply that the left-hand
side of \eqref{eq_stat-phase} has finite length, hence the
same holds for each of the summands on the right-hand side.

(ii)\ \ 
We shall later need also a localized version of the stationary
phase identity. Namely, let $\delta\in|K^\times|$, and suppose
that $F$ is a locally constant $\Lambda$-module of finite type
on $\D(0,\delta)^*_\et$. Let $w\in\bar p{}^{\prime-1}_\eps(\infty')$
be the unique point such that $(\D(0,\delta)\times\{\infty'\})\cup\{w\}$
is the topological closure of $\D(0,\delta)\times\{\infty'\}$, and set
$$
\cF_\psi^{(w,\infty')}(F):=
\Gamma(\{w\},R^1\Psi_{\eta_{\infty'}}
(\bar p{}^*_\eps(\bar\jmath_!F)\otimes_\Lambda
\bar\cL_\psi\La m\Ra_{|X(\eps)}))
$$
where $\bar\jmath:\D(0,\delta)^*\to D$ is the open immersion
($\{w\}$ is here seen as a pseudo-adic space).
Let also $F'_\eps$ denote the complex of $\Lambda$-modules on
$\D(\infty',\eps)_\et$ obtained as extension by zero of
$\cF_\psi((\bar\jmath_!F)_{|A})$, and notice that we still have
$\cH^i(F'_\eps)=0$ for $i\notin\{0,1\}$, hence the proof of theorem
\ref{th_Stationary-phase} still applies to this $\Lambda$-module
$F$, and we get the identity :
$$
\cH^0(F'_\eps)_{\eta_{\infty'}}\isom
\cF^{(0,\infty')}_\psi(F_{\eta_0})\oplus\cF_\psi^{(w,\infty')}(F).
$$
\end{remark}

\sset\subsubsection{}\label{subsec_Garuti}
In order to apply theorem \ref{th_Stationary-phase}, we shall
need a local-to-global extension result for
$\Lambda[\pi(z,\eps)]$-modules, analogous to Gabber's theorem
\cite[Th.1.4.1]{Kat}. This is contained in the following :

\begin{proposition}\label{prop_Garuti}
Let $F$ be a locally constant $\Lambda$-module on
$\D(z,\eps)^*_\et$ of finite type. Then there exists a Zariski
open subset $U\subset D\setminus\{z\}$, and a locally constant
$\Lambda$-module $\tilde F$ on $U_\et$, such that:
\begin{enumerate}
\alphaenu
\item
the restriction $\tilde F_{|\D(z,\eps)^*}$ is isomorphic to $F$;
\item
For every sufficiently small $\delta\in|K^\times|$, and every
$w\in D\setminus(U\cup\{z\})$, the restriction $\tilde F_{|D(w,\delta)}$
is a tamely ramified locally constant $\Lambda$-module
(see \eqref{subsec_tame-ram}).
\end{enumerate}
\end{proposition}
\begin{proof} Denote by $C\subset\D(z,\eps)^*$ the annulus
$\{a\in\ A~|~|(x-z)(a)|_a=\eps\}$, and let $\D'$ be another
disc containing $C$, such that $\D'\cup\D(z,\eps)=D\setminus\{z\}$.
We may find a finite Galois covering $\phi:C'\to C$, such that
$F_{|C'}$ is a constant $\Lambda$-module (see
\eqref{subsec_was-remark}). Then, by \cite[Th.1]{Gar}, $\phi$
extends to a finite Galois covering $C''\to\D'$, with the same
Galois group, and ramified over finitely many points
$w_1,\dots,w_n\in\D'(K)$. It follows easily that $F_{|C}$ extends
to a locally constant $\Lambda$-module $F'$ on the \'etale site of
$\D'\setminus\{w_1,\dots,w_n\}$. By construction, it is clear
that $F'$ has finite monodromy around each of the points
$w_1,\dots,w_n$. We may then glue $F$ and $F'$ along their
common restriction on $C$, to obtain the sought $\tilde F$.
\end{proof}

\sset\subsubsection{}\label{subsec_def-M}
From \eqref{subsec_more-trivialities} and theorem
\ref{th_perversity}(ii), we know that the source
of the isomorphism \eqref{eq_stat-phase} is an object of
$\Lambda[\pi(\infty')]\Mod$; we wish to show that each
direct summand in the target is also naturally an object in
the same category, so that \eqref{eq_stat-phase} is an
isomorphism of $\Lambda[\pi(\infty')]$-modules.

To this aim, let $z\in D(K)$ be any point, $F$ a $\Lambda$-module
on $\D(z,\eps_0)^*_\et$.
For every $\eps,\delta\in|K^\times|$ with $\eps\leq\eps_0$,
choose as usual a fundamental pro-covering
$(C_H~|~H\subset\pi(\infty',\delta))$, and consider the
commutative diagram :
$$
\xymatrix{
\D(\infty',\delta)^* \ddouble
& \D(z,\eps)^*\times_S\D(\infty',\delta)^* \ar[d]_{j_{\eps,\delta}}
\ar[r]^-{q_{\eps,\delta}} \ar[l]_-{q'_{\eps,\delta}} &
\D(z,\eps)^* \ar[d]^-{j_\eps} \\
\D(\infty',\delta)^* &
\D(z,\eps)\times_S\D(\infty',\delta)^*
\ar[r]^-{\bar q_{\eps,\delta}} \ar[l]_-{\bar q{}'_{\eps,\delta}} &
\D(z,\eps) \\
& \D(z,\eps)\times_SC_H \ar[u]^-{j_H} & 
}$$
where :
\begin{itemize}
\item
$j_H$ is obtained by base change from the projection
$C_H\to\D(\infty',\delta)^*$
\item
$q_{\eps,\delta}$, $q'_{\eps,\delta}$, $\bar q_{\eps,\delta}$
and $\bar q{}'_{\eps,\delta}$ are the natural projections
\item
$j_\eps$ and $j_{\eps,\delta}$ are the open immersions.
\end{itemize}
Define 
\set\begin{equation}\label{eq_compleat-def}
\cG(F,\psi,\eps,\delta):=
q^*_{\eps,\delta} F\otimes_\Lambda\cL_\psi\La m_{\eps,\delta}\Ra
\end{equation}
where $m_{\eps,\delta}$ is the restriction of $m$ to
$\D(z,\eps)^*\times_S\D(\infty',\delta)^*$.
Mostly we shall drop the explicit mention of $F$ and $\psi$,
and simply write $\cG(\eps,\delta)$. With this notation, set
$$
\begin{array}{r@{\: := \:}l}
\cG_!(\eps,\delta) & j_{\eps,\delta!}\cG(\eps,\delta) \\
\cG_*(\eps,\delta) & Rj_{\eps,\delta*}\cG(\eps,\delta) \\
M^i(\eps,\delta) & \colim{H\subset\pi(\infty',\delta)}\,
H^i(\D(z,\eps)\times_SC_{H,\et},j_H^*\cG_!(\eps,\delta)).
\end{array}
$$
Notice that the $\Lambda$-module $M^i(\eps,\delta)$ carries
a natural continuous action of $\pi(\infty',\delta)$.
Moreover, if $\eps'\leq\eps$ and $\delta'\leq\delta$,
any choice of a group homomorphism
$\pi(\infty',\delta')\to\pi(\infty',\delta)$ as
in \eqref{subsec_smaller-radius}, induces a
$\pi(\infty',\delta')$-equivariant homomorphism of
$\Lambda$-modules :
$$
M^i(\eps,\delta)\to M^i(\eps',\delta').
$$
By inspecting the definitions, we find a natural isomorphism :
\set\begin{equation}\label{eq_represent-LFT}
\cF_\psi^{(z,\infty')}(F_{\eta_z})\isom
\colim{n\in\N}\,\colim{\delta>0}\,M^1(\eps_n,\delta).
\end{equation}

\sset\subsubsection{}\label{subsec_long-E}
For given $\eps,\delta\in|K^\times|$, denote by $Q_{\eps,\delta}$
the cone of the natural morphism
$\cG_!(\eps,\delta)\to\cG_*(\eps,\delta)$. Clearly the cohomology
sheaves $\cH^iQ_{\eps,\delta}$ are concentrated on
$\{z\}\times_S\D(\infty',\delta)^*$ for every $i\in\Z$,
and they agree on this closed subspace with the restriction
of $R^ij_{\eps,\delta*}\cG(\eps,\delta)$.

After applying the functor $R\bar q{}'_{\eps,\delta*}$ and
inspecting the resulting distinguished triangle, we obtain
the exact sequence :
$$
\cE\quad :\quad
q{}'_{\eps,\delta*}\cG(\eps,\delta)\to
\bar q{}'_{\eps,\delta*}(\cH^0Q_{\eps,\delta})\to
R^1\bar q{}'_{\eps,\delta*}\cG_!(\eps,\delta)\to
R^1q'_{\eps,\delta*}\cG(\eps,\delta)\to
\bar q{}'_{\eps,\delta*}(\cH^1Q_{\eps,\delta}).
$$
On the other hand, for every maximal point $a\in A'$ and every
$\eps\in|K^\times|$, let $w(a,\eps)\in A\times_SA'$ be the unique
point such that $(\E(z,\eps)\times\{a\})\setminus\{w(a,\eps)\}$
is the interior of $\E(z,\eps)\times\{a\}$ in $D\times_S\{a\}$
(notation of \eqref{eq_def-E_i}).
Set $w(\eps):=p(w(a,\eps))$ (which is independent of $a$), and let 
$$
\beta_{\eps,1},\dots,\beta_{\eps,r}\in\Gamma_{w(\eps)}\otimes_\Z\Q
$$
be the breaks of the stalk $F_{w(\eps)}$ (here $r$ is the
$\Lambda$-rank of $F_{w(\eps)}$, and each break is repeated
with multiplicity equal to the rank of the corresponding direct
summand in the break decomposition of this stalk). Denote also by
$\alpha(a,\eps)\in\Gamma_{w(a,\eps)}\otimes_\Z\Q$ the unique
break of the stalk $\cL_\psi\La m\Ra_{w(a,\eps)}$.
(Here $\Gamma_{w(\eps)}$ and $\Gamma_{w(a,\eps)}$ are the valuation
groups of the residue fields of $w(\eps)$, respectively $w(a,\eps)$,
and notice that the projection $p$ induces an injective homomorphism
of ordered groups $\Gamma_{w(\eps)}\to\Gamma_{w(a,\eps)}$; moreover
the break decompositions are invariant under extension from the
base field $(K,|\cdot|)$ to the residue field $(L,|\cdot|_L)$
of $a$ : see \cite[Lemma 3.3.8]{Ram}.)

\begin{lemma}\label{lem_analogous-to-claim}
For every $\eps\in|K^\times|$ with $\eps\leq\eps_0$ we may find
$\delta\in|K^\times|$ such that :
\begin{enumerate}
\item
$\beta_{\gamma,i}>\alpha(a,\gamma)$
\item
$\sharp_\Q\alpha(a,\gamma)=1$
\item
$\alpha(a,\gamma)$ is the unique break of the stalk
$\cG(\eps,\delta)_{w(a,\gamma)}$
\end{enumerate}
for every $\gamma\in|K^\times|\cap[\eps/2,\eps]$, every
$i=1,\dots,r$, and every maximal point $a\in\D(\infty',\delta)^*$.
\end{lemma}
\begin{proof} (i) and (ii) are analogous to claim
\ref{cl_lower-bound}, and they follow likewise from
\cite[Lemma 9.4]{Hu3}, together with the continuity properties
of the breaks (\cite[Lemma 4.2.12]{Ram}). Assertion (iii)
follows immediately from (i) (cp. the proof of theorem
\ref{th_cycles-vanish}).
\end{proof}

Lemma \ref{lem_analogous-to-claim}(ii),(iii) implies that
\set\begin{equation}\label{eq_deduce-that}
q'_{\eps,\delta*}\cG(\eps,\delta)=0
\qquad
\text{for $\eps$ and $\delta$ as in lemma \ref{lem_analogous-to-claim}.}
\end{equation}

\begin{lemma}\label{lem_no-name}
In the situation of \eqref{subsec_def-M}, suppose $\Lambda$ is a field,
$F$ is a locally constant $\Lambda$-module of finite type, with bounded
ramification at the point $z$, and let $\eps\in|K^\times|$ with
$\eps\leq\eps_0$. Then there exists $\delta\in|K^\times|$ such that :
\begin{enumerate}
\item
The $\Lambda$-module $R^iq'_{\eps,\delta!}\cG(\eps,\delta)$
vanishes for $i\neq 1$, and is locally constant of finite type
on $\D(\infty',\delta)^*_\et$ for $i=1$.
\item
The $\Lambda$-module $R^iq'_{\eps,\delta*}\cG(\eps,\delta)$
vanishes for $i\neq 1$, and is locally constant of finite type
on $\D(\infty',\delta)^*_\et$ for $i=1$.
\item
For $i\leq 1$, the $\Lambda$-module $\cH^iQ_{\eps,\delta}$ is
locally constant of finite type on
$\{z\}\times_S\D(\infty',\delta)^*_\et$.
\end{enumerate}
\end{lemma}
\begin{proof}(i):  The vanishing assertion for $i=0$ is clear.
Hence, suppose that $i>0$, and for every
$\gamma,\gamma'\in|K^\times|$ with $\gamma<\gamma'\leq\eps$
denote by $\Lambda(\gamma,\gamma')$ the $\Lambda$-module on the
\'etale site of $\D(z,\eps)^*\times_S\D(\infty',\delta)^*$ obtained
as extension by zero of the constant $\Lambda$-module with stalk
$\Lambda$ on the \'etale site of
$A(z,\gamma,\gamma')\times_S\D(\infty',\delta)^*$, where
$A(z,\gamma,\gamma'):=\D(z,\gamma')\setminus\E(z,\gamma)$
is the annulus centered at $z$ of external radius $\gamma'$
and internal radius $\gamma$. Set
$$
G^i_{\gamma,\gamma'}:=
R^iq'_{\eps,\delta!}
(\cG(\eps,\delta)\otimes_\Lambda\Lambda(\gamma,\gamma')).
$$
By \cite[Th.6.2.2]{Hu2}, the $\Lambda$-module $G^i_{\gamma,\gamma'}$
is constructible on $\D(\infty',\delta)^*_\et$, for every
$\gamma,\gamma'$ as above and every $i\in\N$. Moreover
\set\begin{equation}\label{eq_colimit-gamma}
R^iq'_{\eps,\delta!}\cG(\eps,\delta)=\colim{\gamma>0}\,
G^i_{\gamma,\eps}
\end{equation}
by \cite[Prop.5.4.5(i)]{Hu2}. Let $a\in\D(\infty',\delta)$
be any maximal point, and $(L,|\cdot|_L)$ the algebraic closure of
the residue field of $a$; the stalk $(G^i_{\gamma,\gamma'})_a$
calculates
$$
H^i_c(A(z,\gamma,\gamma')\times_S\Spa(L,L^+),\cG(\eps,\delta))
$$
(\cite[Th.5.4.6]{Hu2}). Pick $\delta\in|K^\times|$ fulfilling
conditions (i),(ii) of lemma \ref{lem_analogous-to-claim};
by Poincar\'e duality (lemma \ref{lem_PD-lambda}(iii)) and lemma
\ref{lem_analogous-to-claim}(ii),(iii) we deduce that
\set\begin{equation}\label{eq_vanish-kaboom}
(G^i_{\gamma,\gamma'})_a=0
\qquad
\text{for every $i\neq1$ and
      every $\gamma'\in|K^\times|\cap[\eps/2,\eps]$.}
\end{equation}
Since the maximal points are dense in the constructible topology
of $\D(\infty',\delta)^*$ (\cite[Cor.4.2]{Hu1}), we see that
$G^i_{\gamma,\eps}=0$ for $i>1$; combining with
\eqref{eq_colimit-gamma} we get the sought vanishing for $i>1$
and every sufficently small $\delta\in|K^\times|$.

To deal with the remaining case $i=1$, we remark :

\begin{claim}\label{cl_over-here}
There exists $\delta\in|K^\times|$ such that
$R^1q'_{\eps,\delta!}\cG(\eps,\delta)$ is overconvergent,
and all its stalks have the same finite length.
\end{claim}
\begin{pfclaim} Indeed, choose $\delta$  fulfilling
conditions (i),(ii) of lemma \ref{lem_analogous-to-claim},
so that \eqref{eq_vanish-kaboom} holds.
In view of \cite[Cor.10.3]{Hu3}, it follows that the natural
morphisms
\set\begin{equation}\label{eq_range}
G^1_{\gamma,\gamma'}\to G^1_{\gamma,\eps}
\end{equation}
induce isomorphisms on the stalks over the maximal points of
$\D(\infty',\delta)^*$, whenever $\gamma<\eps/2\leq\gamma'\leq\eps$.
Then again, since the maximal points are dense in the constructible
topology of $\D(\infty',\delta)^*$, we deduce that \eqref{eq_range}
restrict to isomorphisms on $\D(\infty',\delta)^*_\et$ for the stated
range of $\gamma,\gamma'$. Set
$$
G^1:=\colim{\gamma>0}\,\colim{\gamma'<\eps}\,G^1_{\gamma,\gamma'}.
$$
Combining with \eqref{eq_colimit-gamma} we see that the natural map :
$$
G^1\to R^1q'_{\eps,\delta!}\cG(\eps,\delta)
$$
is an isomorphism. However, notice that $G^1=R^iq'_{\eps,\delta!}
(\cG(\eps,\delta)\otimes_\Lambda\Lambda(0^+,\eps^-))$, where
$$
\Lambda(0^+,\eps^-)=
\colim{\gamma>0}\,\colim{\gamma'<\eps}\Lambda(\gamma,\gamma')
$$
which is an overconvergent sheaf, whose support is partially
proper over $\D(\infty',\delta)^*$; then $G^1$ is overconvergent
(\cite[Cor.8.2.4]{Hu2}), so the same holds for
$R^1q'_{\eps,\delta!}\cG(\eps,\delta)$.

In order to compute the length of the stalk over a given point
$a\in\D(\infty',\delta)^*$, we may therefore assume that $a$ is
maximal. Since the vanishing assertion (i) is already known
for $i\neq 1$, this length is completely determined -- in view of
\cite[Cor.10.3]{Hu3} and lemma \ref{lem_analogous-to-claim} --
by the Swan conductor at the point $(z,a)$ of the restriction
of $\cG(\eps,\delta)$ to the fibre $q^{\prime-1}_{\eps,\delta}(a)$.
Now, if $z\in A(K)$ it is clear that this Swan conductor equals
$\ssw^\natural_z(F,0^+)$, so it does not depend on $a$, whence
the claim. Lastly, if $z=\infty$, proposition \ref{prop_drop-is-finite}(i)
implies that the Swan conductor will also be independent of $a$,
provided $\D(\infty',\delta)^*\cap d_K(F)=\emptyset$, which can
be arranged by further shrinking $\delta$.
\end{pfclaim}

To achieve the proof of (i), it suffices now to invoke
claims \ref{cl_deduce-loc-constance} and \ref{cl_over-here},
together with identity \eqref{eq_colimit-gamma} for $i=1$.

(ii): Set
$\cG(\eps,\delta)^\vee:=\cHom_\Lambda(\cG(\eps,\delta),\Lambda_X)$
and notice that
$\cG(\eps,\delta)^\vee\simeq\cG(F^\vee,\psi^{-1},\eps,\delta)$
(notation of \eqref{eq_compleat-def}) hence assertion (i) holds
as well with $\cG(\eps,\delta)$ replaced by $\cG(\eps,\delta)^\vee$.
Moreover, since $\Lambda$ is a field, the natural morphism
$$
\cG(\eps,\delta)\to R\cHom_\Lambda(\cG(\eps,\delta)^\vee,\Lambda)
$$
is an isomorphism in
$\sD(\D(z,\eps)^*\times_S\D(\infty',\delta)^*_\et,\Lambda)$.
By Poincar\'e duality (lemma \ref{lem_PD-lambda}(ii)), there follows
a natural isomorphism
$$
Rq'_{\eps,\delta*}\cG(\eps,\delta)\isom
R\cHom_\Lambda(Rq'_{\eps,\delta!}\cG(\eps,\delta)^\vee,\Lambda(1)[2])
\isom\cHom_\Lambda(Rq'_{\eps,\delta!}\cG(\eps,\delta)^\vee,\Lambda(1)[2])
$$
whence the contention.

(iii): Suppose first that $z\in A(K)$; then, by smooth base change
(\cite[Th.4.1.1(a)]{Hu2}), it is easily seen that $\cH^iQ_{\eps,\delta}$
is locally isomorphic (in the \'etale topology) to the constant
$\Lambda$-module whose stalk is $(R^ij_{\eps*}F)_z$
(notation of \eqref{subsec_def-M}). It remains therefore only
to show that the latter is a $\Lambda$-module of finite type,
which is clear, since by assumption $F$ has bounded ramification
at the point $z$ (details left to the reader).

In case $z=\infty$, proposition \ref{prop_perverse-FT} implies that
actually $Q_{\eps,\delta}=0$ (details left to the reader).
\end{proof}

\begin{proposition}\label{prop_significant}
In the situation of \eqref{subsec_def-M}, suppose that $F$ is
locally constant of finite type, and with bounded ramification
at the point $z$. Then for every $\eps\in|K^\times|$ with
$\eps\leq\eps_0$ we may find $\delta_\eps\in|K^\times|$ such
that the following holds for every $\delta\leq\delta_\eps$ :
\begin{enumerate}
\item
$M^1(\eps,\delta)$ is a $\Lambda$-module of finite length.
\item
$M^i(\eps,\delta)=0$ for every $i\neq 1$.
\end{enumerate}
\end{proposition}
\begin{proof} By the usual arguments, we may reduce to the case
where $\Lambda$ is a field. We notice :

\begin{claim}\label{cl_we-notice}
In order to prove the proposition, it suffices to find
$\delta_\eps\in|K^\times|$, such that the following holds
for every $\delta\leq\delta_\eps$ :
\begin{enumerate}
\alphaenu
\item
The $\Lambda$-module
$R^1\bar q{}'_{\eps,\delta*}\cG_!(\eps,\delta)$ is locally
constant of finite type on $\D(\infty',\delta)^*_\et$.
\item
$R^i\bar q{}'_{\eps,\delta*}\cG_!(\eps,\delta)=0$ for every
$i\neq 1$.
\end{enumerate}
\end{claim}
\begin{pfclaim} Indeed, for every open subgroup
$H\subset\pi(\infty',\delta)$ let
$q'_H:\D(z,\eps)\times_SC_H\to C_H$ be the projection; given
$\delta_\eps$ as in the claim, the smooth base change theorem
of \cite[Th.4.1.1(a)]{Hu2} implies that
$R^iq'_{H*}j^*_H\cG_!(\eps,\delta)$ vanishes for $i\neq 1$,
and is a locally constant $\Lambda$-module of finite type on
$C_{H,\et}$, for $i=1$ and every $\delta\leq\delta_\eps$.
Then the Leray spectral sequence for the morphism $q'_H$
yields a natural (equivariant) isomorphism :
\set\begin{equation}\label{eq_also-for-later}
\begin{array}{r@{\:\isom\:}l}
M^i(\eps,\delta) & \colim{H\subset\pi(\infty',\delta)}\,
H^{i-1}(C_H,R^1q'_{H*}j^*_H\cG_!(\eps,\delta)) \\
& \colim{H\subset\pi(\infty',\delta)}\,
H^{i-1}(C_H,R^1\bar q{}'_{\eps,\delta*}\cG_!(\eps,\delta)).
\end{array}
\end{equation}
for $i=1$, the target of \eqref{eq_also-for-later} is a
$\Lambda$-module of finite length, since every $C_H$ is
connected, and $\Lambda$ is noetherian. For $i=0$, the
target of \eqref{eq_also-for-later} trivially vanishes.
For $i>1$, we remark that, (up to restricting the
colimit to a cofinal system of open subsets $H$, we may replace
$R^1\bar q{}'_{\eps,\delta*}\cG_!(\eps,\delta)$ by a constant
$\Lambda$-module; then, the sought vanishing for $i=2$ follows
from lemma \ref{lem_complete-later}. Lastly, for $i>2$, we may
argue by Poincar\'e duality, as in the proof of proposition
\ref{prop_vanish-vanishes}, to conclude.
\end{pfclaim}

Now, to begin with, it is clear that
$\bar q{}'_{\eps,\delta*}\cG_!(\eps,\delta)=0$
for every $\eps\leq\eps_0$ and every $\delta\in|K^\times|$.
Next, in view of \eqref{eq_deduce-that}, the exact sequence
$\cE$ of \eqref{subsec_long-E}, and lemma \ref{lem_no-name}(ii,iii),
we see that condition (a) of claim \ref{cl_we-notice} holds.
Similarly, for $i>1$, we consider the exact sequence 
$$
0=R^{i-1}\bar q{}'_{\eps,\delta*}Q_{\eps,\delta}\to
R^i\bar q{}'_{\eps,\delta*}\cG_!(\eps,\delta)\to
R^iq'_{\eps,\delta*}\cG(\eps,\delta)
$$
whose last term vanishes, according to lemma \ref{lem_no-name}(ii);
so also condition (b) of claim \ref{cl_we-notice} holds, and the
proposition follows.
\end{proof}

\sset\subsubsection{}\label{subsec_see-here}
The significance of proposition \ref{prop_significant}(i) is that
it allows to endow the local Fourier transform of $F_{\eta_z}$
with an action of $\pi(\infty',\delta)$, for some sufficiently
small $\delta\in|K^\times|$. Indeed, pick a monotonically
descending sequence $(\delta_n~|~n\in\N)$ of values in
$|K^\times|$ such that $M(\eps_n,\delta_n)$ is a $\Lambda$-module
of finite length for every $n\in\N$.
In view of \eqref{eq_represent-LFT} the $\Lambda$-module
$\cF_\psi^{(z,\infty')}(F_{\eta_z})$ is the colimit of the
filtered system $\cM$ of such $M(\eps_n,\delta_n)$. However,
proposition \ref{prop_Garuti}, together with remark
\ref{rem_not-much}(i) implies that the local Fourier transform
has finite length as a $\Lambda$-module, hence some
$M(\eps_n,\delta_n)$ must surject onto it. Moreover, the image
of $M(\eps_n,\delta_n)$ into $\cF_\psi^{(z,\infty')}(F_{\eta_z})$
is already isomorphic to the image of the same module into
$M(\eps_m,\delta_m)$ for some $m>n$. Since the transition maps
in the filtered system $\cM$ are equivariant, we may in this
way endow $\cF_\psi^{(z,\infty')}(F_{\eta_z})$ with an action of
$\pi(\infty',\delta_m)$. Let $F\to F'$ be a morphism of locally
constant $\Lambda$-modules on $\D(z,\eps_0)^*$ with bounded
ramification at $z$; by repeating the above procedure on $F'$,
we endow $\cF_\psi^{(z,\infty')}(F'_{\eta_z})$ with a
$\pi(\infty',\delta_{m'})$-action for some other $m'\in\N$.
A direct inspection of the construction shows that the induced map
$\cF_\psi^{(z,\infty')}(F_{\eta_z})\to
\cF_\psi^{(z,\infty')}(F'_{\eta_z})$ shall be
$\pi(\infty',\delta_{m''})$-equivariant for some $m''\geq m,m'$.

Now, for any $z\in D(K)$ (and any $z\in D'(K)$), define the
category of {\em $\Lambda[\pi(z)]$-modules with bounded ramification} :
$$
\Lambda[\pi(z)]\bMod
$$
as the full subcategory of $\Lambda[\pi(z)]\Mod$ whose
objects are represented by $\Lambda$-modules on the \'etale site
of $\D(z,\eps)^*$ (for some $\eps\in|K^\times|$) which have bounded
ramification at the point $z$. With this notation, the foregoing
means that the datum of $\cF_\psi^{(z,\infty')}(F_{\eta_z})$ together
with its $\pi(\infty',\delta_m)$-action, is well defined as a functor
$$
\cF^{(z,\infty')}_\psi:
\Lambda[\pi(z)]\bMod\to\Lambda[\pi(\infty')]\bMod
$$
which we call the {\em local Fourier transform from $z$ to $\infty'$}.

\sset\subsubsection{}\label{subsec_furth-signify}
Furthermore, \eqref{eq_also-for-later} implies that the
natural decomposition \eqref{eq_stat-phase} is already
well defined in $\Lambda[\pi(\infty')]\bMod$. Indeed, a direct
inspection shows that the composition of \eqref{eq_stat-phase}
with the projection onto the direct factor indexed by $z$,
is obtained as the colimit of the filtered system of maps
$$
\Gamma(C_H,\cH^1(\cF_\psi(j_!F)))\to
\Gamma(C_H,R^1\bar q{}'_{\eps,\delta*}\cG_!(\eps,\delta))
$$
induced by the natural morphism
$$
R^1p'_*(p^*(j_!F)\otimes_\Lambda\cL_\psi\La m\Ra)_{|\D(\infty',\delta)}
\to R^1\bar q{}'_{\eps,\delta*}\cG_!(\eps,\delta)
$$
which is clearly $\pi(\infty',\delta)$-equivariant.

\begin{remark}\label{rem_dual-phase}
We will also use the ``dual'' functors $\cF_\psi^{(z,\infty)}$,
for $z\in D'(K)$ : they are defined by exchanging everywhere the
roles of $D$ and $D'$ in the foregoing (cp. \cite[D\'ef.2.4.2.3]{Lau2}).
Especially, given a Zariski open immersion $j':U'\to A'$, and
a local system $F'$ on $U'_\et$ with bounded ramification everywhere,
one has a similar stationary phase decomposition for
$b_*\cF'_\psi(j'_*F')_{\eta_\infty}$, in terms of such functors
(here $b:A''\to A$ is the double duality isomorphism). We leave
to the reader the task of spelling out this isomorphism.
\end{remark}

\section{Local analysis of monodromy}\label{chap_monodro}
\subsection{Study of the local Fourier transforms}\label{sec_study}
To begin with, we notice :

\begin{proposition}\label{prop_local-FT-exact}
For every $z\in D(K)$, the local Fourier transform
$\cF_\psi^{(z,\infty')}$ is an exact functor on
$\Lambda[\pi(z)]\bMod$.
\end{proposition}
\begin{proof} Choose a sequence $(\eps_n~|~n\in\N)$ as in
\eqref{subsec_pass-to-lim}. For every continuous
$\Lambda[\pi(z,\eps)]$-module $M$ of finite type, set
$$
T^i(M):=\colim{n\in\N}i_z^*R^i\Psi_{\eta_{\infty'},\eps_n}
(\bar p{}^*_\eps(M_!)\otimes_\Lambda\bar\cL_\psi\La m\Ra_{|X(\eps)})
$$
(notation of \eqref{subsec_as-expla}). It suffices to check :

\begin{claim} Suppose that $M_!$ has bounded ramification at
the point $z$. Then $T^i(M)=0$ for every $i\neq 1$.
\end{claim}
\begin{pfclaim}[] The assertion follows easily from proposition
\ref{prop_significant}(ii).
\end{pfclaim}
\end{proof}

\sset\subsubsection{}
For any $z\in K$, let $\theta_z:A\to A$ the translation by $z$ :
$x\mapsto x+z$ (notation of \eqref{subsec_intro-Fourier}).
Then $\theta_z$ induces an isomorphism of categories :
$$
\theta_{z*}:\Lambda[\pi(0)]\Mod\isom\Lambda[\pi(z)]\Mod
$$
in an obvious way. With this notation, an easy calculation
yields a natural isomorphism :
\set\begin{equation}\label{eq_translate-locFT}
\cF^{(z,\infty')}_\psi(\theta_{z*}M)\isom
\cF_\psi^{(0,\infty')}(M)\otimes_\Lambda
\cL_\psi\La zx'\Ra_{\eta_{\infty'}}
\end{equation}
for all $\Lambda[\pi(0)]$-modules $M$ with bounded ramification,
where $\cL_\psi\La zx'\Ra_{\eta_{\infty'}}$ is the object of
$\Lambda[\pi(\infty')]\bMod$ defined as in
\eqref{subsec_more-trivialities} (and remark
\ref{rem_one-more-obviety}). It follows that the study of
$\cF^{(z,\infty')}_\psi$ for $z\in K$ is reduced to the case
where $z=0$.

Notice as well that a $\Lambda[\pi(z)]$-module $M$ with
bounded ramification has a well defined length and Swan conductor :
$$
\length_\Lambda M
\qquad
\ssw(M)
$$
namely, the length (resp. the Swan conductor at $z$) of any
$\Lambda[\pi(z,\eps)]$-module representing $M$.

Let $U\subset A$ be a Zariski open subset, $j:U\to A$ the open
immersion, and $F$ a Zariski constructible $\Lambda$-module on
$A_\et$, such that $j^*F$ is locally constant on $U_\et$, with
bounded ramification at all the points of $D\setminus U$.
Denote by $\lambda_1,\dots,\lambda_l$ the breaks of $F_{\eta_\infty}$
(repeated with their respective multiplicities). Then, in view
of lemma \ref{lem_Groth-Ogg-Shaf}, theorem \ref{th_perversity}
and remark \ref{rem_perversity}, we deduce the following
identity :
\set\begin{equation}\label{eq_gen-length-FT}
\length_\Lambda\cH^0(\cF_\psi(F))_{\eta_{\infty'}}=
\sum_{x\in A\setminus U}a_x(F)+
\sum_{i=1}^l\max(0,\lambda^\natural_i-1).
\end{equation}
(Details left to the reader.)

\sset\subsubsection{}\label{subsec-tame}
Hence, let $z\in\{0,\infty\}$.
We consider first the case of a {\em tamely ramified} object of
$\Lambda[\pi(z)]\Mod$. Namely, for a given value
$\eps\in|K^\times|$, let $T$ be a coordinate on $\D(z,\eps)$,
such that $T(z)=0$. For every integer $n>0$, the rule $T\mapsto T^n$
determines a finite morphism $\D(z,\eps^{1/n})\to\D(z,\eps)$;
its restriction to $\D(z,\eps^{1/n})^*$ is a finite connected
Galois \'etale covering of $\D(z,\eps)^*$, whose Galois group is
naturally isomorphic to $\bmu_n$, the $n$-torsion in $K^\times$.
This covering corresponds, as in \eqref{subsec_G-sets}, to a
surjective continuous group homomorphism $\pi(z,\eps)\to\bmu_n$.
For variable $n>0$, we obtain an inverse system of group
homomorphisms, whence a continuous group homomorphism :
\set\begin{equation}\label{eq_tame-char}
\pi(z,\eps)\to\hat\Z(1):=\liminv{n>0}\,\bmu_n
\end{equation}
(where the target group is endowed with the profinite topology).
The map \eqref{eq_tame-char} is surjective; to see this, we may
suppose that the image $v$ of the geometric point $\xi$ chosen
in \eqref{sec_vanish} corresponds to the Gauss valuation
``at the border'' of $\D(z,\eps)$, {\em i.e.} :
$$
|f(v)|_v=
\max(|f(a)|\text{ for all $a\in K$ such that $|a-z|\leq\eps$})
$$
for every power series $f(T)$ convergent on $\D(z,\eps)$.
Let $\kappa(v)$ denote the residue field of $v$; then
$\pi_1^\localg(v,\xi)$ is the Galois group of
$\kappa(v)$-automorphisms of any algebraic closure of
$\kappa(v)$. We have a natural group homomorphism :
$$
\pi_1^\localg(v,\xi)\to\pi(z,\eps)
$$
and it suffices to see that its composition with
\eqref{eq_tame-char} is surjective. The latter is a continuous
map of profinite groups, so we reduce to showing that the
induced maps $\pi_1^\localg(v,\xi)\to\bmu_n$
are surjective for every $n>0$, which is left to the reader
({\em e.g.} one may look at a rank two specialization of
$v$, and argue by standard valuation theory).

\sset\subsubsection{}\label{subsec_tamely}
With the notation of \eqref{subsec-tame}, a
$\Lambda[\pi(z,\eps)]$-module $M$ is {\em tamely ramified\/}
relative to the coordinate $T$, if the action of $\pi(z,\eps)$
on $M$ factors through \eqref{eq_tame-char}.
We say that a $\Lambda[\pi(z)]$-module $M$ is
{\em tamely ramified} if it is represented by a tamely ramified
$\Lambda[\pi(z,\eps)]$-module, for some $\eps\in|K^\times|$.
Notice that this notion is independent of the choice of
coordinate $T$. For such a module $M$, we may find $n>0$
large enough, so that the action of $\pi(z,\eps)$ factors
through $\bmu_n$. If the action of $\pi(z,\eps)$ on $M$
is trivial, we say that $M$ is {\em unramified}.

We may study the local Fourier transform $\cF^{(z,\infty)}_\psi(M)$
by a global argument, as in \cite{Lau2}.
Namely, suppose now that $T$ is a global coordinate on $\A^1_K$;
then the rule $T\mapsto T^n$ defines the {\em Kummer covering} of
$\G_m:=\A^1_K\setminus\{0\}$ with Galois group $\bmu_n$. Let
$\chi:\bmu_n\to\Lambda^\times$ be any non-trivial character;
we denote by $\cK_\chi$ the locally constant $\Lambda$-module
of rank one on $(\G_m)^\ad_\et$ associated to this Kummer covering
and the character $\chi$.
If $f:X\to(\G_m)^\ad_\et$ is any morphism of $S$-adic spaces,
we let as usual $\cK_\chi\La f\Ra:=f^*\cK_\chi$.
Especially, the global coordinate $x$ on $A$ (notation of
\eqref{subsec_intro-Fourier}) yields a morphism
$x:A\setminus\{0\}\to(\G_m)^\ad$, whence the local system
$\cK_\chi\La x\Ra$ on $A\setminus\{0\}$. Likewise we have
the $\Lambda$-module $\cK_\chi\La x'\Ra$ on $A'\setminus\{0'\}$.
Let $j:A\setminus\{0\}\to A$ and $j':A'\setminus\{0\}\to A'$ be
the open immersions; the Fourier transform
of $j_*\cK_\psi\La x\Ra=j_!\cK_\chi\La x\Ra$ is calculated by
the following :

\begin{lemma}\label{lem_tame-FT}
There is a natural isomorphism of $\Lambda$-modules :
$$
\cF_\psi(j_*\cK_\chi\La x\Ra)\isom
j'_*\cK_{\chi^{-1}}\La x'\Ra[0]\otimes_\Lambda G(\chi,\psi)
$$
where
$G(\chi,\psi):=H^1_c((\G_m)^\ad_\et,\cK_\chi\otimes_\Lambda\cL_\psi)$
is a free $\Lambda$-module of rank one.
\end{lemma}
\begin{proof} {\em Mutatis mutandi}, the proof of
\cite[Prop.1.4.3.2]{Lau2} can be taken over.
\end{proof}

As a corollary we obtain :

\begin{proposition}\label{prop_tame-FT}
Let $M$ be a tamely ramified $\Lambda[\pi(z)]$-module. We have :
\begin{enumerate}
\item
If $z=0$, then
$$
\length_\Lambda\cF^{(z,\infty')}_\psi(M)=\length_\Lambda M
$$
and the $\Lambda[\pi(\infty')]$-module $\cF^{(z,\infty')}_\psi(M)$
admits a filtration
$$
0=F_0\subset F_1\subset\cdots\subset F_l=\cF^{(z,\infty')}_\psi(M)
$$
whose subquotients are tamely ramified $\Lambda[\pi(\infty')]$-modules.
\item
If $z=\infty$, then $\cF^{(z,\infty')}_\psi(M)=0$.
\end{enumerate}
\end{proposition}
\begin{proof} The module $M$ corresponds to a representation of $\bmu_n$
for some $n>0$; we may then easily reduce to the case of an irreducible
representation of $\bmu_n$, {\em i.e.} $\Lambda$ is a field, and
$\dim_\Lambda M=1$ (see \cite[Partie II, \S2.5, 2.6]{Se2}),
so $M$ is given by a character $\chi:\bmu_n\to\Lambda^\times$.
For the case of the trivial character $\chi$, we have the more
precise :

\begin{claim}\label{cl_verbatim}
For every unramified $\Lambda[\pi(0)]$-module, there is a
natural isomorphism :
$$
\cF_\psi^{(0,\infty')}(M)\isom M.
$$
\end{claim}
\begin{pfclaim} We may take over {\em verbatim\/} the proof of
\cite[Prop.2.5.3.1(i)]{Lau2}.
\end{pfclaim}

In case $\chi$ is non-trivial and $z=0$ (resp. $z=\infty$), $M$
is represented by $\cK_\chi\La x\Ra_{\eta_0}$
(resp. by $\cK_{\chi^{-1}}\La x\Ra_{\eta_\infty}$).
Using lemma \ref{lem_tame-FT} and the stationary phase argument
of \cite[Prop.2.5.3.1(ii)]{Lau2}, we deduce a natural isomorphism :
$$
\cF^{(0,\infty')}_\psi(M)\oplus
\cF^{(\infty,\infty')}_\psi(\cK_{\chi^{-1}}\La x\Ra_{\eta_\infty})
\isom
\cK_{\chi^{-1}}\La x'\Ra_{\eta_{\infty'}}\otimes_\Lambda G(\chi,\psi).
$$
By comparing ranks, we see that
\set\begin{equation}\label{eq_bound-tame}
\dim_\Lambda\cF^{(0,\infty')}_\psi(M)\leq 1
\end{equation}
and both assertions will follow, once we show that
$\cF^{(0,\infty')}_\psi(M)$ does not vanish. To this aim, let
$\bar\cK_{\!\!\chi}\La x\Ra$ denote the $\Lambda$-module
on $D_\et$ which is the extension by zero of
$\cK_\chi\La x\Ra$. The rule $x\mapsto x/(x+1)$ defines
an automorphism $\omega$ of $D$ such that $\omega(0)=0$,
$\omega(\infty)=1$ and $\omega(-1)=\infty$. Let $F$ denote the
restriction to $A_\et$ of $\omega^*\bar\cK_{\!\!\chi}\La x\Ra$.
Clearly $F_{|\D(\infty,\eps)^*}$ corresponds to the trivial
character of $\pi(\infty,\eps)$ (for $\eps$ small enough);
also $F_{\eta_0}$ represents $M$, and $F_{\eta_{-1}}$
is a tamely ramified $\Lambda[\pi(-1)]$-module of
length one. We apply theorem \ref{th_Stationary-phase} : by the
foregoing, the term $\cF_\psi^{(\infty,\infty')}(F_{\eta_\infty})$
vanishes, whence a natural isomorphism :
$$
\cH^0(\cF_\psi(F))_{\eta_{\infty'}}\isom
\cF^{(0,\infty')}_\psi(M)\oplus\cF_\psi^{(-1,\infty')}(F_{\eta_{-1}}).
$$
Computing with \eqref{eq_gen-length-FT} we find that
$\cH^0(\cF_\psi(F))_{\eta_{\infty'}}$ has length equal to 2.
Taking \eqref{eq_bound-tame} and \eqref{eq_translate-locFT}
into account, we conclude that both direct summands on the
right hand-side must have length 1, as required.
\end{proof}

\begin{corollary}\label{cor_infty-infty}
Let $M$ be a $\Lambda[\pi(\infty)]$-module with bounded
ramification, and $\lambda_1,\dots,\lambda_l$ the breaks of $M$.
Suppose that $\lambda_i^\natural\leq 1$ for every $i=1,\dots,l$.
Then $\cF^{(\infty,\infty')}_\psi(M)=0$.
\end{corollary}
\begin{proof} We may extend $M$ to a locally constant
$\Lambda$-module $F$ on a Zariski open subset $U\subset A$,
tamely ramified at the points of $\Sigma:=A\setminus U$, and with
$F_{\eta_\infty}=M$ (proposition \ref{prop_Garuti}).
Let $j:U\to A$ be the open immersion, and set
$F':=\cH^0(\cF_\psi(j_!F))$; according to
\eqref{eq_gen-length-FT} we have
$$
\length_\Lambda F'_{\eta_{\infty'}}=\sharp\Sigma\cdot l
$$
where $\sharp\Sigma$ denotes the cardinality of $\Sigma$, and $l$
is the generic length of $F$. On the other hand,
\eqref{eq_stat-phase}, \eqref{eq_translate-locFT} and
proposition \ref{prop_tame-FT}(i) imply :
$$
\length_\Lambda F'_{\eta_{\infty'}}=\sharp\Sigma\cdot l+
\length_\Lambda\cF^{(\infty,\infty')}_\psi(M)
$$
whence the contention.
\end{proof}

\begin{theorem}\label{th_analyze-breaks}
Let $M$ be any $\Lambda[\pi(0)]$-module with bounded ramification.
Then we have :
\begin{enumerate}
\item
$\length_\Lambda\cF_\psi^{(0,\infty')}(M)=l:=\length_\Lambda M+\ssw(M)$.
\item
$\ssw(\cF_\psi^{(0,\infty')}(M))=\ssw(M)$.
\item
Let $\gamma_1,\dots,\gamma_l$ be the breaks of
$\cF_\psi^{(0,\infty')}(M)$ (repeated with their respective
multiplicities). Then $\gamma_i^\natural\leq 1$ for $i=1,\dots,l$.
\end{enumerate}
\end{theorem}
\begin{proof} We shall borrow an argument from the proof of
\cite[Prop.8.6.2]{Kat2}. By the usual reductions, we may assume
that $\Lambda$ is a field. Moreover, we may assume that $M$ does
not admit any quotient which is a {\em trivial}
$\Lambda[\pi(0)]$-module (meaning, a module represented
by a constant $\Lambda$-module on $\D(0,\eps)^*_\et$, for
some $\eps\in|K^\times|$). Indeed, for such trivial modules,
the theorem is already known, in view of proposition
\ref{prop_tame-FT}(i).

By proposition \ref{prop_Garuti}, we may
find a Zariski open subset $U\subset A\setminus\{0\}$ and a locally
constant $\Lambda$-module $F$ on $U_\et$, with bounded ramification
at every point of $\Sigma:=D\setminus U$, such that $F$ is tamely
ramified at every point of $\Sigma\setminus\{0\}$, and such that
$F_{\eta_0}=M$. Set $G:=\cH^0(\cF_\psi(j_!F))$, where $j:U\to A$ is
the open immersion; from \eqref{eq_gen-length-FT} we deduce easily :
\set\begin{equation}\label{eq_one-punch}
\length_\Lambda G_{\eta_{\infty'}}=
(\sharp\Sigma-1)\cdot\length_\Lambda M+\ssw(M)
\end{equation}
(where $\sharp\Sigma$ denotes the cardinality of $\Sigma$).
On the other hand, theorem \ref{th_Stationary-phase}, together
with proposition \ref{prop_tame-FT}(ii) yields a natural
isomorphism :
\set\begin{equation}\label{eq_punch-three}
G_{\eta_{\infty'}}\isom
\bigoplus_{z\in\Sigma\setminus\{\infty\}}\cF_\psi^{(z,\infty')}(F_{\eta_z})
\end{equation}
which, by virtue of proposition \ref{prop_tame-FT}(i), leads
to the identity :
\set\begin{equation}\label{eq_two-punch}
\length_\Lambda G_{\eta_{\infty'}}=(\sharp\Sigma-2)\cdot\length_\Lambda M+
\length_\Lambda\cF_\psi^{(0,\infty')}(M).
\end{equation}
Assertion (i) follows by comparing \eqref{eq_one-punch} and
\eqref{eq_two-punch}. Next we remark :

\begin{claim}\label{cl_Fourier-type}
$\cF_\psi(j_!F)=G[0]$.
\end{claim}
\begin{pfclaim} Indeed, by remark \ref{rem_perversity}(ii),(iii)
we know already that $\cH^i(\cF_\psi(j_!F))=0$ for $i\notin\{0,1\}$.
Moreover, in light of theorem \ref{th_perversity}(i) (and again
remark \ref{rem_perversity}(ii)), in order to show the same
vanishing for $i=1$, it suffices to verify that
$H^2_c(U_\et,F\otimes_\Lambda\cL_\psi\La ax\Ra)=0$ for every
$a\in K$ (where, as usual, $x$ is our fixed coordinate on $A$).
Since $\Lambda$ is a field, lemma \ref{lem_PD-lambda}(iii)
reduces to checking that
$H^0(U_\et,F^\vee\otimes_\Lambda\cL_\psi\La ax\Ra)=0$
for every $a\in K$, where $F^\vee:=\cHom_\Lambda(F,\Lambda_U)$.
However, since $M$ does not admit trivial quotients, the
$\Lambda[\pi(0)]$-module
$F^\vee_{\eta_0}=(F^\vee\otimes_\Lambda\cL_\psi\La ax\Ra)_{\eta_0}$
does not contain trivial submodules, whence the claim.
\end{pfclaim}

Notice as well that $d_K(F)=\{0\}$ (notation of \eqref{subsec_drop}),
hence $G$ is locally constant on $A'\setminus\{0'\}$.
From claim \ref{cl_Fourier-type} and \eqref{eq_FT-is-selfdual}
we get a natural isomorphism :
$$
\cF'_\psi G\isom a_*j_!F(-1)
$$
where $a:A\to A''$ is $(-1)$-times the double duality isomorphism,
and $(-1)$ denotes the Tate twist. Then we may apply
\eqref{eq_gen-length-FT} to derive the identity :
\set\begin{equation}\label{eq_with-a_0}
\length_\Lambda M=a_0(G)+\sum_{i=1}^l\max(0,\gamma^\natural_i-1).
\end{equation}
Let $\bar 0{}'$ be a geometric point localized at $0'\in A'$;
from claim \ref{cl_Fourier-type} we deduce as well that :
\set\begin{equation}\label{eq_will-need-later}
\length_\Lambda G_{\bar 0{}'}=-\chi_c(A,F)=
(\sharp\Sigma-2)\cdot\length_\Lambda M+\ssw(M).
\end{equation}
Taking \eqref{eq_one-punch} into account, we get :
$$
a_{0'}(G)=\length_\Lambda M+\ssw(G_{\eta_{0'}}).
$$
Comparing with \eqref{eq_with-a_0} we conclude that assertion
(iii) holds, and also that 
\set\begin{equation}\label{eq_extra-bonus}
a_{0'}(G)=\length_\Lambda M
\qquad
\ssw(G_{\eta_{0'}})=0.
\end{equation}
Next, from \eqref{eq_translate-locFT} and proposition
\ref{prop_tame-FT}(i) we deduce that
$\length_\Lambda\cF_\psi^{(z,\infty)}(F_{\eta_z})=\length_\Lambda M$
for every $z\in\Sigma\setminus\{0,\infty\}$. Combining with
\eqref{eq_punch-three} we obtain :
\set\begin{equation}\label{eq_but-you}
\ssw(\cF_\psi^{(0,\infty')}(M))=\ssw(G_{\eta_{\infty'}})-
(\sharp\Sigma-2)\cdot\length_\Lambda M.
\end{equation}
On the other hand, \eqref{eq_extra-bonus} and claim
\ref{cl_Fourier-type} also imply that
$$
0=\chi_c(A',G)=-2\cdot\length_\Lambda G_{\eta_{\infty'}}+
a_{0'}(G)+a_{\infty'}(G)=\length_\Lambda M+\ssw(G_{\eta_{\infty'}})-
\length_\Lambda G_{\eta_{\infty'}}
$$
which, in view of \eqref{eq_one-punch} leads to the identity :
\set\begin{equation}\label{eq_nothing-comp}
\ssw(G_{\eta_{\infty'}})=(\sharp\Sigma-2)\cdot\length_\Lambda M+
\ssw(M).
\end{equation}
Assertion (ii) follows by comparing \eqref{eq_but-you} and
\eqref{eq_nothing-comp}.
\end{proof}

\sset\subsubsection{}
Following \cite[\S2.4]{Lau2}, we wish now to exhibit a left
quasi-inverse for the functor $\cF^{(0,\infty')}_\psi$. Namely,
for any $\eps\in|K^\times|$, set
$$
Y(\eps):=\D(0,\eps)\times_SD'\subset D\times_SD'
$$
and denote by $\bar p_\eps:Y(\eps)\to\D(0,\eps)$,
$\bar p'_\eps:Y(\eps)\to D'$ the projections.
Also, for every $z'\in D'(K)$ and every object $M'$ of
$\Lambda[\pi(z',\eps)]\Mod_\mathrm{f.cont}$, let
$i_{z'}:S\to\bar q^{\prime-1}_\eps(0)$ be the closed immersion
with image $(0,z')$, and denote by $M'_!$ the
$\Lambda$-module on $D'_\et$ which is the extension by zero of
$\tau_{z'}M'$, where $\tau_{z'}$ is the equivalence between
continuous $\Lambda[\pi(z',\eps)]$-modules of finite length
and locally constant $\Lambda$-modules of finite length on
$\D(z',\eps)^*_\et$, as in \eqref{subsec_as-expla}. We define :
$$
\cF^{(z',0)}_\psi(M'):=i^*_{z'}R^1\Phi_{\eta_0}
(\bar p^{\prime*}_\eps(M'_!)\otimes_\Lambda
\bar\cL_\psi\La m\Ra_{|Y(\eps)})
$$
which clearly depends only on the image of $M'$ in
$\Lambda[\pi(z')]\Mod$. In fact, we will see that this definition
is interesting only for $z'=\infty'$, since this functor
vanishes identically for $z'\in A'(K)$. As a first step,
a simple inspection shows that, for $z'\in A'(K)$, the
translation map $\theta_{z'}:A'\to A'$ given by $x'\mapsto x'+z'$
induces an isomorphism :
\set\begin{equation}\label{eq_first-step}
\cF^{(z',0)}_\psi(\theta_{z*}M')\isom\cF_\psi^{(0',0)}(M')
\qquad
\text{for all $\Lambda[\pi(0')]$-modules $M'$}
\end{equation}
(analogous to \eqref{eq_translate-locFT}).
Next, we remark that these functors obey as well a principle
of stationary phase. Indeed, let $F'$ be a $\Lambda$-module on
$D'_\et$, and $j':U'\to D'$ an open immersion such that
$j^{\prime*}F'$ is a locally constant $\Lambda$-module of
finite length on $U'_\et$. Then corollary \ref{cor_vanish-smooth}
implies that :
\set\begin{equation}\label{eq_fake-station}
j^{\prime*}R\Phi_{\eta_0}(\bar p^{\prime*}_\eps F'\otimes_\Lambda
\bar\cL_\psi\La m\Ra_{|Y(\eps)})=0
\qquad
\text{in\ \  $\sD(U'_\et,\Lambda)$.}
\end{equation}

\sset\subsubsection{}\label{subsec_split-longish}
Suppose now $U'\subset A'$ is a Zariski open subset, and that
$F'=j'_!j^{\prime*}F'$. Let $b:A''\isom A$ be the double duality
isomorphism, given by the rule $x''\mapsto x$ (notation of
\eqref{subsec_intro-Fourier}), and set
$F_\eps:=b_*\cF'_\psi(F')_{|\D(0,\eps)}$, which is a complex
of $\Lambda$-modules on $\D(0,\eps)_\et$. Combining
\eqref{eq_fake-station} with proposition
\ref{prop_vanish-properties}(ii) and corollary
\ref{cor_vanish-properties}, we deduce a natural isomorphism :
\set\begin{equation}\label{eq_fake-stationary}
R^0\Phi_{\eta_0}(F_\eps)\isom
\bigoplus_{z'\in D'\setminus U'}\cF^{(z',0)}_\psi(F'_{\eta_{z'}}).
\end{equation}

\begin{remark}\label{rem_not-much-again}
We have as well a localized version of \eqref{eq_fake-stationary},
analogous to remark \ref{rem_not-much}(ii). Namely, let
$\delta\in|K^\times|$, and suppose that $F'$ is a locally constant
$\Lambda$-module of finite type on $\D(\infty',\delta)^*_\et$.
Let $w'\in\bar p{}^{-1}_\eps(0)$ be the unique point such that
$(\{0\}\times\D(\infty',\delta))\cup\{w'\}$ is the topological
closure of $\{0\}\times\D(\infty',\delta)$, and set
$$
\cF_\psi^{(w',0)}(F'):=
\Gamma(\{w'\},R^1\Phi_{\eta_0}
(\bar p{}^{\prime*}_\eps(\bar\jmath_!F')\otimes_\Lambda
\bar\cL_\psi\La m\Ra_{|Y(\eps)}))
$$
where $\bar\jmath:\D(\infty',\delta)^*\to D'$ is the open immersion.
Let also $F_\eps$ be the $\Lambda$-module on $\D(0,\eps)_\et$
that is the restriction of $b_*\cF'_\psi((\bar\jmath_!F')_{|A'})$;
arguing as in \eqref{subsec_split-longish} we get a natural isomorphism :
$$
R^0\Phi_{\eta_0}(F_\eps)\isom
\cF^{(\infty',0)}_\psi(F'_{\eta_{\infty'}})\oplus
\cF_\psi^{(w',0)}(F').
$$
\end{remark}

\sset\subsubsection{}
We are now sufficiently equipped to make some basic computations.
To begin with, in the situation of \eqref{subsec_split-longish},
theorem \ref{th_perversity}(iv) and remark \ref{rem_perversity}(iii)
imply that the sequence of $\Lambda$-modules :
\set\begin{equation}\label{eq_shortish}
0\to\cH^0(F_\eps)_0\to\cH^0(F_\eps)_{\eta_0}\to
R^0\Phi_{\eta_0}(F_\eps)\to\cH^1(F_\eps)_0\to 0
\end{equation}
is exact.  Let us take for $F'$ the extension by zero of the
constant $\Lambda$-module $\Lambda_{A'}$ on $A'_\et$. In this case,
$F_\eps=i_{0*}\Lambda_S(-1)[-1]$, where $i_0:S\to\D(0,\eps)$
is the closed immersion with image $0$. From \eqref{eq_shortish}
and \eqref{eq_fake-stationary} we deduce that :
\set\begin{equation}\label{eq_triv-equipp}
\length_\Lambda\cF^{(\infty',0)}_\psi(\Lambda)=
\length_\Lambda R^0\Phi_{\eta_0}(i_{0*}\Lambda_S[-1])=
\length_\Lambda\cH^0(i_{0*}\Lambda_S)_0=l
\end{equation}
with $l:=\length_\Lambda\Lambda$. Next we have :

\begin{lemma}\label{lem_fake-vanish}
Let $M'$ be a tamely ramified $\Lambda[\pi(z')]$-module. We have :
\begin{enumerate}
\item
If $z'\in A'(K)$, then $\cF^{(z',0)}_\psi(M')=0$.
\item
If $z'=\infty'$, then
$\length_\Lambda\cF^{(z',0)}_\psi(M')=\length_\Lambda M'$.
\end{enumerate}
\end{lemma}
\begin{proof} As usual, we may reduce to the case where
$\Lambda$ is a field, and we may also assume that $M'$
arises from a character $\chi:\bmu_n:\to\Lambda^\times$, for
some $n>0$. 

(i): In view of \eqref{eq_first-step} we may assume that
$z'=0'$. Consider first the case where $\chi$ is not the
trivial character; then $M'$ is represented by
$\cK_\chi\La x'\Ra_{\eta_0}$. Denote by
$\bar\cK_\chi\La x'\Ra$ the $\Lambda$-module on $D'_\et$
which is the extension by zero of $\cK_\chi\La x'\Ra$.
Let $\omega':D'\isom D'$ be the automorphism given by the
rule $x'\mapsto x'/(x'+1)$, and set
$F':=\omega^{\prime*}\bar\cK_\chi\La x'\Ra_{|A'}$.
By theorem \ref{th_perversity}(ii) and remark
\ref{rem_perversity}(ii),(iii) we easily see that
$F:=b_*\cF'_\psi(F')$ (notation of \eqref{subsec_split-longish})
is concentrated in degree zero (since $\chi$ is not trivial),
and is locally constant on $A\setminus\{0\}$.
Set $F_\eps:=F_{|\D(0,\eps)}$; a simple calculation using lemma
\ref{lem_Groth-Ogg-Shaf} yields the identities:
$$
\length_\Lambda\cH^0(F_\eps)_0=1
\qquad
\length_\Lambda\cH^0(F_\eps)_{\eta_0}=2
$$
whence $\length_\Lambda R^0\Phi_{\eta_0}(F_\eps)=1$, by virtue
of \eqref{eq_shortish}. Notice that $F'_{\eta_{\infty'}}=\Lambda$
is the trivial $\Lambda[\pi(\infty')]$-module; hence, combining
with \eqref{eq_fake-stationary} and \eqref{eq_triv-equipp} we
conclude that $\cF^{(0',0)}_\psi(F'_{\eta_{0'}})$ and
$\cF^{(-1',0)}_\psi(F'_{\eta_{-1'}})$ must both vanish in this
case, as stated.

The case where $\chi$ is the trivial character is analogous, but
easier : we may take $F'$ to be the extension by zero of the
trivial $\Lambda$-module, with stalk $\Lambda$, on
$(A'\setminus\{0'\})_\et$. Then the same sort of calculation
yield :
$$
\length_\Lambda\cH^0(F_\eps)_0=\length_\Lambda\cH^1(F_\eps)_0=
\length_\Lambda\cH^0(F_\eps)_{\eta_0}=1
$$
and one may then argue as in the foregoing, to conclude.

(ii): The case where $\chi$ is trivial is known by
\eqref{eq_triv-equipp}. The case where $\chi$ is not trivial
follows easily from lemma \ref{lem_tame-FT},
\eqref{eq_fake-stationary}, and assertion (i) : the details
shall be left to the reader.
\end{proof}

\begin{proposition}\label{prop_reveal-fake}
For $z'\in A'(K)$, let $M'$ be any $\Lambda[\pi(z')]$-module
with bounded ramification. Then $\cF^{(z',0)}_\psi(M')=0$.
\end{proposition}
\begin{proof} In view of \eqref{eq_first-step} we may assume
that $z'=0'$. Moreover, we may assume that $M'$ does not admit
any quotient which is a $\Lambda$-module with trivial
$\pi(0')$-action, since the case of a trivial module is
already covered by lemma \ref{lem_fake-vanish}(i).
By proposition \ref{prop_Garuti}, we may
find a Zariski open subset $U'\subset A'\setminus\{0'\}$ and
a locally constant $\Lambda$-module $F'$ on $U'_\et$, with
bounded ramification at every point of $\Sigma':=D'\setminus U'$,
such that $F'$ is tamely ramified at every point of
$\Sigma'\setminus\{0'\}$, and such that $F'_{\eta_{0'}}=M'$.
Set $G:=\cH^0(\cF'_\psi(j'_!F'))$, where $j':U'\to A'$
is the open immersion. The proof of claim \ref{cl_Fourier-type}
shows that
\set\begin{equation}\label{eq_e-uno}
F_\eps=b_*G[0]_{|\D(0,\eps)}
\end{equation}
(notation of \eqref{subsec_split-longish}). Likewise,
\eqref{eq_one-punch} implies that :
\set\begin{equation}\label{eq_e-due}
\length_\Lambda\cH^0(F_\eps)_{\eta_0}=
(\sharp\Sigma'-1)\cdot\length_\Lambda M'+\ssw(M')
\end{equation}
whereas \eqref{eq_will-need-later} implies that
\set\begin{equation}\label{eq_e-tre}
\length_\Lambda\cH^0(F_\eps)_0=
(\sharp\Sigma'-2)\cdot\length_\Lambda M'+\ssw(M').
\end{equation}
Combining \eqref{eq_e-uno}, \eqref{eq_e-due}, \eqref{eq_e-tre}
with \eqref{eq_shortish} and lemma \ref{lem_fake-vanish}(ii),
we deduce that :
$$
\length_\Lambda R^0\Phi_{\eta_0}(F_\eps)=\length_\Lambda M'=
\length_\Lambda\cF^{(\infty',0)}(F'_{\eta_{\infty'}}).
$$
In view of \eqref{eq_fake-stationary}, the contention follows.
\end{proof}

\sset\subsubsection{}
Let $U'\subset A'$ be a Zariski open subset, $F'$ a locally
constant $\Lambda$-module of finite length on $U'_\et$, and
suppose that $F'$ has bounded ramification at all the points
of $D'\setminus U'$.
In light of proposition \ref{prop_reveal-fake}, the (fake)
decomposition \eqref{eq_fake-stationary} is revealed as a
natural isomorphism :
\set\begin{equation}\label{eq_reveal}
R^0\Phi_{\eta_0}(F_\eps)\isom
\cF^{(\infty',0)}_\psi(F'_{\eta_{\infty'}}).
\end{equation}
Therefore, \eqref{eq_shortish} is a short exact sequence :
\set\begin{equation}\label{eq_shortish-true}
0\to H^1_c(A',j'_!F')\to\cH^0b_*\cF'_\psi(j_!F')_{\eta_0}\to
\cF_\psi^{(\infty',0)}(F'_{\eta_{\infty'}})\to H^2_c(A',j'_!F')\to 0.
\end{equation}
By \eqref{subsec_more-trivialities}  we know that
$R^0\Phi_{\eta_0}(F_\eps)$ is a $\Lambda[\pi(0)]$-module,
and we wish now to show that, for every
$\Lambda[\pi(\infty')]$-module $M'$ with bounded ramification,
also $\cF^{(\infty',0)}_\psi(M')$ is naturally a
$\Lambda[\pi(0)]$-module, in such a way that \eqref{eq_reveal}
is actually an isomorphism in $\Lambda[\pi(0)]\bMod$.
This is achieved as for the functors $\cF^{(z,\infty')}_\psi$,
up to some minor modification (and simplification).

Namely, say that $M'$ is represented by a locally constant
$\Lambda$-module of finite length on the \'etale site of
$\D(\infty',\eps_0)^*$, which we may denote again $M'$.
For any $\eps,\delta\in|K^\times|$, let
$$
j_{\eps,\delta}:\D(0,\delta)\times_S\D(\infty',\eps)^*\to
\D(0,\delta)\times_S\D(\infty',\eps)
$$
be the open immersion. Denote also by $\bar q_{\eps,\delta}$
(resp. $\bar q'_{\eps,\delta}$) the projection of
$\D(0,\delta)\times_S\D(\infty',\eps)$ onto $\D(0,\delta)$
(resp. onto $\D(\infty',\eps)$), and set
$$
q_{\eps,\delta}:=\bar q_{\eps,\delta}\circ j_{\eps,\delta}
\qquad
q'_{\eps,\delta}:=\bar q'_{\eps,\delta}\circ j_{\eps,\delta}.
$$
Choose a fundamental pro-covering $(C_H~|~H\subset\pi(0,\delta))$,
and for any open subgroup $H\subset\pi(0,\delta)$, denote by
$j_H:C_H\times_S\D(\infty',\eps)\to\D(0,\delta)\times_S\D(\infty',\eps)$
the morphism obtained by base change from the covering
$C_H\to\D(0,\delta)^*$. For $\eps<\eps_0$, define
$$
\cG'(\eps,\delta):=q^{\prime*}_{\eps,\delta}M'\otimes_\Lambda
\cL_\psi\La m_{\eps,\delta}\Ra
\qquad
\cG'_!(\eps,\delta):=j_{\eps,\delta!}\cG(\eps,\delta)
\qquad
\cG'_*(\eps,\delta):=Rj_{\eps,\delta*}\cG(\eps,\delta)
$$
where $m_{\eps,\delta}$ is the restriction of $m$ to
$\D(0,\delta)\times_S\D(\infty',\eps)^*$. Finally set
$$
N^i(\eps,\delta):=\colim{H\subset\pi(0,\delta)}
H^i(C_H\times_S\D(\infty',\eps)_\et,j^*_H\cG'_!(\eps,\delta)).
$$
Clearly $N^i(\eps,\delta)$ is a $\Lambda[\pi(0,\delta)]$-module,
and just as in \eqref{eq_represent-LFT}, we have a natural
isomorphism:
$$
\cF_\psi^{(\infty',0)}(M')\isom
\colim{n\in\N}\,\colim{\delta>0}\,N^1(\eps_n,\delta).
$$

\begin{lemma}\label{lem_flesh-out}
For every $\eps\leq\eps_0$ there exists $\delta\in|K^\times|$
such that $N^1(\eps,\delta)$ is a $\Lambda$-module of finite
length, and $N^i(\eps,\delta)=0$ for $i\neq 1$.
\end{lemma}
\begin{proof} This is analogous to the proof of proposition
\ref{prop_significant}, so we just provide an outline, which
the industrious reader is invited to flesh out.

To begin with, the usual d\'evissage allows to assume that $\Lambda$
is a field. Next, arguing as in the proof of claim \ref{cl_we-notice},
we reduce to showing that, for a given $\eps$, there exists
$\delta$ such that $R^i\bar q_{\eps,\delta*}\cG'_!(\eps,\delta)$
is locally constant of finite type on $\D(0,\delta)^*_\et$,
for $i=1$, and vanishes for $i\neq 1$.
However, proposition \ref{prop_perverse-FT} implies that the
natural map $\cG'_!(\eps,\delta)\to\cG'_*(\eps,\delta)$
is an isomorphism, hence it suffices to prove that
$R^iq_{\eps,\delta*}\cG'(\eps,\delta)$ is locally constant
of finite type for $i=1$, and vanishes for $i\neq 1$. Using
Poincar\'e duality as in the proof of lemma \ref{lem_no-name}(ii),
this will follow from the :

\begin{claim} For every $\eps\leq\eps_0$ there exists
$\delta\in|K^\times|$ such that
$R^iq_{\eps,\delta!}\cG'(\eps,\delta)$ vanishes for $i\neq 1$,
and is locally constant of finite type for $i=1$.
\end{claim}
\begin{pfclaim}[] For every maximal point $a\in\D(0,\delta)$,
let $\cG(\eps,\delta,a)$ denote the restriction of $\cG(\eps,\delta)$
to $\{\bar a\}\times_S\D(\infty',\eps)^*$, where $\bar a$ is a
geometric point localized at $a$. Moreover, for every
$\gamma\in|K^\times|$ with $\gamma\leq\delta$, let $w(a,\gamma)$
be the unique point of $\{\bar a\}\times_S D'$
such that $(\{\bar a\}\times_S\E(\infty',\gamma))\cup\{w(a,\gamma)\}$
is the topological closure of $\{\bar a\}\times_S\E(\infty',\gamma)$
in $\{\bar a\}\times_S D'$.

By direct inspection we find that, for every $\eps\leq\eps_0$,
there exists $\delta\in|K^\times|$ such that :
\begin{enumerate}
\alphaenu
\item
the stalk $\cG(\eps,\delta,a)_{w'(a,\gamma)}$ has Swan conductor
equal to the Swan conductor of $M'$, for every maximal point
$a\in\D(0,\delta)$ and every $\gamma\in|K^\times|\cap[\eps/2,\eps]$
\item
the Swan conductor of the $\Lambda[\pi(0)]$-module
$\cG(\eps,\delta,a)_{\eta_{\infty'}}$ is a constant function
on the set of maximal points $a$ of $\D(0,\delta)^*$, and all
the breaks of these modules are strictly greater than $0$ (more
precisely, if $\lambda$ is such a break, then $\lambda^\natural\geq 1$).
\end{enumerate}
Indeed, to ensure (b) it suffices, in view of proposition
\ref{prop_drop-is-finite}(i), to take $\delta$ sufficiently
small, so that $\D(0,\delta)^*\cap d_K(M')=\emptyset$.
With these conditions (a) and (b), one may easily redo the
proof of lemma \ref{lem_no-name}(i), and thus obtain the claim.
\end{pfclaim}
\end{proof}

\sset\subsubsection{}\label{subsec_c}
Having lemma \ref{lem_flesh-out} at our disposal, we may now
repeat the considerations of \eqref{subsec_see-here} and
\eqref{subsec_furth-signify}. The upshot is that we have a
well defined functor
$$
\cF^{(\infty',0)}_\psi:
\Lambda[\pi(\infty')]\bMod\to\Lambda[\pi(0)]\bMod
$$
and \eqref{eq_reveal} is indeed $\pi(0)$-equivariant.
Moreover, we can argue as in the proof of proposition
\ref{prop_local-FT-exact}, to deduce that $\cF^{(\infty',0)}_\psi$
is an exact functor.
Next, denote by $c:A\to A$ the automorphism given by the
rule : $x\mapsto -x$; it induces an automorphism $c_*$ of
$\Lambda[\pi(0)]\bMod$. 

\begin{theorem}\label{th_last-local-FT}
With the notation of \eqref{subsec_c}, we have :
\begin{enumerate}
\item
Let $M'$ be a $\Lambda[\pi(\infty')]$-module with bounded
ramification, $\lambda_1\leq\cdots\leq\lambda_l$ the breaks
of $M'$ (repeated with their multiplicities), and set
$k:=\max\{i\leq l~|~\lambda_i^\natural<1\}$.
Then
$$
\length_\Lambda\cF^{(\infty',0)}_\psi(M')=
k-\sum_{i=1}^k\lambda_i^\natural.
$$
\item
The composition of functors
$$
c_*\cF^{(\infty',0)}_\psi(1)\circ\cF^{(0,\infty')}_\psi:
\Lambda[\pi(0)]\bMod\to\Lambda[\pi(0)]\bMod
$$
is naturally isomorphic to the identity functor. (Here $(1)$
denotes the Tate twist.)
\end{enumerate}
\end{theorem}
\begin{proof} (i): Pick as usual a Zariski open immersion
$j':U'\subset A'$ and a locally constant $\Lambda$-module $F'$ on
$U'_\et$, tamely ramified at the points of $\Sigma':=A'\setminus U'$,
and such that $F'_{\eta_{\infty'}}=M'$. Set $F:=\cF'_\psi(j'_!F')$;
according to \eqref{eq_gen-length-FT} we have :
$$
\length_\Lambda\cH^0(F)_{\eta_0}=
\sharp\Sigma'\cdot l+\sum_{i=1}^l\max(0,\lambda_i-1)
$$
(where $l$ is the generic length of $F'$). On the other hand,
lemma \ref{lem_Groth-Ogg-Shaf} shows that
$$
\length_\Lambda\cH^0(F)_0-\length_\Lambda\cH^1(F)_0=
\sharp\Sigma'\cdot l+\sum_{i=1}^l(\lambda_i-1).
$$
Combining with the exact sequence \eqref{eq_shortish} and
with \eqref{eq_reveal}, we deduce the assertion.

(ii): The existence of such an isomorphism is shown as in
\cite[\S 2.5.6]{Lau2} : for a given $\Lambda[\pi(0)]$-module
$M$ with bounded ramification, find a Zariski open subset
$U\subset A\setminus\{0\}$, and a locally constant $\Lambda$-module
$F$ on $U_\et$, with bounded ramification at all the points of
$\Sigma:=D\setminus U$, tamely ramified at the points of
$\Sigma\setminus\{0\}$ and such that $F_{\eta_0}$ is isomorphic
to $M$ (proposition \ref{prop_Garuti}).
Let $j:U\to A$ be the open immersion, set $F':=\cF_\psi(j_!F)$
and notice that $\cH^0(F')$ is locally constant on
$U':=A'\setminus\{0\}$ (theorem \ref{th_perversity}(ii)).
Let $j':U'\to A'$ be the open immersion, and set
$G:=j'_!j^{\prime*}\cH^0(F')[0]$; by theorem \ref{th_perversity}(i),
the cone $C$ of the natural morphism $G\to F'$ is supported
on $\{0'\}$, {\em i.e.} $C=i_{0'*}i^*_{0'}C$, where
$i_{0'}:S\to A'$ is the closed immersion with image $0'$.

In view of remark \ref{rem_perversity}(i) and
\eqref{eq_FT-is-selfdual}, we obtain a distinguished triangle :
$$
\cF'_\psi(G)\to a_*F(-1)\to\pi^*i^*_{0'}C[1]
$$
where $\pi:A'\to S$ is the structure morphism, and
$a=b^{-1}\circ c:A\isom A''$ is $(-1)$-times the double duality
isomorphism. From corollary \ref{cor_vanish-smooth} and
\eqref{eq_reveal}, there follow natural isomorphisms :
$$
\cF^{(\infty',0)}_\psi(G_{\eta_{\infty'}})\isom
R^0\Phi_{\eta_0}b_*\cF'_\psi(G)_{|\D(0,\eps)}\isom
c_*F(-1)_{\eta_0}=c_*M(-1).
$$
On the other hand, \eqref{eq_stat-phase} and proposition
\ref{prop_tame-FT}(ii) yield a natural decomposition :
$$
G_{\eta_{\infty'}}\isom\bigoplus_{z\in\Sigma\setminus\{\infty\}}
\cF^{(z,\infty')}_\psi(F_{\eta_z}).
$$
Moreover, by proposition \ref{prop_tame-FT}(i) and
\eqref{eq_translate-locFT}, for every $z\in\Sigma\setminus\{0,\infty\}$,
the $\Lambda[\pi(\infty')]$-module
$\cF^{(z,\infty')}_\psi(F_{\eta_z})$ has just one break, and this
break equals $1$. Taking (i) into account, we deduce the isomorphism
$$
\cF^{(\infty',0)}_\psi(G_{\eta_{\infty'}})\isom
\cF^{(\infty',0)}_\psi(\cF^{(0,\infty')}_\psi(M))
$$
as required. However, we still have to worry about the naturality
of the isomorphism, {\em i.e.} we have to check that it does not
depend on the choice of $F$. To this aim, let $\tilde F$ be the
$\Lambda$-module on $A_\et$ obtained as extension by zero of
$F_{|\D(0,\eps)^*}$, where $\eps$ is chosen small enough so that
$F$ restricts to a locally constant $\Lambda$-module on
$\D(0,\eps)^*_\et$. Set $\tilde F':=\cF_\psi(\tilde F)$ and
$\tilde G:=j'_!j^{\prime*}\cH^0(\tilde F')[0]$, and notice that
$\cH^i(\tilde F')=0$ for every $i<0$, hence we have again a natural
morphism $\tilde G\to\tilde F'$, as well as a commutative diagram :
$$
\xymatrix{
R^0\Phi_{\eta_0}b_*\cF'_\psi(\tilde G)_{|\D(0,\eps)}
\ar[r]^-\alpha \ar[d] & 
c_*\tilde F(-1)_{\eta_0} \ddouble \\
R^0\Phi_{\eta_0}b_*\cF'_\psi(G)_{|\D(0,\eps)} \ar[r]^-\sim &
c_*F(-1)_{\eta_0}
}$$
whose right vertical arrow is induced by the natural map
$\tilde F\to F$. From remark \ref{rem_not-much}(ii) we get
a natural map
\set\begin{equation}\label{eq_doonesbury}
\cF_\psi^{(0,\infty')}(F_{\eta_0})\to\tilde G_{\eta_{\infty'}}.
\end{equation}
We may find $\delta\in|K^\times|$ small enough, and a locally
constant $\Lambda$-module $T$ of finite type on
$\D(\infty',\delta)^*_\et$, with a $\pi(\infty',\delta)$-equivariant
isomorphism $T_{\eta_{\infty'}}\isom\cF_\psi^{(0,\infty')}(F_{\eta_0})$,
and after shrinking $\delta$, we may also find a morphism
$\omega:T\to\tilde G_{|\D(\infty',\delta)}$ of $\Lambda$-modules
on $\D(\infty',\delta)^*_\et$, such that $\omega_{\eta_{\infty'}}$
equals \eqref{eq_doonesbury} (since $\Lambda$ is noetherian,
$T_{\eta_{\infty'}}$ is a finitely presented $\Lambda$-module).
The morphism $\omega$ is not unique, but any two such morphisms
agree on $\D(\infty',\delta')^*_\et$, for some $\delta'\leq\delta$.

Let $T_!$ be the $\Lambda$-module on $A'_\et$ obtained as extension
by zero of $T$; the morphism $\omega$ extends to a map $T_!\to\tilde G$,
whence an induced map
\set\begin{equation}\label{eq_gogogo}
R^0\Phi_{\eta_0}b_*\cF'_\psi(T_!)_{|\D(0,\eps)}\to
R^0\Phi_{\eta_0}b_*\cF'_\psi(\tilde G)_{|\D(0,\eps)}.
\end{equation}
Furthermore, from remark \ref{rem_not-much-again} we obtain
a natural map :
\set\begin{equation}\label{eq_youre-there}
\cF^{(\infty',0)}_\psi(G_{\eta_{\infty'}})\to
R^0\Phi_{\eta_0}b_*\cF'_\psi(T_!)_{|\D(0,\eps)}
\end{equation}
and by inspecting the construction, it is easily seen that
$\beta:=\eqref{eq_gogogo}\circ\eqref{eq_youre-there}$ does not
depend on the choice of $\omega$, and lastly, $\beta\circ\alpha$
is the isomorphism obtained in the foregoing, so the proof of the
theorem is complete.
\end{proof}

\subsection{Break decomposition}\label{sec_main-theor}
Let $j:U\to A$ be a Zariski open immersion, such that
$0\in A\setminus U$, and denote by $i:A\setminus U\to A$
the closed immersion. Let also $F$ be a locally constant
$\Lambda$-module on $U_\et$, with bounded ramification at
all the points of $D\setminus U$. We consider the diagram
of $\Lambda[\pi(\infty')]$-modules :
\set\begin{equation}\label{eq_bar-local-FT}
{\diagram
(j_*F)_0 \ar[r] \ar[d] & \cH^0\cF_\psi(j_!F)_{\eta_{\infty'}} \ar[d] \\
\cF_\psi^{(0,\infty')}((j_*F)_0) \ar[r] &
\cF_\psi^{(0,\infty')}(F_{\eta_0})
\enddiagram}
\end{equation}
where $(j_*F)_0$ is regarded as an unramified
$\Lambda[\pi(\infty')]$-module. The bottom horizontal arrow is
deduced from the natural map $(j_*F)_0\to F_{\eta_0}$; the left
vertical arrow is the isomorphism of claim \ref{cl_verbatim};
the right vertical arrow is the projection obtained from the
stationary phase decomposition \eqref{eq_stat-phase}, and the
top horizontal arrow comes from the boundary map :
$$
\cF_\psi(i_*i^*j_*F)[-1]\to\cF_\psi(j_*F)
$$
associated to the short exact sequence
$0\to j_!F\to j_*F\to i_*i^*j_*F\to 0$ (see remark
\ref{rem_perversity}(i)).

\begin{lemma}\label{lem_bar-local-FT}
Diagram \eqref{eq_bar-local-FT} commutes.
\end{lemma}
\begin{proof} The detailed verification shall be left as an
exercise for the reader. The main point is that also the left
vertical arrow is defined in terms of a boundary map, just as
for the top horizontal arrow.
\end{proof}

\sset\subsubsection{}\label{subsec_bar-local-FT}
We denote by $\Lambda[\pi(\infty')]\bMod^\dagger$ the full
subcategory of $\Lambda[\pi(\infty')]\bMod$ whose objects
are all the $\pi(\infty')$-modules $M'$ such that :
\begin{enumerate}
\alphaenu
\item
the $\natural$-breaks of $M'$ are different from $1$ 
\item
$M'$ does not admit unramified quotients.
\end{enumerate}
For every object $M'$ of $\Lambda[\pi(\infty')]\bMod^\dagger$,
we define :
$$
\begin{array}{r@{\: := \:}l}
\cG_<(M') &
\cF{}^{(0,\infty')}_\psi\circ c_*\cF^{(\infty',0)}_\psi(M') \\
\cG_>(M') &
\cF_\psi^{(\infty,\infty')}\circ c_*\cF_\psi^{(\infty',\infty)}(M')
\end{array}
$$
(where $c_*$ is defined as in \eqref{subsec_c}). We will exhibit
a natural isomorphism :
\set\begin{equation}\label{eq_separate-slopes}
M'(-1)\isom\cG_<(M')\oplus\cG_>(M').
\end{equation}
To this aim, pick as usual a Zariski open immersion $j':U'\subset A'$
and a locally constant $\Lambda$-module $F'$ on $U'_\et$, tamely
ramified at the points of $A'\setminus U'$, and such that
$F'_{\eta_{\infty'}}=M'$. Set $F:=a^*\cF'_\psi(j'_!F')$, where
$a:A\to A''$ is $(-1)$-times the double duality isomorphism;
from condition (a) and remark \ref{rem_perversity}(ii) we know that
$F$ restricts to a locally constant $\Lambda$-module on
$U:=A\setminus\{0\}$, and condition (b) implies that $\cH^1F_0=0$
(cp. the proof of claim \ref{cl_vanish-H-2}), hence $\cH^1F=0$, by
theorem \ref{th_perversity}(i). In this situation,
\eqref{eq_shortish-true} reduces to the short exact sequence
\set\begin{equation}\label{eq_combine-below}
0\to H^1_c(A',j'_!F')\to\cH^0F_{\eta_0}\to
c_*\cF_\psi^{(\infty',0)}(M')\to 0
\end{equation}
and the stationary phase decomposition (see remark \ref{rem_dual-phase})
yields a natural isomorphism :
\set\begin{equation}\label{eq_combine-also}
\cH^0F_{\eta_\infty}\isom c_*\cF_\psi^{(\infty',\infty)}(M')\oplus R
\end{equation}
where $R$ is a $\Lambda[\pi(\infty)]$-module whose $\natural$-breaks
are all $\leq 1$, by (the dual of) \eqref{eq_translate-locFT} and
proposition \ref{prop_tame-FT}(i).
Let $j:U\to A$ be the open immersion, set $G:=j_!j^*F[0]$, and
let $C$ denote the cone of the natural map $G\to F$; clearly
$C=i_*i^*C$, where $i:\{0\}\to A$ is the closed immersion. From
\eqref{eq_stat-phase} we get as well, a natural isomorphism :
\set\begin{equation}\label{eq_nat-project}
\cH^0\cF_\psi(G)_{\eta_{\infty'}}\isom
\cF_\psi^{(0,\infty')}(\cH^0F_{\eta_0})\oplus
\cF_\psi^{(\infty,\infty')}(\cH^0F_{\eta_\infty}).
\end{equation}
Combining with \eqref{eq_combine-below}, \eqref{eq_combine-also},
and corollary \ref{cor_infty-infty}, we deduce a short exact sequence :
$$
0\to\cF_\psi^{(0,\infty')}(H^1_c(A',j'_!F'))\to
\cH^0\cF_\psi(G)_{\eta_{\infty'}}\xrightarrow{\beta}
\cG_<(M')\oplus\cG_>(M')\to 0.
$$
On the other hand, in view of \eqref{eq_FT-is-selfdual} and remark
\ref{rem_perversity}(i) we have a short exact sequence :
\set\begin{equation}\label{eq_take-into-acc}
0\to(C_0)_A\to\cH^0\cF_\psi(G)\to j'_!F'(-1)\to 0
\end{equation}
and taking into account lemma \ref{lem_bar-local-FT}, we see that
$\beta$ factors through $j'_!F'(-1)_{\eta_{\infty'}}=M'(-1)$, and
induces the sought isomorphism \eqref{eq_separate-slopes} of
$\Lambda[\pi(\infty')]$-modules.
The naturality of \eqref{eq_separate-slopes} can be shown as
in the proof of theorem \ref{th_last-local-FT}(ii) : the details
shall be left to the reader.

\sset\subsubsection{}\label{subsec_clearly-the-breaks}
Clearly, the breaks of $\cG_<(M')$ are also breaks of $M'$. On
the other hand, theorem \ref{th_analyze-breaks}(iii) implies that
the breaks $\lambda$ of $M'$ such that $\lambda^\natural>1$ are
all found (with their multiplicities) in $\cG_>(M')$.
Notice also the natural isomorphism :
$$
\cF_\psi^{(\infty',0)}\circ\cG_<(M')\isom\cF_\psi^{(\infty',0)}(M')(-1)
$$
obtained from theorem \ref{th_last-local-FT}(ii). Taking into account
theorem \ref{th_analyze-breaks}(iii), we deduce :
$$
\length_\Lambda\cF^{(\infty',0)}_\psi(M')=
\length_\Lambda\cG_<(M')-\ssw(\cG_<(M')).
$$
From this, in view of theorem \ref{th_last-local-FT}(i), we
easily deduce that the breaks $\lambda$ of $M'$ such that
$\lambda^\natural<1$ are all found (with their multiplicities)
in $\cG_<(M')$. In other words, by means of the natural isomorphism
\eqref{eq_separate-slopes}, we have been able to separate
the breaks $\lambda$ such that $\lambda^\natural<1$ from those
such that $\lambda^\natural>1$.

\begin{remark} A more detailed analysis might allow to drop
condition (b) of \eqref{subsec_bar-local-FT} from the above
discussion. However, this would probably entail some messy
verifications; instead, we will resort to a trick : given
any $\Lambda[\pi(\infty')]$-module, we shall twist it by a
tamely ramified character of order high enough; this operation
does not alter the distribution of breaks, but it will remove
all unramified quotients.
\end{remark}

\sset\subsubsection{}\label{subsec_trace-on-bmod}
For every $n\in\Z$, let $\phi_n:A'\to A'$ be the map given by
the rule : $x'\mapsto x^{\prime n}$ (where $x'$ is a global coordinate
on $A$). It induces exact endomorphisms $\phi_n^*$ and $\phi_{n*}$
of the category $\Lambda[\pi(\infty')]\bMod$, in the obvious way.
Moreover, the trace map induces a natural transformation
$$
\phi_{n*}\circ\phi^*_n M'\to M'
\qquad
\text{for every object $M'$ of $\Lambda[\pi(\infty')]\bMod$.}
$$
This trace map is an epimorphism whenever $n\in\Lambda^\times$.
For every
$s\in\Q_+$, write $s=ab^{-1}$, with $a,b\in\Z$ and $(a,b)=1$,
and set
$$
\phi_s^*:=\phi_{b*}\circ\phi_a^*
\qquad
\phi_{s*}:=\phi_{a*}\circ\phi_b^*.
$$
Notice the natural map :
\set\begin{equation}\label{eq_trace-counit}
\phi_{s*}\circ\phi^*_sM'\to\phi_{a*}\circ\phi^*_aM'\to M'
\end{equation}
where the first arrow is given by the counit of the adjoint
pair $(\phi^*_b,\phi_{b*})$, and the second arrow is the trace
map; especially, \eqref{eq_trace-counit} is an epimorphism,
provided $a\in\Lambda^\times$. Furthermore, since
$\phi_{b*}=\phi_{b!}$, it is easily seen that $\phi^*_s$
is left adjoint to $\phi_{s*}$.

Suppose that $M'$ is represented by a $\Lambda$-module on
$\D(\infty',\eps)^*_\et$ (for some $\eps\in|K^\times|$), and
recall that the break decomposition of \eqref{subsec_brk-decomp}
is defined on a connected open subset $U\subset\D(\infty',\eps)^*$
such that $U\cap\D(\infty',\delta)^*\neq\emptyset$ for every
$\delta\in|K^\times|$. Then, for every given integer $n>0$, the
subset $\phi^{-1}_nU\subset\D(\infty',\eps^{1/n})^*$ is also connected,
and still intersects every punctured disc centered at $\infty'$.
Likewise, we may find an open connected subset
$V\subset\D(\infty',\eps^n)^*$ such that $\phi^{-1}_n(V)\subset U$,
and $V$ intersects every punctured disc centered at $\infty'$;
thus, $\phi_{n*}(M^{\prime\natural}(r))$ can be seen as a $\Lambda$-module
on $V_\et$, and with this notation we may state the following :

\begin{lemma}\label{lem_mr-penn} 
For every object $M'$ of $\Lambda[\pi(\infty')]\bMod$,
every $n\in\Z$ with $(p,n)=1$, and every $r\in\Q_+$, we
have natural identifications :
$$
\phi_n^*(M^{\prime\natural}(r))=(\phi^*_nM')^\natural(nr)
\qquad
\phi_{n*}(M^{\prime\natural}(r))=(\phi_{n*}M')^\natural(r/n)
$$
(notation of \eqref{subsec_brk-decomp}).
\end{lemma}
\begin{proof} Let $G$ be the automorphism group of a finite
Galois extension of the residue field of the point $\eta(1)\in A'$
(notation of \cite[\S 2.2.10]{Ram}), and suppose that $H\subset G$
is a normal subgroup such that $G/H$ is cyclic of order $n$.
Let us first remark :

\begin{claim}\label{cl_reduce-to-G}
Let $M$ (resp. $N$) be any $\Lambda[G]$-module (resp. a
$\Lambda[H]$-module), and suppose that $M$ (resp. $N$) admits
a single break. Then the same holds for the $\Lambda[H]$-module
$\mathrm{Res}^G_H(M)$ (resp. for the $\Lambda[G]$-module
$\mathrm{Ind}^G_H(N)$).
\end{claim}
\begin{pfclaim} Since $G$ and $H$ have the same $p$-Sylow subgroup
$G^{(p)}$, and since the break decomposition depends only on the
underlying $\Lambda[G^{(p)}]$-module (see \cite[\S4.1.17]{Ram}),
we see more generally that, for any $\Lambda[G]$-module $M$, the
restriction $\mathrm{Res}^G_H(M)$ has the same number of breaks as $M$.
Next, recall that
$\mathrm{Res}^G_H(\mathrm{Ind}^G_H(N))=\bigoplus_{g\in G/H}{}^g\!N$,
so the assertion for the induced module is reduced to the foregoing
case.
\end{pfclaim}

Now, the lemma is easily reduced to the corresponding assertion
for $\Lambda[G]$-modules, where $G$ is a Galois group as above.
For such modules, we may assume that $M'=M^{\prime\natural}(r)$,
in which case -- in view of claim \ref{cl_reduce-to-G} -- it
suffices to compare the Swan conductors of $M'$ and $\phi^*_nM'$
(resp. $\phi_{n*}M'$). The assertion for $\phi_{*n}M'$ follows
straightforwardly from \cite[Lemma 8.1(ii)]{Hu3}. To deal with
$\phi^*_nM'$, we are then reduced to computing the Swan conductor
of $\phi_{*n}(\phi^*_nM')$. The latter is the same as
$M'\otimes_\Lambda\phi_{*n}(\phi^*_n\Lambda)\simeq
M'\otimes_\Lambda\phi_{*n}(\Lambda)$, whence the contention
(details left to the reader).
\end{proof}

\begin{lemma}\label{lem_twist-and-shout}
Let $M'$ be any object of $\Lambda[\pi(\infty')]\bMod$. There
exists a finite set $\Sigma_{M'}\subset\N$, such that
$1\notin\Sigma_{M'}$, and the following holds :
\begin{enumerate}
\item
If $\chi:\bmu_n\to\Lambda^\times$ is a non-trivial character,
and $M'\otimes_\Lambda\cK_\chi\La x'\Ra_{\eta_{\infty'}}$ admits
an unramified quotient, then the order of $\chi$ lies in
$\Sigma_{M'}$. (Notation of \eqref{subsec_tamely}.)
\item
For a given integer $n\neq 0$, suppose that $M'$ does not
admit unramified quotients, but $\phi_n^*M'$ admits an
unramified quotient. Then an element of\/ $\Sigma_{M'}$ divides $n$.
\item
If $M'$ does not admit unramified quotients, then the same
holds for $\phi_{n*}M'$, for every integer $n\neq 0$.
\end{enumerate}
\end{lemma}
\begin{proof} (i): Since the only simple quotient of $\Lambda$ is
its residue field, we may assume from start that $\Lambda$ is a
field. Suppose that we have found non-trivial characters
$\chi_1,\dots,\chi_r$, of orders respectively $a_1,\dots,a_r$,
and epimorphisms
$t_i:M'\otimes_\Lambda\cK_{\chi_i}\La x'\Ra_{\eta_{\infty'}}\to\Lambda$,
for $i=1,\dots,r$, and suppose that $a_i\neq a_j$ for all
$i<j\leq r$. Now, $t_i$ is the same as the datum of an epimorphism
$\tau_i:M'\to\cK_{\chi^{-1}_i}\La x'\Ra$, and we claim that
$\length_\Lambda\bigcap_{i=1}^j\Ker\,\tau_j\leq\length_\Lambda M'-j$
for every $j\leq r$. This is the same as showing that :
\set\begin{equation}\label{eq_notsubset}
\bigcap_{i=1}^{j-1}\Ker\,\tau_i\not\subset\Ker\,\tau_j
\qquad
\text{for every $j\leq r$.}
\end{equation}
By way of contradiction, let $j$ be the smallest index
for which \eqref{eq_notsubset} fails (so $j\geq 2$); then
it is easily seen that the induced map
$$
\tau:M'\to\bigoplus_{i=1}^{j-1}\cK_{\chi^{-1}_i}\La x'\Ra
$$
is an epimorphism, and $\tau_j$ factors through $\tau$,
therefore we get a non-zero map
$\cK_{\chi_i^{-1}}\La x'\Ra\to\cK_{\chi_j^{-1}}$ for some $i<j$,
or what is the same, a non-zero map
$\Lambda\to\cK_{\chi_i\chi^{-1}_j}$. Under our assumptions,
it is easily seen that no such map can exist, whence the
claim. Hence, we must have $r\leq\length_\Lambda M'$, and
the contention follows easily.

(ii): Again, we may assume that $\Lambda$ is a field, say of
characteristic $l$, and moreover that $\Lambda$ is algebraically
closed. Suppose that we have found an epimorphism
$t:\phi^*_nM'\to\Lambda$, for some integer $n>1$.
We may assume that $(l,n)=1$ (since we can assume that our
candidate $\Sigma_{M'}$ contains $l$). This $t$ is the same
as the datum of a map $\tau:M'\to\phi_{n*}\Lambda$. Under our
assumptions, $\phi_{n*}\Lambda$ is the direct sum of
$\Lambda[\pi(\infty')]$-modules of the type
$\cK_\chi\La x'\Ra_{\eta_{\infty'}}$, where $\chi$ ranges over
the characters $\bmu_n\to\Lambda^\times$. Moreover, $\tau$
cannot factor through the trivial character of $\bmu_n$, since
$M'$ does not have unramified quotients. Hence, after composing with
a projection $\phi_{n*}\Lambda\to\cK_\chi\La x'\Ra_{\eta_{\infty'}}$,
we deduce a non-zero map
$M'\otimes_\Lambda\cK_{\chi^{-1}}\La x'\Ra_{\eta_{\infty'}}\to\Lambda$,
where $\chi$ is non-trivial, and its order divides $n$. Then the
assertion follows easily from (i).

(iii): Since $\phi_{n*}=\phi_{n!}$, a non-zero map $\phi_{n*}M'\to N'$
to an unramified $\Lambda[\pi(\infty')]$-module corresponds by
adjunction, to a non-zero map $M'\to\phi_n^*N'$, whence the
contention.
\end{proof}

The following is the main result of this paper.

\begin{theorem}\label{th_main}
Let $M'$ be any object of $\Lambda[\pi(\infty')]\bMod$.
There exists a unique decomposition
$$
M'=\bigoplus_{r\in\Q_+}M'_r
$$
such that each $M'_r$ is an object of $\Lambda[\pi(\infty')]\bMod$,
whose only $\natural$-break is $r$.
\end{theorem}
\begin{proof} For every $\natural$-break $r$ of $M'$, we will exhibit a
decomposition $M'=M'_{[0,r]}\oplus M'_{]r,+\infty[}$, such that the
$\natural$-breaks of $M'_{[0,r]}$ (resp. of $M'_{]r,+\infty[}$ are
all $\leq r$ (resp. $>r$).
If such submodules exist, it is easily seen that they are unique
(since their lengths are determined by $M'$); then the theorem
will follow straightforwardly. Moreover, notice that the theorem
holds for $M'$ if and only if it holds for
$M'\otimes_\Lambda\cK_\chi\La x'\Ra_{\eta_{\infty'}}$, where
$\chi:\bmu_n\to\Lambda^\times$ is any character; by lemma
\ref{lem_twist-and-shout}(i), we may then assume that $M'$
does not admit unramified quotients (recall that $\Lambda^\times$
contains a subgroup isomorphic to $\mu_{p^\infty}$, hence
we may find $\chi$ of arbitrarily high order).

Now, pick a subset $\Sigma_{M'}\subset\N$ as in lemma
\ref{lem_twist-and-shout}(ii), let $r'\in\Q_+$ be the smallest
$\natural$-break of $M'$ such that $r<r'$, and choose
$s=ab^{-1}\in\Q_+$ such that :
$$
r<s^{-1}<r'
\qquad
a\in\Lambda^\times
\qquad
(p,a)=(p,b)=(a,b)=1
$$
and moreover, $(c,a)=1$ for every $c\in\Sigma_{M'}$.
Notice that $\phi^*_s(M')$ does not have any breaks $\lambda$
with $\lambda^\natural=1$ (by lemma \ref{lem_mr-penn}); moreover,
$\phi^*_s(M')$ does not admit any quotients of pure break zero
(by lemma \ref{lem_twist-and-shout}(ii,iii)).
Especially, $\phi^*_s(M')$ is an object of the subcategory
$\Lambda[\pi(\infty')]\bMod^\dagger$, hence the discussion
of \eqref{subsec_bar-local-FT} and \eqref{subsec_clearly-the-breaks},
yields a decomposition
\set\begin{equation}\label{eq_mr-penn}
\phi^*_s(M')(-1)=\cG_<(\phi^*_s(M'))\oplus\cG_>(\phi^*_s(M'))
\end{equation}
such that all the $\natural$-breaks of the first (resp. the second)
summand are $<1$ (resp. $>1$). Whence, a decomposition
$\phi_{s*}\eqref{eq_mr-penn}(1)$ of $\phi_{s*}\phi^*_s(M')$,
and again by lemma \ref{lem_mr-penn}, we see that the first
(resp. second) summand of the latter have $\natural$- breaks in
the range $[0,r]$ (resp. $[r',+\infty[$).

Lastly, using \eqref{eq_trace-counit} we may project
$\phi_{s*}\eqref{eq_mr-penn}(1)$ onto a decomposition of $M'$
with the same separation of $\natural$-breaks, as required
(details left to the reader).
\end{proof}

We may now complete the study of the local Fourier transforms.

\begin{lemma}\label{lem_unram-infty}
Let $M'$ be any unramified $\Lambda[\pi(\infty')]$-module.
There is a natural isomorphism :
$$
\cF_\psi^{(\infty',0)}(M')\isom M'(-1).
$$
\end{lemma}
\begin{proof} Let $F'$ be the constant $\Lambda$-module on $A'$,
such that $F'_{\eta_{\infty'}}=M'$; in this situation, the exact
sequence \eqref{eq_shortish-true} reduces to the isomorphism :
$$
\cF_\psi^{(\infty',0)}(M')\isom H^2_c(A',F')\isom
H^2_c(A',\Lambda)\otimes_\Lambda M'
$$
whence the lemma.
\end{proof}

\begin{proposition}\label{prop_gogogo}
Let $M'$ be any object of $\Lambda[\pi(\infty')]\bMod$, whose
$\natural$-breaks are all $<1$.
Then $\ssw(\cF_\psi^{(\infty',0)}(M'))=\ssw(M')$.
\end{proposition}
\begin{proof} In view of lemma \ref{lem_unram-infty}, we may
assume that $M'$ does not admit unramified quotients, therefore
$M'$ is an object of $\Lambda[\pi(\infty')]\bMod^\dagger$.
Then, define $F'$, $F$ and $G$ as in \eqref{subsec_bar-local-FT};
the short exact sequence \eqref{eq_combine-below} shows that
$$
\ssw(\cF_\psi^{(\infty',0)}(M'))=\ssw(\cH^0F_{\eta_0})=
\ssw(\cH^0G_{\eta_0})=\ssw(\cF_\psi^{(0,\infty')}(G_{\eta_0}))
$$
where the last identity holds by virtue of theorem
\ref{th_analyze-breaks}(ii). Recall that $G$ is locally constant
on $A\setminus\{0\}$, whence a natural stationary phase
decomposition :
$$
\cH^0\cF_\psi(G)_{\eta_{\infty'}}\isom
\cF_\psi^{(0,\infty')}(G_{\eta_0})\oplus
\cF^{(\infty,\infty')}_\psi(G_{\eta_\infty}).
$$
However, $G_{\eta_\infty}=F_{\eta_\infty}$ admits as well a
stationary phase decomposition, whose factors are of the form
$\cF_\psi^{(z,\infty)}(F'_{\eta_z})$, for various $z\in D'(K)$.
Especially, $\cF_\psi^{(\infty',\infty)}(F'_{\eta_{\infty'}})=0$, by
corollary \ref{cor_infty-infty}; on the other hand, the
contributions with $z\neq 0,\infty'$ (resp. with $z=0$) are
all pure of $\natural$-break equal to $1$ (resp. equal to $0$),
by \eqref{eq_translate-locFT} and proposition \ref{prop_tame-FT}(i).
Hence, applying again corollary \ref{cor_infty-infty},
we conclude that $\cF_\psi^{(\infty,\infty')}(G_{\eta_\infty})=0$,
and therefore,
$\ssw(\cF_\psi^{(\infty',0)}(M'))=
\ssw(\cH^0\cF_\psi(G)_{\eta_{\infty'}})$. Taking into account
$\eqref{eq_take-into-acc}_{\eta_{\infty'}}$, the assertion follows.
\end{proof}

\sset\subsubsection{}
For every subset $I\subset\Q_+$, denote by
$\Lambda[\pi(\infty')]\bMod_I$ the full subcategory of the
category $\Lambda[\pi(\infty')]\bMod$, whose objects are the
modules that have all their $\natural$-breaks in $I$.
The following corollary refines theorem \ref{th_analyze-breaks}(iii).

\begin{corollary}\label{cor_used-in-det}
Let $M$ be any $\Lambda[\pi(0)]$-module with bounded ramification.
We have :
\begin{enumerate}
\item
All the $\natural$-breaks of $\cF_\psi^{(0,\infty')}(M)$ are $<1$.
\item
More precisely, $\cF^{(0,\infty)}_\psi$ is an equivalence
$\Lambda[\pi(0)]\bMod\isom\Lambda[\pi(\infty')]\bMod_{[0,1[}$,
whose quasi-inverse is the restriction of the functor
$c_*\cF_\psi^{(\infty',0)}(1)$.
\item
If $M$ has length $l$, and is pure of\/ $\natural$-break $r$,
then $\cF^{(0,\infty)}_\psi(M)$ has length $(1+r)l$, and is
pure of\/ $\natural$-break $r/(1+r)$.
\end{enumerate}
\end{corollary}
\begin{proof} (i): In view of theorem \ref{th_analyze-breaks}(iii)
and theorem \ref{th_main}, we may write
$\cF_\psi^{(0,\infty')}(M)=N'_{[0,1[}\oplus N'_1$, where all the
$\natural$-breaks of $N'_{[0,1[}$ (resp. of $N'_1$) are $<1$
(resp. are equal to $1$). By theorem \ref{th_analyze-breaks}(ii)
we have $\ssw(M)=\ssw(N'_{[0,1[})+\length_\Lambda N'_1$.
On the other hand, theorem \ref{th_last-local-FT} and proposition
\ref{prop_gogogo} imply that
$$
\ssw(M)=\ssw(\cF_\psi^{(\infty',0)}(N'_{[0,1[}))+
\ssw(\cF_\psi^{(\infty',0)}(N'_1))=
\ssw(\cF_\psi^{(\infty',0)}(N'_{[0,1[}))=\ssw(N'_{[0,1[}).
$$
The contention is an immediate consequence.

(ii): This follows directly from (i), in view of theorem
\ref{th_last-local-FT}(ii) and the discussion of
\eqref{subsec_clearly-the-breaks}.

(iii): We may assume that $M'$ is irreducible, in which case
the same holds for $\cF^{(0,\infty')}_\psi(M)$, in view of (ii).
Then the latter has a unique $\natural$-break (by theorem
\ref{th_main}), and the contention follows easily from
theorem \ref{th_analyze-breaks}(i,ii).
\end{proof}

\begin{theorem}
Let $M$ be any object $\Lambda[\pi(\infty)]\bMod_{]1,+\infty[}$.
We have :
\begin{enumerate}
\item
There exists a natural isomorphism of $\Lambda[\pi(\infty)]$-modules :
$$
M(-1)\isom c_*\cF_\psi^{(\infty',\infty)}\circ\cF_\psi^{(\infty,\infty')}(M)
$$
(where $(-1)$ denotes the Tate twist and $c_*$ is defined as in
\eqref{subsec_c}).
\item
More precisely, $\cF_\psi^{(\infty,\infty')}$ restricts to an
equivalence of categories
$$
\Lambda[\pi(\infty)]\bMod_{]1,+\infty[}\isom
\Lambda[\pi(\infty')]\bMod_{]1,+\infty[}
$$
whose quasi-inverse is the restriction of
$c_*\cF_\psi^{(\infty',\infty)}(1)$.
\item
If $M$ has length $l$ and is pure of\/ $\natural$-break $r$,
then $\cF_\psi^{(\infty,\infty')}$ has length $(r-1)l$ and
is pure of\/ $\natural$-break $r/(r-1)$.
\end{enumerate}
\end{theorem}
\begin{proof}(i): As usual, pick a Zariski open immersion $j:U\to A$,
and a locally constant $\Lambda$-module $F$ on $U_\et$, such that
$F_{\eta_\infty}=M$, and $F$ is tamely ramified at all the points
of $A\setminus U$. Set $F':=\cF_\psi(j_!F)$; according to
remark \ref{rem_perversity}(ii,iii), the complex $F'$ is concentrated
in degree zero, and the $\Lambda$-module $\cH^0F'$ is locally constant
on $A'$. In this situation, the stationary phase decomposition
\eqref{eq_stat-phase} yields the natural isomorphisms
$$
\begin{array}{r@{\:\isom\:}l}
\cH^0F'_{\eta_{\infty'}} & \cF_\psi^{(\infty,\infty')}(M)\oplus R \\
a^*\cH^0\cF'_\psi(F')_{\eta_{\infty}} &
c_*\cF_\psi^{(\infty',\infty)}(\cH^0F'_{\eta_{\infty'}})
\end{array}
$$
where $R$ is a $\Lambda[\pi(\infty')]$-module whose $\natural$-breaks
are all $\leq 1$, and whose length is
$\sharp(A\setminus U)\cdot\length_\Lambda M$, by virtue of
\eqref{eq_translate-locFT} and proposition \ref{prop_tame-FT}(i).
The sought isomorphism follows easily, taking into account corollary
\ref{cor_infty-infty}.

(ii): This is an immediate consequence of (i) and the discussion
in \eqref{subsec_clearly-the-breaks}.

(iii): Keep the notation of the proof of (i); on the one hand,
\eqref{eq_gen-length-FT} tells us that
$$
\length_\Lambda\cH^0F'_{\eta_{\infty'}}=
(\sharp(A\setminus U)-1)\cdot\length_\Lambda M+\ssw(M).
$$
On the other hand, the foregoing implies that
$$
\length_\Lambda\cH^0F'_{\eta_{\infty'}}=
\length_\Lambda\cF_\psi^{(\infty,\infty')}(M)+
\sharp(A\setminus U)\cdot\length_\Lambda M.
$$
Therefore :
\set\begin{equation}\label{eq_nico-bella}
\length_\Lambda\cF_\psi^{(\infty,\infty')}(M)=
\ssw(M)-\length_\Lambda M.
\end{equation}
Dually, this also proves that
$\length_\Lambda\cF_\psi^{(\infty',\infty)}(M')=
\ssw(M')-\length_\Lambda M'$ for every object $M'$ of
$\Lambda[\pi(\infty')]\bMod_{]1,+\infty[}$. Letting
$M':=\cF_\psi^{(\infty,\infty')}(M)$, we deduce
-- in view of (i) -- that :
$$
\length_\Lambda M=
\ssw(\cF_\psi^{(\infty,\infty')}(M))-\ssw(M)+\length_\Lambda M
$$
whence :
\set\begin{equation}\label{eq_aleale}
\ssw(\cF_\psi^{(\infty,\infty')}(M))=\ssw(M).
\end{equation}
The contention follows easily by combining \eqref{eq_nico-bella}
and \eqref{eq_aleale}.
\end{proof}

\subsection{Modules of break zero}\label{sec_break-zero}
The most glaring omission in this work concerns the
$\Lambda[\pi(0)]$-modules that are pure of break zero,
or equivalently, those of Swan conductor equal to zero.
The obvious conjecture is that all such modules
are tamely ramified (see \eqref{subsec_tamely});
however, I do not know how to show this. If one can prove
this conjecture, then one can apply a general tannakian
argument due to Y.Andr\'e, to derive an analogue, for
modules with bounded ramification, of the so-called
``local monodromy theorems'' that one encounters in
various situations : see \cite{And}.

In this final section, we present some evidence supporting
our conjecture; namely, proposition \ref{prop_line-of-attack},
a very partial result, which however might suggest a line of
attack for the general case.

\sset\subsubsection{}\label{subsec_swan-function}
We need some preliminaries from \cite{Ram}. Let $R$ be any
artinian $\Z[1/p]$-algebra, $F$ a locally constant, constructible
$R$-module on $\D(0,1)^*_\et$, say with stalks of length $d$.
To $F$ we have associated in \cite[\S4.2.10]{Ram}, a set of 
continuous piecewise linear maps, the {\em break functions}
$$
0\leq f_1(\rho)\leq f_2(\rho)\leq\cdots\leq f_d(\rho)
$$
defined for $\rho\in\R_{\geq 0}$, and with values in $\R_{\geq 0}$
(actually, in {\em loc.cit.} we assume as well that the
stalks of $F$ are free $R$-modules, and then $d$ is taken
to be the common rank of the stalks; but the whole discussion
can adapted easily to our current setting).
The sum :
$$
\delta_F(\rho):=\sum_{i=1}^df_i(\rho)
$$
is therefore a piecewise linear funtion $\R_{\geq 0}\to\R_{\geq 0}$,
and we know that the right derivative of $\delta_F$ at a given
point $\rho:=-\log r$ (for $r\in|K^\times|\cap[0,1]$) is the integer
$$
\ssw^\natural(F,r^+)
$$
called the Swan conductor of $F$ {\em at the point $\eta(r)$} (see
\cite[Prop.4.1.15]{Ram}; the point $\eta(r)$ is the one introduced
in \cite[\S2.2.10]{Ram}).

\sset\subsubsection{}\label{subsec_disc-function}
On the other hand, let $f:X\to\D(0,1)^*$ be a Galois, locally
algebraic \'etale covering; to $f$ we attach a {\em normalized
discriminant function}, as follows. For every
$r\in|K^\times|\cap[0,1]$, let $X_r$ denote a connected component
of $f^{-1}(A(r,1))$, where $A(r,1)$ denotes the affinoid annulus
of internal radius $r$, and external radius $1$; the restriction
$f_r:X_r\to\D(0,r)^*$ is a finite \'etale Galois covering, whose
Galois group we denote $G_r$. In such situation, \cite[\S2.3.12]{Ram}
assigns to $f_r$ a discriminant function
$$
\delta_{f_r}:[0,-\log r]\cap|K^\times|\to\R_{\geq 0}.
$$
Let $s\in|K^\times|\cap[0,1]$, with $s>r$; then we may choose
$X_s$ as one of the connected components of $f_r^{-1}(A(s,1))$,
so that $G_s\subset G_r$, and by inspecting the definitions, it
is easily seen that the restriction of $\delta_{f_r}$ to
$[0,-\log s]\cap|K^\times|$ agrees with $\delta_{f_s}\cdot[G_r:G_s]$.
Therefore :
$$
\delta_f(\rho):=
[G_r:1]^{-1}\delta_{f_r}(\rho)=[G_s:1]^{-1}\delta_{f_s}(\rho)
\qquad
\text{for every $\rho\in[0,-\log s]\cap|K^\times|$.}
$$
By \cite[Th.2.3.25]{Ram}, this function extends to a continuous
piecewise linear function $\delta_f:\R_{\geq 0}\to\R_{\geq 0}$,
which is convex and with rational slopes.

\sset\subsubsection{}\label{subsec_over-Y}
Clearly, the considerations of \eqref{subsec_swan-function}
(resp. of \eqref{subsec_disc-function}) can be repeated for
a locally constant $R$-module $F$ on the \'etale site
of an annulus $A(a,b)$ (resp. for a locally algebraic Galois
covering $f:X\to A(a,b)$). In either case, we obtain a piecewise
linear function defined on $[\log 1/b,\log1/a]$, with values
in $\R_{\geq 0}$. Likewise, the whole discussion can be repeated
for modules (resp. coverings) on $A\setminus\{0\}$. We shall
denote by $Y$ any of these adic spaces, and by $I_Y$ the domain
of the functions $\delta_F$ and $\delta_f$.

Especially, to a $\Lambda$-module $F$ on $Y_\et$, we may associate
the locally algebraic Galois covering $f_F:\cIsom_\Lambda(M_Y,F)\to Y$
(where $M$ is any fixed stalk of $F$ : see \eqref{subsec_was-remark}).

I owe the following simple observation to Ofer Gabber.

\begin{lemma}\label{lem_Gabber-observ}
With the notation of \eqref{subsec_over-Y}, the following
estimates hold for every $\rho\in I_Y$ :
$$
\begin{array}{c}
\delta_F(\rho)/d\leq f_d(\rho)\leq\delta_F(\rho) \\
f_d(\rho)\cdot\left(1-\frac{1}{p}\right)\leq
\delta_{f_F}(\rho)\leq f_d(\rho).
\end{array}
$$
\end{lemma}
\begin{proof} The inequalities for $\delta_F(\rho)$ are obvious.
To show the other two inequalities, we may assume that
$Y=A(a,b)$, and we may replace $f_F$ by its restriction
$f:X\to Y$ to a connected component. Now, pick any
$x\in f^{-1}(\eta(r))$, and let $\mathrm{St}_x\subset G:=\Aut(X/Y)$
be the stabilizer of $x$; the group $\mathrm{St}_x$ admits
an upper numbering ramification filtration :
$$
P^{\flat,\gamma_m}\subset P^{\flat,\gamma_{m-1}}\subset\cdots\subset
P^{\flat,\gamma_1}\subset\mathrm{St}_x
$$
by non-zero normal $p$-subgroups, indexed by certain values
$\gamma_1>\dots>\gamma_{m-1}>\gamma_m\in|K^\times|$, called the
$\flat$-breaks of the filtration (see \cite[\S4.1.26]{Ram}).
Let $(f_i~|~i=1,\dots,d)$ be the break functions of $F$; by
inspecting the definitions, it is easily seen that
\set\begin{equation}\label{eq_sloppy-slope}
f_d(-\log r)=-\log\gamma_m.
\end{equation}
Moreover, we may decompose
$$
\C[G]=\Ind^G_{\mathrm{St}_x}(T)\oplus\C[G/P^{\flat,\gamma_m}]
$$
where $T$ is a direct sum of irreducible
$\C[\mathrm{St}_x]$-modules, on each of which $P^{\flat,\gamma_m}$
acts non-trivially; it follows that the break functions
of the $\C$-module $H_1$ on $Y_\et$ arising from
$\Ind^G_{\mathrm{St}_x}(T)$ all take the value $-\log\gamma_m$,
when evaluated for $\rho=-\log r$. On the other hand, the
break functions of the $\C$-module $H_2$ on $Y_\et$ corresponding
to $\C[C/P^{\flat,\gamma_m}]$ all take values strictly less than
$-\log\gamma_m$, when evaluated for $\rho=-\log r$. Furthermore,
from \cite[Lemma 3.3.10]{Ram} we derive :
$$
\delta_f=[G:1]^{-1}\cdot(\delta_{H_1}+\delta_{H_2}).
$$
Lastly, notice that
$\dim_\C\C[G/P^{\flat,\gamma_m}]\leq p^{-1}\cdot[G:1]$, since
$P^{\flat,\gamma_m}$ is a non-trivial $p$-group. Taking
\eqref{eq_sloppy-slope} into account, the sought
inequalities for $\delta_f(\rho)$ follow straightforwardly.
\end{proof}

\begin{proposition}\label{prop_line-of-attack}
Let $F$ be a locally constant and constructible
$\Lambda$-module on the \'etale site of $A\setminus\{0\}$,
such that the $\Lambda[\pi(0)]$-module $F_{\eta_0}$ (resp.
the $\Lambda[\pi(\infty)]$-module $F_{\eta_\infty}$) has
Swan conductor equal to $0$ (resp. is tamely ramified).
Then $F_{\eta_0}$ is tamely ramified.
\end{proposition}
\begin{proof} Let $\phi_n:A\to A$ be as in
\eqref{subsec_trace-on-bmod}; we may choose $n>0$ such
that $(\phi^*F)_{\eta_\infty}$ is an unramified module, and
clearly it suffices to show that $\phi^*_nF$ is a constant
$\Lambda$-module on the \'etale site of $A\setminus\{0\}$.
Thus, we may replace $F$ by $\phi^*_nF$, and assume from
start that $F$ is locally constant on the \'etale site
of $D^*:=D\setminus\{0\}$. Let $f_F:\cIsom_\Lambda(M_{D^*},F)\to D^*$
(where $M:=F_\infty$) be as in \eqref{subsec_over-Y}.
It suffices to show that the restriction of $f_F$ to each
connected component, is an isomorphism onto $D^*$.

To this aim, notice that, since $F$ is unramified around
$\infty$, the assertion holds for the restriction of $f_F$
to some open subset of the type $f_F^{-1}(\D(\infty,\eps))$.
Consequently, the normalized discriminant function $\delta$
of $X$ (which is defined on the whole of $\R$) vanishes
identically on a half-line $(-\infty,a]$. Moreover, recall
that
$$
0=\ssw^\natural_0(F,0^+)=\lim_{r\to 0^+}\ssw^\natural(F,r^+)
$$
(\cite[Cor.4.1.16]{Ram}). Therefore all the break functions
of $F$ are bounded everywhere on $\R$. In view of lemma
\ref{lem_Gabber-observ}, it follows that $\delta$ is
bounded as well. Since the latter is also convex, we deduce
that $\delta$ actually vanishes identically on $\R$.
Now, for given $r\in|K^\times|$, let $X_r$ be a connected
component of $f^{-1}_F(\D(\infty,r))$, and denote by
$\fd^+_r$ the discriminant of $B_r^+:=\cO^+_X(X_r)$ over
the ring $A(r)^+$ (notation as in \cite[\S2.2.7]{Ram},
and notice that $\fd^+_r$ is well defined, in view of
\cite[Prop.2.3.5(i)]{Ram}). Then $\fd^+_r$ is a unit in
$A(r):=A(r)^+\otimes_{K^+}K$, and therefore it is of the
form $c\cdot u$, where $c\in K^\times$ and $u$ is a unit
in $A(r)^+$. However, since $\delta(-\log r)=0$, the constant
$c$ must actually be a unit of $K^+$, {\em i.e.} $\fd^+_r$
is a unit of $A(r)^+$, therefore $B^+_r$ is an \'etale
$A(r)^+$-algebra. By standard arguments, this implies
that, for every $s\in|K^\times|$ with $s<r$, the restriction
of $f_F$ to the preimage of $\D(\infty,s)$ in $X_r$, splits
as a disjoint union of discs, each of which maps
isomorphically onto $\D(\infty,s)$. The contention
is an immediate consequence.
\end{proof}

\begin{remark}\label{rem_case-of-rank-one}
We should also mention that, if $M$ is a $\Lambda[\pi(0)]$-module
of length one and of Swan conductor zero, then we do know that $M$
is tamely ramified. Indeed, in this case we may extend $M$ to a
locally constant $\Lambda$-module of generic length one on some
Zariski open subset $U\subset A$; such $\Lambda$-modules have been
classified completely in \cite{Ram0}. The assertion follows by a
direct inspection of this classification.
\end{remark}

\subsection{Application : the determinant of cohomology}
\label{subsec_application}
In this section we assume that $K$ is the completion of the
algebraic closure of a field $K_0$ which is a finite extension
of $\Q_p$. Also, we shall take $\Lambda:=\bar\F_\ell$, the
algebraic closure of the field with $\ell$ elements, where
$\ell$ is a prime number different from $p$.

\sset\subsubsection{}
Notice first that the $\mu_{^\infty}$-torsor $\cL$ of \eqref{sec_Fourier}
is already defined over $\Q_p$. The associated $\Lambda$-module
$\cL_\psi$ of \eqref{subsec_Lubin-Tate-torsor} is not necessarily
defined on the \'etale site of $(\A^1_{\Q_p})^\ad$, but we wish to
show that there exists an integer $N\in\N$, depending only on $\ell$,
such that $\cL_\psi$ carries a natural action of
$\Gal(\bar\Q_p/\Q_p[\mu_{p^N}])$. Indeed, set
$$
N:=\max(1,v_p(\ell-1))
\qquad
K_N:=K_0[\mu_{p^N}]
\qquad
S_N:=\Spa(K_N,K_N^+)
$$
where $v_p:\Z\to\N\cup\{+\infty\}$ is the $p$-adic valuation.
It is easily seen that there exist natural isomorphisms of
Galois groups
\set\begin{equation}\label{eq_both-natural}
G_p:=\Gal(\Q_p(\mu_{p^\infty})/\Q_p[\mu_{p^N}])\isom\Z_p\isom
G_\ell:=\Gal(\F_\ell[\mu_{p^\infty}]/\F_\ell[\mu_{p^N}]).
\end{equation}
Let $\hat\Z$ denote the profinite completion of $\Z$; under
\eqref{eq_both-natural}, the surjection
\set\begin{equation}\label{eq_profinite}
\hat\Z\isom\Gal(\Lambda/\F_\ell[\mu_{p^N}])\to G_\ell
\end{equation}
is identified with the natural projection $\hat\Z\to\Z_p$,
and the latter admits a natural (continuous) splitting
$\Z_p\to\hat\Z$. This means that there exists a (continuous)
group homomorphism
\set\begin{equation}\label{eq_already-too-many}
G_p\to\Gal(\Lambda/\F_\ell[\mu_{p^N}])
\qquad
\sigma\mapsto\bar\sigma
\end{equation}
whose composition with \eqref{eq_profinite} equals
\eqref{eq_both-natural}. Let us denote by $\bLambda$
the field $\Lambda$, endowed with the $G_p$-action given
by $\omega$. Then the constant $\Lambda$-module $\bLambda_S$
on the \'etale site of $S$ is correspondingly endowed with a
system of isomorphisms
$$
\sigma_\bLambda:
\sigma_S^*\bLambda_S\isom\bLambda_S
\qquad
\text{for all $\sigma\in G_p$}
$$
(where $\sigma_S:S\isom S$ is the isomorphism deduced from the
automorphism $\sigma$ of $K$) that amount to an action of $G_p$
on this sheaf, such that
$$
\sigma_\bLambda(a\lambda)=\bar\sigma(a)\cdot\sigma_\bLambda(\lambda)
$$
for every $\sigma\in G_p$, every \'etale morphism $U\to S$,
every $\lambda\in\sigma_S^*\bLambda_S(U)$, and every $a\in\Lambda$.
In other words, the action of $G_p$ on $\bLambda_S$ is {\em semilinear}.

On the other hand, $\cL$ is endowed with a corresponding action
of $\Gal(\Q_p[\mu_{p^\infty}]/\Q_p)$, {\em i.e.} a system of
isomorphisms
$$
\sigma_\cL:\sigma_{\A^1_K}^*\cL\isom\cL
\qquad
\text{for all $\sigma\in G_p$}
$$
(where $\sigma_{\A^1_K}:(\A^1_K)^\ad\isom(\A^1_K)^\ad$ is induced
by $\sigma$) compatible with the same action on $\mu_{p^\infty}$,
{\em i.e.} such that
$$
\sigma_\cL(\zeta s)=\sigma(\zeta)\cdot\sigma_\cL(s)
$$
for every $\zeta\in\mu_{p^\infty}$ and every local section $s$
of $\sigma_{\A^1_K}^*\cL$. It follows that we may endow
$\cL_\psi$ with the action
$$
\sigma_\cL\overset{\mu_{p^\infty}}{\times}\sigma_\bLambda:
\sigma_{\A^1_K}^*\cL_\psi\isom\cL_\psi
\qquad
\text{for all $\sigma\in G_p$}
$$
which again, shall be semilinear, relative to the action of
$G_p$ on $\bLambda$.

These consideration lead us to make the following :

\begin{definition}\label{def_semilin-act}
Let $X_N$ be any $S_N$-adic space, and set $X:=X_N\times_{S_N}S$.
The group $H_K:=\Gal(\bar\Q_p/K_0[\mu_{p^N}])$ acts on $X$, and
for every $\sigma\in H_K$, let $\sigma_X:X\isom X$ denote the
corresponding automorphism. A {\em $\bLambda$-module on $X_\et$}
is the datum of a $\Lambda$-module $F$ on $X_\et$ together with a
{\em semilinear action\/} of $H_K$, {\em i.e.} a system of
isomorphisms of $\Lambda$-modules
$$
\sigma_F:\sigma^*_XF\isom F
\qquad
\text{for all $\sigma\in H_K$}
$$
such that
$$
\sigma_F(as)=\bar\sigma(a)\cdot\sigma_F(s)
$$
for every $a\in\Lambda$ and every local section $s$ of $\sigma^*_XF$
(here $\bar\sigma$ denotes the image under \eqref{eq_already-too-many}
of the projection of $\sigma$ in $G_p$), and such that
$$
\tau_F\circ\tau_X^*(\sigma_F)=(\tau\circ\sigma)_F
\qquad
\text{for every $\tau,\sigma\in H_K$.}
$$
The category of locally constant $\bLambda$-modules of finite
type on $X_\et$ shall be denoted $\bLambda_X\Mod_\mathrm{loc}$.
\end{definition}

\begin{example}\label{ex_base-chance-to-Lambda}
In the situation of definition \ref{def_semilin-act}, let $F_N$ be
any $\F_\ell$-module on $X_{N,\et}$, and denote by $\pi_X:X\to X_N$
and $\pi_S:X\to S$ the natural projections. Then the sheaf
$\pi_X^*F_N\otimes_{\F_\ell}\pi_S^*\bLambda_S$ carries a natural
structure of $\bLambda$-module (details left to the reader).
\end{example}

\sset\subsubsection{}\label{subsec_Omegas}
Let $X_N$ and $X$ be as in definition \ref{def_semilin-act}.
Since we shall be working with vanishing cycles, we need to
add Galois equivariance to our geometric constructions of
sections \ref{sec_loc-alg-cov} and \ref{sec_vanish}.
We use a method that works under the assumption that $X$
admits a $K_N$-rational point; this is sufficient for our
purposes. Namely, for every locally algebraic covering
$f:Y\to X$, and every $\sigma\in H_K$, denote $Y^\sigma$ the
fibre product in the cartesian diagram :
$$
\xymatrix{ Y^\sigma \ar[r]^-{\sigma_Y} \ar[d]_{f^\sigma} &
           Y \ar[d]^f \\
           X \ar[r]^-{\sigma_X} & X.
}$$
Clearly $f^\sigma$ is a locally algebraic covering as well,
and the rule $f\mapsto f^\sigma$ defines an automorphism of
the category $\bCov^\localg(X)$. We fix now a geometric point
$\xi:S\to X$, whose image is a $K_N$-rational point. From
$\xi$ we deduce a cartesian diagram :
$$
\xymatrix{
Y^\sigma\times_XS \ar[rr]^-{\sigma_Y\times_X\sigma_S} \ar[d]
& & Y\times_XS \ar[d] \\
S \ar[rr]^-{\sigma_S} & & S
}$$
whence a natural bijection :
$$
\omega_{Y,\sigma}:F_\xi(Y^\sigma)\isom F_\xi(Y)
$$
such that
$$
\omega_{Y,\tau}\circ\omega_{Y^\tau,\sigma}=\omega_{Y,\tau\circ\sigma}
\qquad
\text{for every $\tau,\sigma\in H_K$}.
$$
Next, for every $g\in\pi_1(X,\xi)$, let $g_Y:F_\xi(Y)\isom F_\xi(Y)$
be the corresponding bijection; we define an automorphism
$$
\Omega_\sigma:\pi_1^\localg(X,\xi)\isom\pi_1^\localg(X,\xi)
$$
by the rule :
$$
\Omega_\sigma(g)_Y:=
\omega_{Y,\sigma}\circ g_{Y^\sigma}\circ\omega^{-1}_{Y,\sigma}
$$
for every $g\in\pi_1^\localg(X,\xi)$ and every locally algebraic
covering $Y\to X$. Notice that
$$
\Omega^{-1}_\sigma H(Y,\bar y)=H(Y^\sigma,\omega^{-1}_{Y,\sigma}(\bar y))
$$
for every open subgroup $H(Y,\sigma)$ as in \eqref{subsec_fix-a-max};
especially $\Omega_\sigma$ is a continuous map, for every $\sigma\in H_K$,
and the reader may check that the rule $\sigma\mapsto\Omega_\sigma$
defines a group homomorphism :
$$
\Omega:H_K\to\Aut(\pi_1^\localg(X,\xi))
$$
and we define the arithmetic (locally algebraic) fundamental
group of $X_N$ as the semi-direct product of groups arising
from $\Omega$ :
$$
\pi_1^\localg(X_N,\xi):=\pi_1^\localg(X,\xi)\rtimes_\Omega H_K.
$$

\sset\subsubsection{}
Let $\Sigma$ be a set, and $\rho:\pi_1^\localg(X_N,\xi)\to\Aut(\Sigma)$
a representation. We shall say that {\em $\rho$ is continuous\/}
if the same holds for the restriction of $\rho$ to the subgroup
$\pi_1^\localg(X,\xi)$. Likewise we may define continuous
representations of $\pi_1^\localg(X_N,\xi)$ in $\bLambda$-modules :
namely, these are continuous in the foregoing sense, and moreover
they are semilinear for the action of $\pi_1^\localg(X_N,\xi)$
on $\bLambda$ induced via the projection onto $H_K$. The category
of continuous representations of $\pi_1^\localg(X_N,\xi)$ on
$\bLambda$-modules of finite type shall be denoted
$$
\bLambda[\pi_1^\localg(X_N,\xi)]\Mod_\mathrm{f.cont}.
$$

\begin{lemma}\label{lem_arithmo-fund-equiv}
Let $X_N$ and $X$ be as in definition {\em\ref{def_semilin-act}},
and suppose furthermore, that $X$ is connected. We have :
\begin{enumerate}
\item
Let $F$ be any locally constant $\bLambda$-module on $X_\et$.
The $\pi_1^\localg(X,\xi)$-action on $F_\xi$ extends naturally
to a (continuous) semilinear action of $\pi_1^\localg(X_N,\xi)$.
\item
More precisely, the rule $F\mapsto F_\xi$ defines an equivalence
of categories :
$$
\bLambda_X\Mod_\mathrm{loc}\isom
\bLambda[\pi^\localg_1(X_N,\xi)]\Mod_\mathrm{f.cont}.
$$
\end{enumerate}
\end{lemma}
\begin{proof} More generally, let $\Sigma$ be any (discrete) set,
endowed with a continuous left action $\rho$ of $\pi_1^\localg(X,\xi)$.
By proposition \ref{prop_pi-sets}, the set $\Sigma$ corresponds to a
locally algebraic covering $Y\to X$. Suppose now that $Y$ is endowed
with an action of $H_K$ covering the action on $X$; this means
that we have a system of isomorphisms of $X$-adic spaces :
$$
\sigma_Y:Y^\sigma\isom Y
\qquad
\text{for every $\sigma\in H_K$}
$$
compatible as usual with the composition law of $H_K$. Under the
identification $\omega_{Y,\sigma}$, the left action $\rho_\sigma$
of $\pi_1^\localg(X,\xi)$ on $F_\xi(Y^\sigma)$ is none else than
$\rho^{\Omega_\sigma}$, {\em i.e.}
$$
\rho_\sigma(g)=\rho\circ\Omega_\sigma(g)
\qquad
\text{for every $g\in\pi_1^\localg(X,\xi)$.}
$$
On the other hand, $\sigma_Y$ corresponds to an isomorphism
$t_\sigma:F_\xi(Y^\sigma)\isom F_\xi(Y)$ of sets with left
$\pi_1^\localg(X,\xi)$-action; under the identification
$\omega_{Y,\sigma}$, the bijection $t_\sigma$ becomes an
automorphism of the set $\Sigma$, such that
$$
t_\sigma\circ\rho(g)=\rho_\sigma(g)\circ t_\sigma
\qquad
\text{for every $g\in\pi_1^\localg(X,\xi)$.}
$$
Therefore, we obtain a continuous action of $\pi_1^\localg(X_N,\xi)$
on $\Sigma$ by the rule :
$$
(g,\sigma)\mapsto\rho(g)\circ t_\sigma
\qquad
\text{for every $\sigma\in H_K$ and $g\in\pi_1^\localg(X,\xi)$}.
$$
By reversing the construction, it is clear that from such an action
we may recover a locally algebraic covering $Y$ of $X$, together
with a compatible system of isomorphisms $\sigma_Y$ as above.
In view of the equivalence \eqref{eq_belated}, the same argument
applies unchanged to $\bLambda$-modules, whence the lemma.
\end{proof}

We need also to check the effect of a change of base point.
This is the same as in the algebraic geometric case; namely,
we have :

\begin{lemma}\label{lem_change-of-base-pts}
Let $X$ and $X_N$ be as in lemma {\em\ref{lem_arithmo-fund-equiv}}.
Let also $\xi,\xi':S\to X$ be two geometric points whose images
are $K_N$-rational points. Then there exists an isomorphism of
arithmetic fundamental groups :
$$
\pi_1^\localg(X_N,\xi)\isom\pi_1^\localg(X_N,\xi').
$$
\end{lemma}
\begin{proof} In light of claim \ref{cl_iso-of-fibre-fctrs}, we
have an isomorphism $b:F_\xi\isom F_{\xi'}$ of fibre functors,
inducing an isomorphism
$$
\beta:\pi_1^\localg(X,\xi)\isom\pi_1^\localg(X,\xi').
$$
We claim that the rule
$$
(g,\sigma)\mapsto(\beta(g),\sigma)
\qquad
\text{for every $g\in\pi_1^\localg(X,\xi)$ and $\sigma\in H_K$}
$$
defines an isomorphism as sought. Indeed, let
$$
\Aut(\pi_1^\localg(X,\xi'))\xleftarrow{\ \Omega'\ } H_K
\xrightarrow{\ \Omega\ }\Aut(\pi_1^\localg(X,\xi))
$$
denote the homomorphisms as in \eqref{subsec_Omegas}, associated
to $\xi$ and respectively $\xi'$. The assertion comes down to
the following :

\begin{claim} $\Omega'_\sigma\circ\beta=\beta\circ\Omega_\sigma$
for every $\sigma\in H_K$.
\end{claim}
\begin{pfclaim}[] By inspecting the constructions, we get
the following (non-commutative!) diagram of sets :
$$
\xymatrix{F_{\xi'}(Y^\sigma) \ar[rrr]^-{\omega'_{Y,\sigma}}
\ar[ddd]_{\beta(g_{Y^\sigma})}
& & & F_{\xi'}(Y) \ar[ddd]^{\beta(g_Y)} \\ 
& F_\xi(Y^\sigma) \ar[r]^-{\omega_{Y,\sigma}} \ar[lu] \ar[d]_{g_{Y^\sigma}}
& F_\xi(Y) \ar[ru] \ar[d]^{g_Y} \\
& F_\xi(Y^\sigma) \ar[r]_{\omega_{Y,\sigma}} \ar[dl] & F_\xi(Y) \ar[dr] \\
          F_{\xi'}(Y^\sigma) \ar[rrr]^-{\omega'_{Y,\sigma}} & & & F_{\xi'}(Y)}
$$
(where $\omega'_{Y,\sigma}$ is the natural bijection that yields
$\Omega'$) for every $g\in\pi_1^\localg(X,\xi)$, and every locally
algebraic covering $Y\to X$. The diagonal maps are given by the
isomorphism $b$. The right and left subdiagrams commute, and
the assertion will follow, once we know that the upper (hence also
the the lower) subdiagram commutes as well. However, arguing as in
the proof of claim \ref{cl_iso-of-fibre-fctrs} we may find a
quasi-compact connected open subset $U_N\subset X_N$ such that
$U:=U_N\times_{S_N}S$ contains the images of $\xi$ and $\xi'$,
and we may assume that $b$ comes from an isomorphism of fibre
functors $F_{U,\xi}\isom F_{U,\xi'}$ for the categories of locally
algebraic coverings of $U$ (notice that $U$ is connected, since
$U_N$ contains $K_N$-rational points). Hence, we may replace $X_N$
by $U_N$, and therefore assume from start that $X_N$ is quasi-compact,
in which case $\bCov^\localg(X)$ is the category of ind-objects of
$\bCov^\alg(X)$ (remark \ref{rem_loc-alg-covs}(v)), and under
this identification, both $F_\xi$ and $F_{\xi'}$ are natural
extensions of fibre functors on $\bCov^\alg(X)$.

By the general theory of \cite[Exp.V]{SGA1}, we may find a
fundamental pro-object
$$
P_\bullet:=(P_H~|~H\subset\pi_1^\alg(X,\xi))
$$
indexed by the family of open normal subgroups $H$ of
$\pi_1^\alg(X,\xi)$, and an automorphism $b_\bullet:P_\bullet\isom P_\bullet$
inducing $b$. The latter condition means that $F_\xi$ and $F_{\xi'}$
are both isomorphic to the functor
$$
(Y\to X)\mapsto\Hom_X(P_\bullet,Y):=\colim{H}\Hom_X(P_H,Y)
$$
and under this identification, $b$ becomes the automorphism
of $\Hom_X(P_\bullet,Y)$ given by the rule :
$$
(\phi:P_H\to Y)\mapsto\phi\circ b_H
\qquad
\text{for every $H\subset\pi_1^\alg(X,\xi)$.}
$$
On the other hand, both $\omega_{Y,\sigma}$ and $\omega'_{Y,\sigma}$
get identified with the mapping :
$$
\Hom_X(P_\bullet,Y^\sigma)\isom\Hom_X(P_\bullet,Y)
\qquad
(\phi:P_H\to Y^\sigma)\mapsto(\sigma_Y\circ\phi:P_{\Omega_\sigma(H)}\to Y).
$$
Thus, the commutativity of the upper subdiagram comes down to
the assertion that composition on the left commutes with composition
on the right, which is obvious.
\end{pfclaim}
\end{proof}

\begin{remark}
(i)\ \ 
It is also easily seen that any two ``geometric'' isomorphisms
of our arithmetic fundamental groups ({\em i.e.} those arising
from isomorphisms of fibre functors, as in the proof of lemma
\ref{lem_change-of-base-pts}), differ by an inner automorphism.

(ii)\ \
A purist might prefer to define the arithmetic fundamental
group more directly, as the automorphism group of a fibre
functor for the category of locally algebraic coverings of
$X_N$. This is the approach taken in \cite{deJ}. However,
our method allows to prove rapidly the basic properties one
needs, and avoids certain technical issues. For instance,
basically by decree, our fundamental group sits in a short
exact sequence
$$
1\to\pi_1^\localg(X,\xi)\to\pi_1^\localg(X_N,\xi)\to H_K\to 1.
$$
As de Jong points out, the corresponding sequence for his
fundamental groups is right exact, but it is not clear whether
it is also left exact (\cite[Rem.2.15]{deJ}).

(iii)\ \
On the other hand, \cite{deJ} obtains automatically a topology
on the arithmetic fundamental group, whereas we do not try to
define a topology for our groups (only the kernel of the projection
onto $H_K$ is endowed with its usual topology). Alternatively,
one may say that our approach replaces the profinite topology
of $H_K$ by the discrete topology.

(iv)\ \
Notice that the ``geometric'' isomorphisms of lemma
\ref{lem_change-of-base-pts} are compatible with the
projection onto $H_K$, in the obvious fashion (details left
to the reader).
\end{remark}

\sset\subsubsection{}
Let $A$ and $D$ be as \eqref{subsec_intro-Fourier}, and fix as
usual a global coordinate $x$ on $\A^1_{K_0}$. If $z\in D(K_N)$
is any $K_N$-rational point, there is also a version of definition
\ref{def_semilin-act} for $\Lambda[\pi(z)]$-modules : indeed,
for such $z$, and every $\eps>0$, we have
$$
\D(z,\eps)^*=\D_N(z_0,\eps)^*\times_{S_N}S
\quad\text{where}\quad
\D_N(z_0,\eps)^*:=\{t\in(\A^1_{K_N})^\ad~|~0<|x(t)-z_0|_t\leq\eps\}
$$
(and $z_0\in\A^1_{K_N}$ is the projection of $z$). We have then
the group
$$
\pi(z_0,\eps):=\pi_1^\localg(\D_N(z_0,\eps)^*,\xi)
$$
(for some choice of $K_N$-rational geometric point $\xi$ of
$\D(z,\eps)$) and the corresponding category
$\bLambda[\pi(z_0,\eps)]\Mod_\mathrm{f.cont}$ of continuous
semilinear representations of $\pi(z_0,\eps)$ into $\bLambda$-modules
of finite type. For $\eps'<\eps$ we may pick (thanks to lemma
\ref{lem_change-of-base-pts}) a group homomorphism
$\pi(z_0,\eps')\to\pi(z_0,\eps)$, well defined up to inner
automorphisms, and we may define the category
$$
\bLambda[\pi(z_0)]\Mod
$$
as the $2$-colimit of the system of categories of
$\bLambda[\pi(z_0,\eps)]$-modules of finite type (where the
transition functors are induced by the foregoing group
homomorphisms), as well as its subcategory $\bLambda[\pi(z_0)]\bMod$
of $\bLambda[\pi(z_0)]$-modules with bounded ramification. The usual
operations for modules extend to $\bLambda[\pi(z_0)]$-modules of
finite type; especially, we shall denote
$$
\det:\bLambda[\pi(z_0)]\Mod\to\bLambda[\pi(z_0)]\Mod
$$
the {\em determinant} functor, that assigns to any object
$M$ of $\bLambda[\pi(z_0)]\Mod$ its highest exterior power
$\Lambda^d_{\bar\F_{\ell}}M$ (where $d$ is the length of $M$).

\sset\subsubsection{}
We may now add Galois equivariance to our nearby and vanishing
cycle functors. Indeed, let $X_N$ and $X$ be as in definition
\ref{def_semilin-act}. Also, let $f_N:X_N\to\D_N(z_0,\eps)$ be
a morphism of $S_N$-adic spaces, and set $f:=f_N\times_{S_N}S$.
We use the notation of \eqref{subsec_vanish}, and set as well
$X_H:=C_H\times_{\D(z,\eps)}X$ for every open subgroup
$H\subset\pi(z,\eps)$. By inspecting the definitions, we find,
for every such $H$, and every $\sigma\in H_K$, a commutative
diagram
\set\begin{equation}\label{eq_vanish-galois}
{\diagram X_{\sigma^{-1}H\sigma} \ar[r] \ar[d] & X_H \ar[d] \\
           X \ar[r]^-{\sigma_X} & X
\enddiagram}
\end{equation}
(whose vertical arrows are the projections) as well as a natural
isomorphism of $X$-adic spaces :
$$
X^\sigma_H\isom X_{\sigma^{-1}H\sigma}
$$
and under this isomorphism, the top horizontal arrow of
\eqref{eq_vanish-galois} is identified with $\sigma_{X_H}$
(notation of \eqref{subsec_Omegas}). Let now $F$ be any
sheaf on $X_\et$, and set $F_\sigma:=\sigma^*_XF$; we define
the sheaf $\tilde F$ on $X_\et$ as in \eqref{subsec_vanish},
and by repeating the construction with $F$ replaced by
$F_\sigma$, we obtain likewise $\tilde F_\sigma$. Considering
the diagrams \eqref{eq_vanish-galois} for varying $H$, and
taking the colimit over the system of all open subgroups $H$
of $\pi(z_0,\eps)$, we deduce a natural isomorphism :
$$
\tilde F_\sigma\isom\sigma^*_X\tilde F
\qquad
\text{for every $\sigma\in H_K$}.
$$
Especially, if $F$ is endowed with a system of compatible
isomorphisms $(\sigma_F:\sigma^*_XF\isom F~|~\sigma\in H_K)$,
then $\tilde F$ is endowed with the system of compatible
isomorphisms
$$
(\tilde\sigma_F:\sigma^*_X\tilde F\isom\tilde F~|~\sigma\in H_K).
$$
Next, for any $g\in\pi(z_0,\eps)$, we have a commutative
diagram :
\set\begin{equation}\label{eq_terrible-notation}
{\diagram
X_{\sigma^{-1}H\sigma} \ar[rr]^-{\sigma_{X_H}}
\ar[d]_{(\sigma^{-1}g\sigma)_{X_{\sigma^{-1}H\sigma}}} &&
X_H \ar[d]^{g_{X_H}} \\
X_{\sigma^{-1}g^{-1}Hg\sigma} \ar[rr]^-{\sigma_{X_{g^{-1}Hg}}} &&
X_{g^{-1}Hg}
\enddiagram}
\end{equation}
where $g_{X_H}$ is deduced from the morphism $g_H$ as in
\eqref{subsec_vanish}, and likewise for the left vertical
arrow in \eqref{eq_terrible-notation}. Recall that $\tilde F$
carries a natural continuous action
$$
\rho_F:\pi(z,\eps)\to\Aut(F).
$$
By unwinding the definitions, we deduce from \eqref{eq_terrible-notation}
that the isomorphism $\tilde\sigma_F$ is $\pi(z,\eps)$-equivariant,
provided we endow $\sigma^*_X\tilde F$ with the action
$$
(\sigma^*_X\rho)^{\Omega_\sigma}:=(\sigma^*_X\rho)\circ\Omega_\sigma
$$
where $\Omega_\sigma$ is the automorphism of $\pi(z,\eps)$
defined as in \eqref{subsec_Omegas}. Arguing as in the proof
of lemma \ref{lem_arithmo-fund-equiv}, we deduce that the
stalks of $\tilde F$ over the $K_N$-rational points, inherit
a natural continuous action of $\pi(z_0,\eps)$. Lastly, if $F$
is any $\bLambda$-module on $X_\et$, all the discussion extends
to the stalks of $R^i\Psi_{\eta,\eps}F$, for every $i\in\Z$, in
the usual way.

\sset\subsubsection{}\label{subsec_U_s-and-V_s}
Let $U_N\subset\A^1_{K_N}$ be a (Zariski) open subset, and
$z_0\in\Sigma_0:=\A^1_{K_N}\setminus U_N$ any point. Let also
$\pi_A:A\to(\A^1_{K_N]})^\ad$ be the natural projection
(where $A$ is as in \eqref{subsec_intro-Fourier}), and denote by
$\eps(z_0)$ the minimum of $|z-z'|$, where $(z,z')$ ranges
over the pair of distinct elements of $\pi_A^{-1}(z_0)$; for
every non-zero $\eps<\eps(z_0)$, consider the morphism
$$
p_{z_0,\eps}:V(z_0,\eps):=
\bigcup_{z\in\pi_A^{-1}(z_0)}\D(z,\eps)^*\to\D(0,\eps)^*
$$
whose restriction to $\D(z,\eps)^*$ is the translation $x\mapsto x-z$,
for every $z\in\pi_A^{-1}(z_0)$. It is easily seen that $V(z_0,\eps)$
is invariant under the action of $H_K$ on $A$, and $p_{z_0,\eps}$
is $H_K$-equivariant. Let now $F$ be a $\bLambda$-module on the
\'etale site of $U:=U_N\times_{\Spec\,K_N}S$; it follows easily that
$$
G:=p_{z_0,\eps*}(F_{|V(z_0,\eps)})
$$
inherits from $F$ a natural $\bLambda$-module structure. Suppose
furthermore, that $F$ has bounded ramification at every point of
$\pi_A^{-1}(z_0)$; then it follows that the $\bLambda[\pi(0)]$-module
$$
p_{z_0*}(F)_{\eta_0}:=G_{\eta_0}
$$
has also bounded ramification (clearly this module is independent
of the choice of $\eps$).

\sset\subsubsection{}\label{subsec_mutiple-trace}
Keep the situation of \eqref{subsec_U_s-and-V_s}, and let $j:U\to A$
be the open immersion. Suppose moreover, that $F$ is locally
constant and constructible on $U$, and with bounded ramification
at every point of $D\setminus U$; from \eqref{eq_stat-phase} and
\eqref{eq_translate-locFT} we get the natural decomposition
$$
\cH^0(\cF_\psi(j_!F))_{\eta_{\infty'}}\isom
\cF^{(\infty,\infty')}_\psi(F_{\eta_\infty})\oplus
\bigoplus_{z_0\in\Sigma_0}
\left(\bigoplus_{z\in\pi^{-1}_A(z_0)}
\cF^{(0,\infty')}_\psi(\theta_{z*}F_{\eta_z})\otimes_\Lambda
\cL_\psi\La zx'\Ra_{\eta_{\infty'}}\right).
$$
Since $\cL_\psi$ is a $\bLambda$-module on $A_\et$, the left-hand
side is naturally a $\bLambda[\pi(\infty')]$-module. Likewise,
it is easily seen that $\cF^{(\infty,\infty')}_\psi(F_{\eta_\infty})$
is an object of $\bLambda[\pi(\infty')]\bMod$, and the same holds
for each of the direct sums in parenthesis. Moreover, the decomposition
is equivariant. Taking determinants on both sides, we deduce -- in view
of theorem \ref{th_analyze-breaks}(i) -- the natural isomorphism of
rank one $\bLambda[\pi(\infty')]$-modules :
$$
\det\cH^0(\cF_\psi(j_!F))_{\eta_{\infty'}}\isom
\det\cF^{(\infty,\infty')}_\psi(F_{\eta_\infty})\otimes\!\!
\bigotimes_{z_0\in\Sigma_0}\!\!
\left(\det\cF^{(0,\infty')}_\psi(p_{z_0*}(F)_{\eta_0})\otimes
\cL_\psi\La\delta^F_{z_0}x'\Ra_{\eta_{\infty'}}\right)
$$
where
$$
\delta^F_{z_0}:=\sum_{z\in\pi_A^{-1}(z_0)}a_z(F)\cdot z
\qquad
\text{for every $z_0\in\Sigma_0$}
$$
(notation of \eqref{subsec_Groth-Ogg-Sha}; notice that $a_z(F)$
depends only on $\pi_A(z)$, hence $\delta^F_{z_0}$ is an integer
multiple of the trace $\sum_{z\in\pi_A^{-1}(z_0)}z$).

\begin{remark}\label{rem_dets-are-tame}
(i)\ \
Notice now that the (unique) $\natural$-break of the rank one
$\bLambda[\pi(\infty')]$-module
$\det\cF^{(0,\infty')}_\psi(p_{z_0*}(F)_{\eta_0})$ is $<1$
for each $z_0\in\Sigma_0$ (corollary \ref{cor_used-in-det}(i)),
hence this $\natural$-break must equal $0$, by
\eqref{eq_break-and-swan}, and therefore all these modules are
tamely ramified (remark \ref{rem_case-of-rank-one}).

(ii)\ \
Let $M$ be any tamely ramified $\bLambda[\pi(\infty')]$-module.
Then $M$ admits a natural extension to a $\bLambda$-module on
the \'etale site of $A'\setminus\{0'\}$. Indeed, let $x'$
be a global coordinate on $A'$ as in \eqref{subsec_intro-Fourier},
and set
$$
U'_N:=\Spec\,K_N[x',1/x']
\qquad
U':=U'_N\times_{K_N}\bar\Q_p.
$$
For every $\eps>0$ we have a commutative diagram with exact rows :
$$
\xymatrix{
1 \ar[r] & \pi_1^\alg(\D(\infty',\eps)^*,\xi) \ar[r] \ar[d]_\beta &
\pi_1^\alg(\D_{K_N}(\infty',\eps)^*,\xi) \ar[r] \ar[d] &
H_K \ar[r] \ddouble & 1 \\
1 \ar[r] & \pi_1(U',\xi) \ar[r] &
\pi_1(U'_N,\xi) \ar[r] & H_K \ar[r] & 1
}$$
where $\xi$ is a chosen $K_N$-rational point of $\D(\infty,\eps)^*$.
The closed immersion of $\{\xi\}$ into the $S_N$-adic space
$\D_{K_N}(\infty,\eps)^*$ induces a splitting for both short exact
sequences, and a $\bLambda$-module on $\D(\infty',\eps)^*_\et$ with
finite geometric monodromy corresponds to a continuous representation
$\rho$ of $\pi_1^\alg(\D_{K_N}(\infty',\eps)^*,\xi)$, semilinear with
respect to the action of $H_K$ on $\bLambda$. Then it is clear that
$\rho$ is determined uniquely by its restrictions to
$\pi_1^\alg(\D(\infty',\eps)^*,\xi)$ and to the image of $H_K$
under the fixed splitting. Our tamely ramified $M$ is represented
by a $\bLambda$-module $G$ on $\D(\infty',\eps)^*_\et$ for some
$\eps>0$, and if $\eps$ is small enough, the geometric monodromy
of $G$ factors through the image of the surjective map $\beta$.
Hence $G$ yields naturally both a representation of $\pi_1(U',\xi)$
and a semilinear representation of $H_K$, that assemble into a unique
continuous semilinear representation of $\pi_1(U'_N,\xi)$. The latter
is the sought canonical extension of $M$.

(iii)\ \
Let $M$ be as in (ii), and denote by $G$ the canonical extension
of $M$ to a $\bLambda$-module on $(A'\setminus\{0\})_\et$. Then
we may consider the stalk $G_t$ over the $K_N$-rational point $t$
such that $x'(t)=1$. We shall denote this $\bLambda$-module by
$M_{1'}$.

(iv)\ \
Also, in the situation of \eqref{subsec_mutiple-trace} suppose
that, for some $K_N$-rational point $z\in D\setminus U$, the
$\bLambda[\pi(z)]$-module $F_{\eta_z}$ is unramified, so it
extends to a locally constant $\Lambda$-module $G$ on
$U':=U\cup\{z\}$; in this case, it is easily seen that $G$ is
a $\bLambda$-module on $U'$, and the stalk $G_z$ is a
$\bLambda$-module. With a slight abuse of notation, we shall
denote by $F_z$ this $\bLambda$-module. With this notation,
we may state the following :
\end{remark}

\begin{theorem}\label{th_epsilon-facts}
In the situation of \eqref{subsec_mutiple-trace}, suppose that
$F_{\eta_\infty}$ is an unramified $\bLambda[\pi(\infty)]$-module.
Then there exists a natural isomorphism of rank one $\bLambda$-modules :
$$
\det(R\Gamma_c(U,F)[1])\otimes_{\bLambda}\det(F_\infty(-1))
\isom\bigotimes_{z_0\in\Sigma_0}
\det(\cF^{(0,\infty')}_\psi(p_{z_0*}(F)_{\eta_0}))_1.
$$
\end{theorem}
\begin{proof} We proceed as in the proof of \cite[Th.3.4.2]{Lau2}.
Taking into account \eqref{eq_reveal} (or better, its dual for
$\Lambda$-modules $F$ on $U$) and lemma \ref{lem_unram-infty},
we find a short exact sequence of $\bLambda[\pi(0')]$-modules :
\set\begin{equation}\label{eq_force_unram}
0\to H_c^1(U,F)\to\cH^0\cF_\psi(j_!F)_{\eta_{0'}}\to F_\infty(-1)
\to H^2_c(U,F)\to 0.
\end{equation}
By theorem \ref{th_perversity}(ii), the $\bLambda$-module
$\cH^0\cF_\psi(j_!F)$ is locally constant on $U':=A'\setminus\{0'\}$;
set $F':=\cH^0\cF_\psi(j_!F)_{|U'}$. It follows from
\eqref{eq_force_unram} that $(\det F')_{\eta_{0'}}$ is
unramified, and there is a natural isomorphism of $\bLambda$-modules :
$$
(\det F')_{0'}\isom
\det(R\Gamma_c(U,F)[1])\otimes_{\bLambda}\det(F_\infty(-1))
$$
(notation of remark \ref{rem_dets-are-tame}(iv)). On the other
hand, in view of \eqref{subsec_mutiple-trace} and proposition
\ref{prop_tame-FT}(ii), we have the isomorphism of
$\bLambda[\pi(\infty')]$-modules :
$$
(\det F')_{\eta_{\infty'}}\isom\cL_\psi\La\delta^Fx'\Ra_{\eta_{\infty'}}
\otimes\bigotimes_{z_0\in\Sigma_0}
\det\cF^{(0,\infty')}_\psi(p_{z_0*}(F)_{\eta_0})
$$
where $\delta^F:=\sum_{z_0\in\Sigma_0}\delta^F_{z_0}$. Let
$j':U'\to A'$ be the open immersion; we deduce that
$j'_*(\det F')\otimes_\bLambda\cL_\psi\La-\delta^Fx'\Ra$
is locally constant on $A'$, and tamely ramified at $\infty'$
(remark \ref{rem_dets-are-tame}(i)). Hence the latter is
actually (geometrically) constant on $A'$, and since we have
a natural isomorphism of $\bLambda$-modules :
$$
\cL_\psi\La-\delta^Fx'\Ra\isom\bLambda
$$
the theorem follows.
\end{proof}

\begin{example}\label{ex_algebraic-case}
Let $U_N\subset\A^1_{K_N}$ be a (Zariski) open subset, $F$ a
locally constant $\F_\ell$-module on the \'etale site of the
scheme $U_N$, and set $U:=U_N\times_{K_N}K$. We may then
consider the $\F_\ell$-module $F^\ad$ on $U^\ad_\et$ obtained
by pullback of $F$ along the morphism of sites
$U^\ad_\et\to U_{N,\et}$, and $F^\ad\otimes_{\F_\ell}\bLambda$
is a $\bLambda$-module on $U^\ad_\et$ (example
\ref{ex_base-chance-to-Lambda}).
It is well known that the cohomology of $F^\ad$ on $U^\ad$
is naturally isomorphic to that of $F$ on $U$; hence theorem
\ref{th_epsilon-facts} yields an equivariant decomposition of
$$
\det(R\Gamma_c(U,F))\otimes_{\F_\ell}\det(F_\infty(-1))
\otimes_{\F_\ell}\bLambda
$$
as a tensor product of rank one $\bLambda$-modules attached to
the points of $\A^1_{K_N}\setminus U$, and it is easily seen
that each of the factors in this decomposition is determined
functorially by the local monodromy of $F$ at the corresponding
point.
\end{example}

\begin{remark}\label{rem_unusual}
(i)\ \
The most unusual feature of example \ref{ex_algebraic-case}
is the appearance of semilinear Galois representations in
our decomposition. This might be an artifact of our method,
in which case a different method will eventually be discovered,
yielding a more standard decomposition in terms of linear
characters.

Or else -- and more intriguingly -- it might be a manifestation
of some new hybrid phenomenon, partaking of both the $\ell$-adic
and $p$-adic worlds. Indeed, on the one hand, one encounters
such semilinear representations in the study of $p$-adic Hodge
theory; on the other hand, both Hodge theory and Fourier
transforms find their common roots in harmonic analysis, and the
same can be said of Witten's proof of the Morse inequalities
via stationary phase.

(ii)\ \
According to this interpretation, our ring $\bLambda$
(with its $H_K$-action) should be viewed as a sort of
ring of $\ell$-adic (or better, $\ell$-torsion) periods,
analogous perhaps to the field $\C_p$, the simplest
of rings of $p$-adic periods. It also suggests that one
should redo the whole theory for $p$-adic \'etale coefficients,
and then there should be a comparison functor from $p$-adic
\'etale to deRham $\eps$-factors, in the style of \cite{Bei2}.
However, the technical obstacles involved seem to me too
daunting, at present.

(iii)\ \
If one is not interested in hybridizations of this sort,
one can restrict the Galois actions to the closed subgroup
$\Gal(K/K_0[\mu_{p^\infty}])$, that acts trivially on
$\bLambda$, and hence linearly on all our constructions.
\end{remark}

\begin{remark}\label{rem_end-of-story}
(i)\ \ 
One can generalize example \ref{ex_algebraic-case} to the case
of a locally constant $k$-module on the \'etale site of $U_\et$,
where $k$ is a finite extension of $\F_\ell$. The only difference
is that, in this case, the tensor decomposition shall be equivariant
only for the action of a certain open subgroup of $H_K$ (details
left to the reader).

(ii)\ \
On the other hand, I expect that theorem \ref{th_epsilon-facts}
can be generalized by replacing $\bLambda$ with more general
coefficient rings; especially, with truncated Witt vector
$W_n(\bLambda)$, and even with the completion of
$\Q_\ell(\mu_{p^\infty})$. This generalization seems to be mostly
a technical refinement, and would not bring any new insight, so
I prefer to leave it to a more motivated reader.
\end{remark}

\end{document}